\documentclass[11pt, twoside]{article}

\usepackage[text={390pt, 592pt}, centering, headheight=13.75pt]{geometry}
\usepackage{fancyhdr}
\usepackage{indentfirst}
\usepackage{hyperref}
\usepackage{breakurl} 
\usepackage{array}
\usepackage{amsmath, amssymb}
\usepackage{stmaryrd}
\usepackage{graphicx}
\usepackage[labelsep=none]{caption}

\hypersetup{pdfstartview=FitH, pdfpagemode=UseNone}

\newcounter{lemma}
\newenvironment{lemma}{\refstepcounter{lemma}\emph{Lemma \thelemma.}}{}
\newcounter{theorem}
\newenvironment{theorem}{\refstepcounter{theorem}\emph{Theorem \thetheorem.}}{}
\newenvironment{proof}{\emph{Proof.}}{$\square$}
\newcounter{corollary}
\newenvironment{corollary}{\refstepcounter{corollary}\emph{Corollary \thecorollary.}}{}

\allowdisplaybreaks

\newcommand{\figuresScale}{0.85}
\newcommand{\stillFramesScale}{0.30}

\begin{document}

\vspace*{\bigskipamount}

\noindent{\LARGE \bfseries The Flexibility and Rigidity\\ of Leaper Frameworks}\bigskip

\noindent{\scshape Nikolai Beluhov}\bigskip

\begin{center} \parbox{352pt}{\setlength{\parindent}{12pt} \footnotesize \emph{Abstract}. A leaper framework is a bar-and-joint framework whose joints are integer points forming a rectangular grid and whose bars correspond to all moves of a given leaper within that grid. We study the flexibility and rigidity of leaper frameworks.

Let $p$ and $q$ be positive integers such that the $(p, q)$-leaper $L$ is free. J\'{o}zsef Solymosi and Ethan White conjectured in 2018 that the leaper framework of $L$ on the square grid of side $2(p + q) - 1$, and so on all larger grids, is rigid.

We prove this conjecture. We also prove that Solymosi and White's conjecture is, in a sense, sharp. Namely, the leaper framework of $L$ on the rectangular grid of sides $2(p + q) - 2$ and $2(p + q) - 1$, and so on all smaller grids (except for, trivially, the $1 \times 1$ grid), is flexible. In particular, we completely resolve the flexibility and rigidity question for leaper frameworks on square grids. We establish a number of related results as well.} \end{center}

\section{Introduction} \label{intro}

A framework is a structure built out of bars and joints. All bars that come together at a joint may revolve freely and independently about it, to the extent that the rest of the structure does not limit their movement. Other than that, bars are perfectly rigid; they do not stretch and they do not bend. We think of the joints as points and of the bars as straight-line segments between these points.

(For this introduction, we are only going to give brief, intuitive explanations of concepts. Precise technical definitions are forthcoming in the next section.)

Let $\mathcal{F}$ be some framework. Suppose that it is possible to deform $\mathcal{F}$ continuously so that all bars in it remain straight and retain their length, and yet the framework as a whole changes its shape. Then we say that $\mathcal{F}$ is flexible. Otherwise, if it is not flexible, we say that $\mathcal{F}$ is rigid. For example, every framework in the shape of a triangle is rigid and every framework in the shape of a square is flexible.

The next concept we introduce comes from a completely different part of mathematics.

A leaper is a fairy chess piece which generalises the knight. The $(p, q)$-leaper moves by leaping $p$ units away along one coordinate axis and $q$ units away along the other. For example, when $p = 1$ and $q = 2$, we obtain the familiar knight from orthodox chess. Traditionally, a handful of other leapers are given special names as well. Among them, the $(1, 4)$-leaper is known as the giraffe and the $(2, 3)$-leaper is known as the zebra.

These two threads come together as follows.

A leaper framework is a framework whose joints are integer points which form a rectangular grid and whose bars correspond to all moves of a given leaper within that grid. Thus the underlying graph of a leaper framework is a leaper graph.

\begin{figure}[p] \centering \includegraphics[scale=\figuresScale]{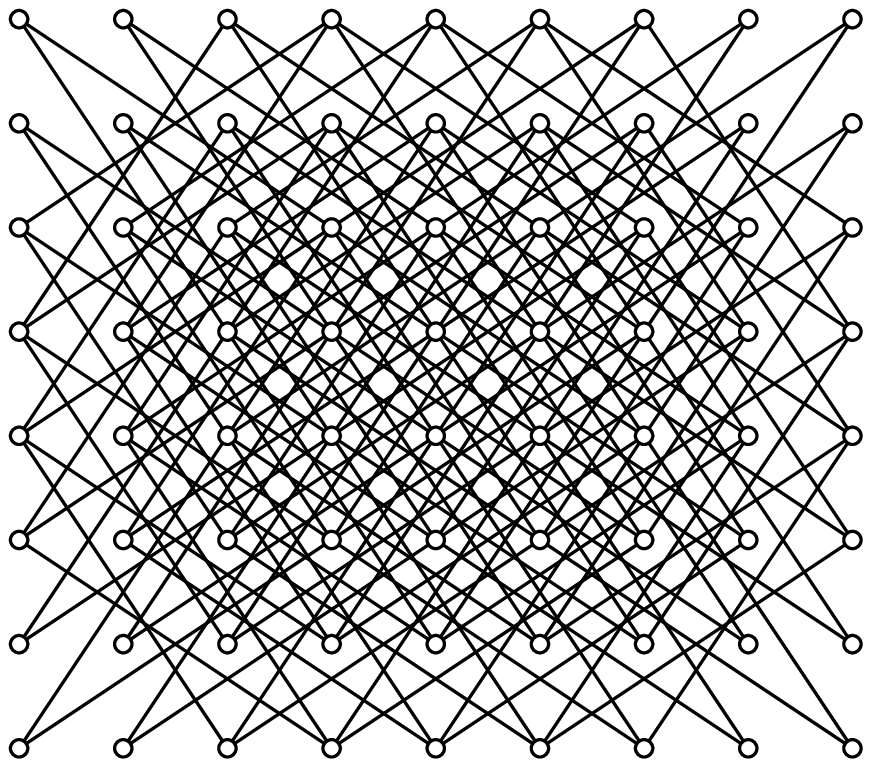} \caption{} \label{zebra-psi} \end{figure}

In particular, in a leaper framework, each joint is common to at most eight bars and all bars are of the same length. For example, Figure \ref{zebra-psi} shows the leaper framework of the zebra on the $8 \times 9$ grid.

The flexibility and rigidity of leaper frameworks is a fascinating question which, as
we shall soon see, gives rise to a lot of beautiful mathematics.

Leaper frameworks were introduced by J\'{o}zsef Solymosi and Ethan White in 2018, in their work \cite{SW}. White's subsequent work \cite{W} revisits the same ideas. They arrive at leaper frameworks in the wider context of bipartite unit-bar frameworks.

Before we go on, we need to motivate one important convention.

From the point of view of flexibility and rigidity, frameworks which consist of multiple independent subframeworks are not particularly interesting. They are flexible in a trivial and unexciting manner, by moving the separate subframeworks relative to each other. Thus it makes sense for us to limit our considerations to connected frameworks only.

As far as leaper frameworks are concerned, this means that we should require the underlying leaper graph to be connected. However, not all leapers are capable of producing connected leaper graphs.

A free leaper is one whose leaper graph is connected on the integer lattice. If a leaper is free, then its leaper graph is connected on all sufficiently large grids as well. And, if a leaper is not free, then its leaper graph is not connected on any grid except for, trivially, the grid of size $1 \times 1$. Thus we are only going to work with free leapers.

Let $p$ and $q$ be two positive integers such that the $(p, q)$-leaper $L$ is free. For this, it is necessary and sufficient that $p + q$ is odd and $p$ and $q$ are relatively prime.

In \cite{SW}, Solymosi and White conjecture that the leaper framework of $L$ is rigid on the square grid of side $2(p + q) - 1$, and so on all larger grids as well. They confirm this conjecture, computationally, in a large but finite number of cases. They also give a short, human-friendly proof for the special case of the knight. This special case is central to their construction of bipartite unit-bar frameworks possessing certain desirable properties.

In the present work, we prove Solymosi and White's conjecture. First, in Section \ref{slope}, we establish a number of purely combinatorial results about leaper graphs. Then, in Section \ref{rigid}, building on these results, we establish the conjecture as well.

\begin{figure}[p] \centering \includegraphics[scale=\figuresScale]{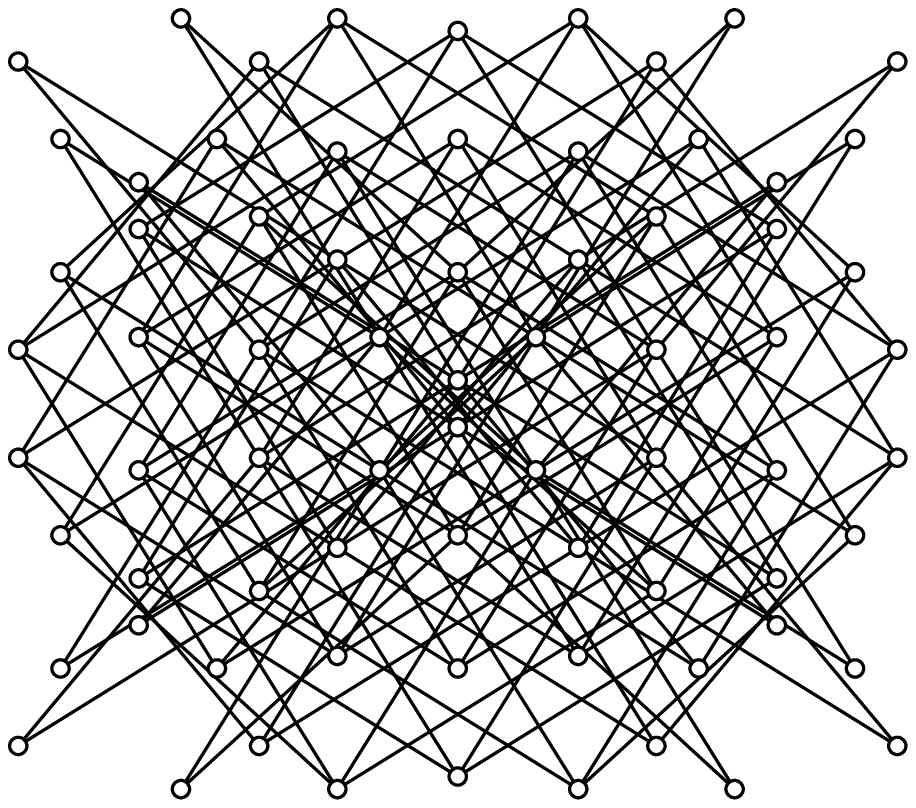} \caption{} \label{zebra-psi-flex} \end{figure}

Delightfully, it turns out that Solymosi and White's conjecture is sharp, in the following sense: Shortening any side of the grid by any amount at all is enough to break rigidity. In other words, the leaper framework of $L$ is flexible on the rectangular grid of sides $2(p + q) - 2$ and $2(p + q) - 1$, and so on all smaller grids as well. (With the trivial exception of the $1 \times 1$ grid.) For example, the zebra framework in Figure \ref{zebra-psi} is flexible, and Figure \ref{zebra-psi-flex} shows it in the process of flexing.

This theorem is somewhat more difficult to prove than the original conjecture. We begin by showing one weaker form of flexibility in Section \ref{inf-flex}. Then, in Section \ref{method}, we develop one reasonably general method for establishing the flexibility of leaper frameworks. The rest of the proof of full flexibility consists of two applications of this method, and spans Sections \ref{flex-i} and \ref{flex-ii}.

These two theorems, on rigidity and on flexibility, are the main results of the present work. Taken together, they completely resolve the question of the flexibility and rigidity of leaper frameworks on square grids. Explicitly, the answer to this question is as follows: Let $n \ge 2$. If $n \le 2(p + q) - 2$, then the leaper framework of $L$ on the $n \times n$ grid is flexible. And, otherwise, if $n \ge 2(p + q) - 1$, then it is rigid.

Finally, in Section \ref{further} we sketch in broad strokes some directions in which we believe the study of the flexibility and rigidity of leaper frameworks could continue.

\section{Preliminaries} \label{prelim}

In this section, we put the vague notions from the introduction on firm formal footing. First we review some basic concepts from the theory of frameworks, then we do the same thing for leapers, and in the end we tie these two threads together.

A lot of what we are going to do generalises immediately to any higher number of dimensions. However, since we are only going to work with two-dimensional leapers and leaper frameworks in the plane (except for, briefly, in Section \ref{further}), we are going to give all definitions for two dimensions only.

We begin with frameworks. For terminology and notation regarding this subject, we mostly follow \cite{W}.

Let $G$ be any finite simple graph.

We denote the edge of $G$ that joins vertices $a$ and $b$ by $ab$. Occasionally, we are going to have to orient some of the edges of otherwise simple graphs. We denote the oriented edge of $G$ that points from vertex $a$ to vertex $b$ by $a \to b$.

Suppose that we have fixed some enumeration $a_1$, $a_2$, \ldots, $a_k$ of the vertices of $G$. Let $\mathbf{c} = (c_1, c_2, \ldots, c_k) \in (\mathbb{R}^2)^k$ be any ordered $k$-tuple of points in the plane. Consider the embedding of $G$ into the plane in which vertex $a_i$ maps onto point $c_i$ for all $i$, and every edge of $G$ maps onto a straight-line segment. In this context, we call $\mathbf{c}$ a \emph{placement} of $G$.

We already gave an informal description of a \emph{framework} in the introduction. Formally, a framework $\mathcal{F}$ is determined by a graph $G$ together with a placement $\mathbf{c}$ of $G$.

We go on to capture the intuitive notion of continuous deformation of a framework.

We say that two placements of $G$ are \emph{equivalent} if every edge of $G$ is represented by segments of the same length in the two placements. Formally, placements $\mathbf{c'}$ and $\mathbf{c''}$ of $G$ are equivalent if $|c'_i - c'_j| = |c''_i - c''_j|$ for all $i$ and $j$ such that $a_i a_j$ is an edge of $G$.

We say that two placements of $G$ are \emph{congruent} if they are congruent figures in the sense of Euclidean geometry, with corresponding elements in the two figures representing the same elements of $G$. Formally, placements $\mathbf{c'}$ and $\mathbf{c''}$ of $G$ are congruent if $|c'_i - c'_j| = |c''_i - c''_j|$ for all $i$ and $j$. This differs from the definition of equivalence only in that we have dropped the restriction to measure distances exclusively along edges of $G$.

Let $F$ be a continuous mapping from the real interval $[0; 1]$ to $(\mathbb{R}^2)^k$. Thus $F$ maps a real parameter, say, $t$, onto a placement $F(t)$ of $G$. Intuitively, $F$ specifies exactly how $\mathcal{F}$ changes its shape as we deform it. Of $F$, we require the following:

(a) $F(0) = \mathbf{c}$. That is, it is in fact $\mathcal{F}$ that we start with.

(b) $F(t)$ is equivalent to $\mathbf{c}$ for all $t$. That is, as we deform $\mathcal{F}$, all of its bars remain straight and retain their length.

(c) $F(t)$ is not congruent to $\mathbf{c}$ for any $t$ except for, as per (b), $t = 0$. That is, our framework does in fact change its shape as we deform it. We are not allowed to simply translate or rotate $\mathcal{F}$.

Then, if $F$ satisfies all three of (a), (b), and (c), we say that $F$ is a \emph{flexion} of $\mathcal{F}$. If $\mathcal{F}$ admits any flexion, then we say that it is \emph{flexible}. Otherwise, if it is not flexible, we say that $\mathcal{F}$ is \emph{rigid}.

These are the notions of flexibility and rigidity which best reflect human intuition. However, other notions of flexibility and rigidity are found in the literature as well. (In this wider context, the notions that we use here are sometimes called \emph{continuous flexibility} and \emph{continuous rigidity} for clarity.) We proceed to an overview of two of them which will be important to us later on.

We say that $\mathcal{F}$ is \emph{globally rigid} if all placements of $G$ equivalent to $\mathbf{c}$ are also congruent to $\mathbf{c}$. Otherwise, if it is not globally rigid, we say that $\mathcal{F}$ is \emph{globally flexible}.

There is no continuity requirement here. Thus global rigidity is a much stronger property than rigidity, and implies it immediately. The converse is false; counterexamples are readily available in the literature, say, in \cite{W}. Intuitively, frameworks which are rigid but not globally rigid can ``snap'' discontinuously between shapes.

Regarding global rigidity, however, there is one important caveat which we should discuss in detail.

Most authors define placements with the additional requirement that the $c_i$ are pairwise distinct. This requirement does not matter much for our purposes, except for in the case of global rigidity.

If we allow some of the $c_i$ to coincide, then of course no bipartite framework in which all bars are of the same length could ever be globally rigid. (And all leaper frameworks possess both of these properties.) To see this, simply place all joints in one part of the underlying bipartite graph at some point $c'$ and all joints in the other part at some other point $c''$ such that the distance between $c'$ and $c''$ equals the bar length. In essence, this placement ``collapses'' the framework.

Thus, for leaper frameworks, we are only going to consider global rigidity with distinctness. For the most part, we are going to say just ``global rigidity'' for simplicity. However, the distinctness convention is important, and we will remind readers of it whenever necessary.

The other alternative notion of flexibility and rigidity that we are going to consider comes from mathematical analysis.

First, a bit of motivation. Suppose that $F$ is some flexion of $\mathcal{F}$, and that the initial velocities of the joints of $\mathcal{F}$ under $F$ are well-defined. These initial velocities are vectors in $\mathbb{R}^2$, and it turns out that they satisfy one remarkably nice system of constraints.

Thus an \emph{infinitesimal motion} of $\mathcal{F}$ is any mapping $I$ from the joints of $\mathcal{F}$ (equivalently, from the vertices of $G$) to $\mathbb{R}^2$ which satisfies that system. Explicitly, the system is as follows: Vectors $c_i - c_j$ and $I(c_i) - I(c_j)$ must be orthogonal for all $i$ and $j$ such that $a_i a_j$ is an edge of $G$.

Some infinitesimal motions of $\mathcal{F}$ are trivial, in the following sense. The continuous motions of $\mathcal{F}$, that is, its translations and rotations, yield solutions to that system just as well as the actual flexions of $\mathcal{F}$ do. Thus we say that $I$ is \emph{nontrivial} if it is not also an infinitesimal motion of the plane. (Strictly speaking, the restriction of an infinitesimal motion of the plane to the joints of $\mathcal{F}$.)

Lastly, we say that $\mathcal{F}$ is \emph{infinitesimally flexible} if there exists some nontrivial infinitesimal motion of $\mathcal{F}$. Otherwise, if it is not infinitesimally flexible, we say that $\mathcal{F}$ is \emph{infinitesimally rigid}.

Flexibility implies infinitesimal flexibility. (This is less obvious than it looks. What if $F$ is not analytic and the initial velocities of the joints of $\mathcal{F}$ under $F$ are not well-defined? We refer interested readers to \cite{W} for details.) Thus also infinitesimal rigidity implies rigidity.

The converses of these statements are false. Moreover, none of global rigidity and infinitesimal rigidity implies the other. Equivalently, none of global flexibility and infinitesimal flexibility implies the other, either. Just as with rigidity and global rigidity, counterexamples are readily available in the literature, say, in \cite{W}.

We may regard every infinitesimal motion $I$ of $\mathcal{F}$ as a vector in $\mathbb{R}^{2k}$. Then all constraints imposed on $I$ become linear and homogeneous, and so all infinitesimal motions of $\mathcal{F}$ form a vector space. We denote this vector space by $\operatorname{Inf}(\mathcal{F})$.

The vector space $\operatorname{Inf}(\mathbb{R}^2)$ of all infinitesimal motions of the plane is a subspace of $\operatorname{Inf}(\mathcal{F})$. (Strictly speaking, its restriction to the joints of $\mathcal{F}$ is.) It is well-known that $\operatorname{Inf}(\mathbb{R}^2)$ is of dimension three; once again, we refer interested readers to \cite{W} for details. Therefore, $\mathcal{F}$ is infinitesimally flexible if and only if the dimension of $\operatorname{Inf}(\mathcal{F})$ is at least four.

This observation suggests the following argument. Every edge of $G$ imposes one homogeneous linear constraint on $I$. The vector space $\mathbb{R}^{2k}$, of which $\operatorname{Inf}(\mathcal{F})$ is a subspace, is of dimension $2k$. Every homogeneous linear constraint that we impose on $I$ reduces the dimension of $\operatorname{Inf}(\mathcal{F})$ by at most one. Therefore, if $G$ has few edges, then the dimension of $\operatorname{Inf}(\mathcal{F})$ will be large.

Thus we arrive at the following lemma.

\medbreak

\begin{lemma} Suppose that $G$ has at most $2k - 4$ edges. Then $\mathcal{F}$ is infinitesimally flexible. \label{linear-bound} \end{lemma}

\medbreak

Both this lemma and its proof are well-known; see \cite{W} for details.

Solymosi and White originally stated their conjecture in terms of infinitesimal rigidity rather than continuous rigidity. We prove this stronger form of the conjecture as well.

This completes our discussion of frameworks in general, and we turn to leapers. To begin with, we introduce the language and notation that we will use to talk about grids.

Formally, a \emph{grid} is the Cartesian product of two intervals of integers.

Thus a grid is a set of ordered pairs of integers. Sometimes, for example in the context of leapers, we are going to consider these pairs as purely combinatorial objects. Then we will refer to them as \emph{cells} to emphasise this. Other times, for example in the context of frameworks, we are going to consider these pairs as integer points in the plane; so, as purely geometric objects. Then we will refer to them as \emph{points} to emphasise that.

Consider the grid $A = I_X \times I_Y = [x + 1; x + w_A] \times [y + 1; y + h_A]$. We refer to $w_A = |I_X|$ as the \emph{width} of $A$ and to $h_A = |I_Y|$ as its \emph{height}. We always give grid sizes in the form ``height times width''. Thus, for example, the size of $A$ becomes $h_A \times w_A$.

Collectively, we refer to $h_A$ and $w_A$ as the \emph{sides} of $A$. A \emph{square} grid is one whose sides are equal. When we want to emphasise the fact that a grid is not necessarily square, or that it is in fact not square, we occasionally say that it is \emph{rectangular}.

We call the number of cells in $A$ its \emph{area}. It equals $h_A w_A$, and, consistently with standard set-theoretic notation, we denote it $|A|$.

Occasionally, we are going to speak of ``the'' grid of some size, as if there were only one such grid. We do this when the property in question depends only on the size of the grid. All of the properties which will be important to us in the present work will be of this kind. In fact, many (though not all, as in Section \ref{slope}) of them will not even depend on which side is taken to be the height of the grid and which one is taken to be its width.

\emph{Column} $i$ of $A$ is the subset $\{x + i\} \times I_Y$ of $A$, and \emph{row} $j$ of $A$ is defined analogously. Occasionally, we use the term \emph{sides} for the outermost rows and columns of $A$ as well.

If the integer intervals $I'_X$ and $I'_Y$ satisfy $I'_X \subseteq I_X$ and $I'_Y \subseteq I_Y$, then we say that the grid $I'_X \times I'_Y$ is a \emph{subgrid} of $A$.

We say that grid $A$ is \emph{smaller} than grid $B$ if $h_A \le h_B$, $w_A \le w_B$, and at least one of the two inequalities is strict. Then we also say that grid $B$ is \emph{larger} than grid $A$. Not all pairs of grids can be compared in this way. When we speak of ``all sufficiently large grids'', what we mean is all grids larger than or equal to some fixed grid.

Occasionally, it will be more convenient to allow rotation in our grid comparisons. We say that grid $A$ \emph{fits inside} grid $B$ if $A$ is congruent, in the sense of Euclidean geometry, to some subgrid of $B$. Equivalently, if at least one of $A$ and the grid of size $w_A \times h_A$ is smaller than or equal to $B$.

Suppose that $h \le h_A$ and $w \le w_A$. Then the \emph{lower left} subgrid of $A$ of size $h \times w$ is the subgrid $[x + 1, x + w] \times [y + 1, y + h]$ of $A$. The \emph{lower right}, \emph{upper left}, and \emph{upper right} subgrids of $A$ of a given size are defined analogously. Subgrids of this form will be important to us. In cases when $h = h_A$, we refer simply to the ``left'' or ``right'' subgrid of $A$ of size $h \times w$, as there is no possibility for confusion. Same goes for cases when $w = w_A$.

The stage is set and we are ready to introduce leapers.

We already gave an informal description of a \emph{leaper} in the introduction. Formally, a leaper $L$ is determined by two nonnegative integers $p$ and $q$.

An \emph{edge} of $L$ is any pair (ordered or unordered depending on the edge being oriented or unoriented) of cells $a = (x_a, y_a)$ and $b = (x_b, y_b)$ such that $\{|x_a - x_b|, |y_a - y_b|\} = \{p, q\}$. A \emph{move} of $L$ is an oriented edge of $L$. Consistently with the general notation for graphs that we introduced earlier, we denote the move of $L$ which leads from cell $a$ to cell $b$ by $a \to b$.

The \emph{leaper graph} of $L$ on the integer lattice is the graph $\mathcal{L}$ whose vertices are all cells of $\mathbb{Z} \times \mathbb{Z}$ and whose edges are all edges of $L$.

Let $C$ be some set of cells. The leaper graph of $L$ on $C$ is the restriction of $\mathcal{L}$ to $C$, which we denote $\mathcal{L}_C$. In this work, $C$ is usually (though not always, Section \ref{flex-ii} being the exception) going to be a grid. For leaper graphs on grids, in cases when only the size $h \times w$ of the grid matters, we occasionally write $\mathcal{L}_{h \times w}$ instead.

As we said in the introduction, $L$ is \emph{free} if $\mathcal{L}$ is connected. Moreover, the necessary and sufficient condition for this is that $p + q$ is odd and $p$ and $q$ are relatively prime. Interested readers may find a proof, as well as detailed references which include earlier statements of this criterion, in Donald Knuth's fundamental work on leaper graphs \cite{K}. There, Knuth also gives a complete description of all grids $A$ such that the leaper graph of a given free leaper $L$ on $A$ is connected. They are all grids $A$ such that the grid of size $(p + q) \times 2q$ fits inside $A$, together with, trivially, the $1 \times 1$ grid.

We already stated in the introduction that we are only going to work with free leapers, and we explained why. There are two more important conventions that we need to introduce.

One of them concerns the $(0, 1)$-leaper, traditionally known as the wazir. It is free, and it is also the only free leaper for which $p$ and $q$ are not both positive. However, its leaper framework is clearly flexible on all grids except for the three grids of sizes $1 \times 1$, $1 \times 2$, and $2 \times 1$. For this reason, we require that both of $p$ and $q$ are positive integers. (As we did in the introduction, without drawing much attention to it.)

The other one is as follows. The behaviour of $L$ depends only on the unordered pair of $p$ and $q$. Moreover, when $p$ and $q$ are equal we do not obtain a free leaper. Therefore, we may assume without loss of generality that $p < q$. This is an important convention, and many of our constructions are not going to make much sense without keeping it in mind.

In short: In all that follows, $p$ and $q$ will positive integers with $p < q$, and $L$ will be free.

We define the \emph{direction} of the move of $L$ from $a$ to $b$ to be vector $b - a$. There are eight possible directions that the moves of $L$ may point in. (Note that, without our conventions that $p$ and $q$ are nonzero and distinct, this claim no longer holds.) They are $(\pm p, \pm q)$ and $(\pm q, \pm p)$, for all choices of the four $\pm$ signs.

We define the \emph{slope} of a move of $L$ of direction $(d', d'')$ to be $\frac{d''}{d'}$, and the slope of an edge of $L$ to be the common slope of its two orientations. There are four possible slopes that occur among the moves and edges of $L$. They are $\frac{p}{q}$, $-\frac{p}{q}$, $\frac{q}{p}$, and $-\frac{q}{p}$.

We are, at long last, properly equipped to define leaper frameworks.

The leaper framework of leaper $L$ on set of cells $C$ is the framework of graph $\mathcal{L}_C$ whose placement is the identity mapping. In other words, the placement of a leaper framework is given simply by reinterpreting cells as points. We call this the \emph{canonical} placement of $\mathcal{L}_C$.

We denote the leaper framework of $L$ on $C$ by $\mathcal{F}_C$. Similarly to leaper graphs, $C$ is usually, though not always, going to be a grid; and, when only the size $h \times w$ of that grid matters, we occasionally write $\mathcal{F}_{h \times w}$ instead.

Every leaper framework is bipartite and all of its bars are of the same length.

The following lemma shows that the flexibility and rigidity of leaper frameworks are both monotone in the size of the grid, one downwards and the other upwards.

\medbreak

\begin{lemma} Let $A$ be any grid. Suppose that $A$ is larger than the $1 \times 1$ grid and that the leaper framework of $L$ on $A$ is rigid. Then the leaper framework of $L$ is rigid on all larger grids as well. Conversely, suppose that the leaper framework of $L$ on $A$ is flexible. Then the leaper framework of $L$ is flexible on all smaller grids as well, with the trivial exception of the $1 \times 1$ grid. \label{monotonicity} \end{lemma}

\medbreak

In \cite{W}, White states (implicitly) and proves this lemma in the special case of the knight. His argument generalises without trouble, and we give the generalisation here for completeness.

\medbreak

\begin{proof} It suffices to prove that if at least one of $h$ and $w$ is greater than one and $\mathcal{F}_{h \times w}$ is rigid, then so is $\mathcal{F}_{h \times (w + 1)}$.

Suppose, for the sake of contradiction, that the leftmost column of the $h \times (w + 1)$ grid contains some cell of degree at most one in $\mathcal{L}_{h \times (w + 1)}$. Then the leftmost column of the $h \times w$ grid also contains some cell of degree at most one in $\mathcal{L}_{h \times w}$, and so $\mathcal{F}_{h \times w}$ cannot be rigid. We have arrived at a contradiction. Therefore, all cells in the leftmost column of the $h \times (w + 1)$ grid are of degree at least two in $\mathcal{L}_{h \times (w + 1)}$.

It follows that every joint in the leftmost column of $\mathcal{F}_{h \times (w + 1)}$ is joined by at least two pairwise nonparallel bars to joints in the rest of the framework.

On the other hand, though, ``the rest of the framework'' is an isomorphic copy $\mathcal{F}'_{h \times w}$ of $\mathcal{F}_{h \times w}$, which we know is rigid. Thus every joint in the leftmost column of $\mathcal{F}_{h \times (w + 1)}$ is fixed into place by the bars that join it to $\mathcal{F}'_{h \times w}$, and the complete framework $\mathcal{F}_{h \times (w + 1)}$ is rigid as well. \end{proof}

\medbreak

Some concrete grids will be of great importance to us, and so we introduce names for them. We list these grids here solely for reference; we will give their definitions once again as we come across them naturally in the course of our discussion.

The grid $\Gamma$ is the rectangular grid of height $2p + q - 1$ and width $p + 2q - 1$.

The grid $\Phi$ is the square grid of side $2(p + q) - 1$. This is the grid which plays the main role in Solymosi and White's conjecture.

The grid $\Psi$ is the rectangular grid of height $2(p + q) - 2$ and width $2(p + q) - 1$. Thus $\Psi$ is just one row shorter than $\Phi$.

Lastly, the grid $\Theta$ is the rectangular grid of height $p + 2q$ and width $2(p + q)$. It makes a brief but significant appearance in Sections \ref{rigid} and \ref{further}.

Most of the symbols we use are going to denote different things in different sections, or even across different parts of the same section. However, the meanings of the symbols $p$, $q$, $L$, $\mathcal{L}$, $\Gamma$, $\Phi$, $\Psi$, and $\Theta$ will be the same throughout our work.

\section{Forbidden Slopes} \label{slope}

In this section, we study the connectedness of leaper graphs in which some slopes have been forbidden, and all edges of these slopes have been omitted. The content of this section is going to be purely combinatorial.

It is far from obvious how this topic ties in with flexibility and rigidity. However, the connection will become clear soon enough, in Section \ref{rigid}.

Let $S$ be some set of slopes. We write $\mathcal{L} \restriction S$ for the subgraph of $\mathcal{L}$ in which only the edges whose slopes are in $S$ have been retained. Dually, we write $\mathcal{L} \setminus S$ for the subgraph of $\mathcal{L}$ obtained by removing all edges whose slopes are in $S$. So, for example, if $S'$ and $S''$ form a partitioning of the set of all slopes, then $\mathcal{L} \restriction S'$ and $\mathcal{L} \setminus S''$ are two different ways to denote the same graph. We omit the curly braces when $S$ contains just one element. Thus, for example, $\mathcal{L} \restriction \{\frac{p}{q}\}$ and $\mathcal{L} \restriction \frac{p}{q}$ are two different ways to denote the same graph.

As with ordinary leaper graphs, given some set of cells $C$, we write $\mathcal{L}_C \restriction S$ and $\mathcal{L}_C \setminus S$ for the restrictions of $\mathcal{L} \restriction S$ and $\mathcal{L} \setminus S$ to $C$, respectively.

For simplicity, throughout this section we are going to state and prove all of our results in terms of concrete sets of slopes. Then we can reflect and rotate these results as needed so as to obtain true statements about other sets of slopes. For example, suppose that we have managed to prove somehow that $\mathcal{L}_{h \times w} \setminus \frac{p}{q}$ is connected, where $h$ and $w$ are some expressions in terms of $p$ and $q$. (As we do in, say, Lemma \ref{three-slopes-connectedness}.) Rotating this statement through $90^\circ$ then tells us that $\mathcal{L}_{w \times h} \setminus -\frac{q}{p}$ is connected as well. When we apply the results of this section later on, we are going to do these conversions silently, without commenting on them every time.

First we examine leaper graphs in which only one slope has been retained. Say that this slope is $\frac{p}{q}$.

\begin{figure}[t] \centering \includegraphics[scale=\figuresScale]{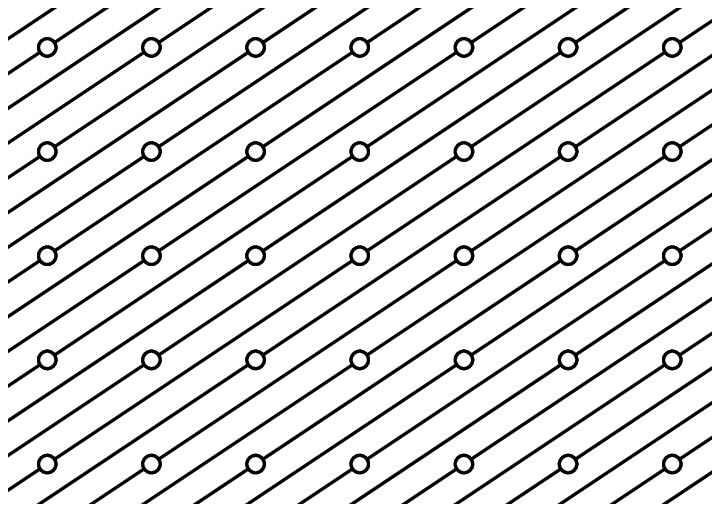} \caption{} \label{zebra-single-slope} \end{figure}

Visually, the graph $\mathcal{L} \restriction \frac{p}{q}$ looks like a pencil of parallel lines spaced evenly at a distance of $\frac{1}{\sqrt{p^2 + q^2}}$ apart. It has infinitely many connected components, and each one of them is a path infinite in both directions. For example, Figure \ref{zebra-single-slope} shows a portion of this graph for the zebra.

Single-slope leaper graphs are not too complicated. Still, one lemma about them will be useful to us.

\medbreak

\begin{lemma} Suppose that $p \le h_A$ and $q \le w_A$. Then the number of connected components of the graph $\mathcal{L}_A \restriction \frac{p}{q}$ is $pw_A + qh_A - pq$. \label{single-slope} \end{lemma}

\medbreak

\begin{proof} Let $A_\texttt{UR}$ be the upper right subgrid of $A$ of size $(h_A - p) \times (w_A - q)$. Then $A \setminus A_\texttt{UR}$ contains precisely one representative out of every connected component of $\mathcal{L}_A \restriction \frac{p}{q}$, and $|A| - |A_\texttt{UR}|$ works out just as it should. \end{proof}

\medbreak

We move on to leaper graphs in which two slopes have been retained and two slopes have been forbidden. There are three geometrically distinct ways to choose the two retained slopes. However, only one of them will be of interest to us here.

Suppose that the two retained slopes are $\frac{p}{q}$ and $-\frac{p}{q}$.

\begin{figure}[t] \centering \includegraphics[scale=\figuresScale]{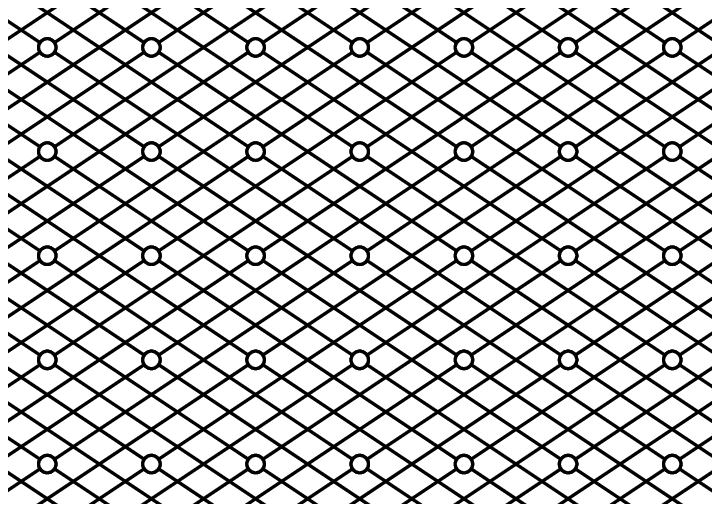} \caption{} \label{zebra-couple-slopes} \end{figure}

\begin{figure}[t] \centering \includegraphics[scale=\figuresScale]{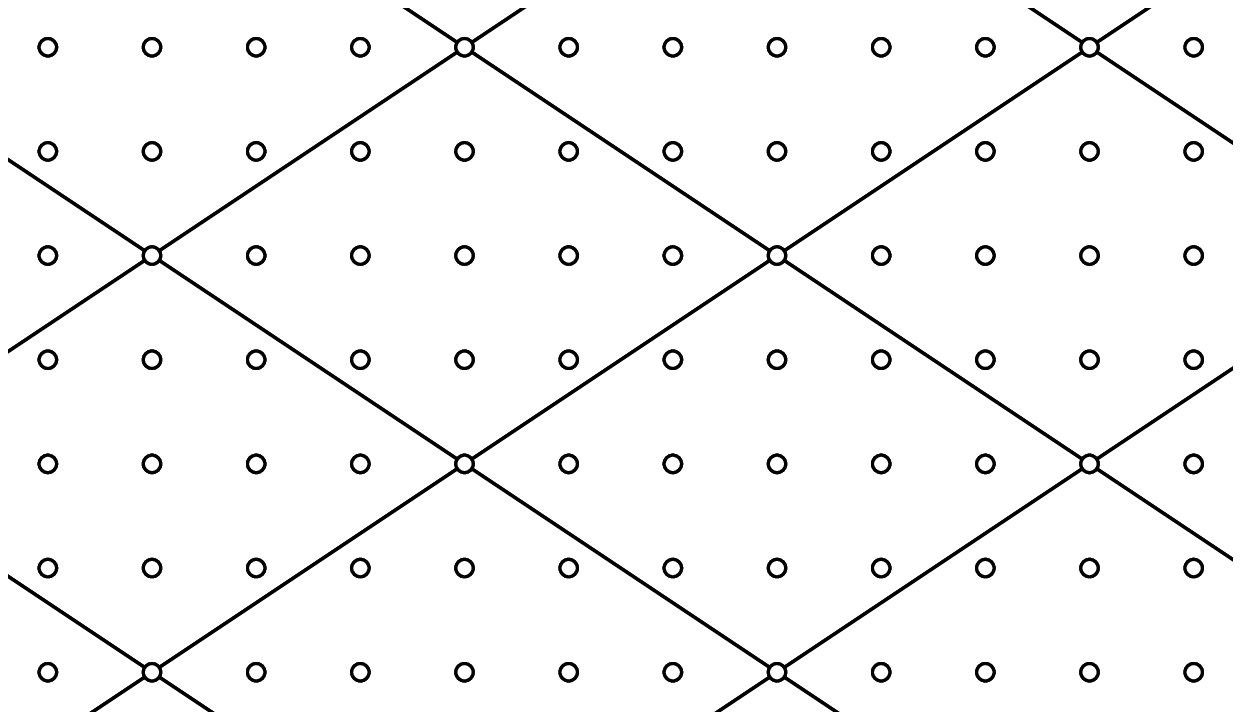} \caption{} \label{zebra-net} \end{figure}

Visually, the graph $\mathcal{L} \restriction \{\frac{p}{q}, -\frac{p}{q}\}$ looks like a dense mesh of little rhombuses, each one of side $\frac{\sqrt{p^2 + q^2}}{2pq}$, formed by two evenly spaced pencils of parallel lines. It has $2pq$ connected components. For example, Figure \ref{zebra-couple-slopes} shows a portion of this graph for the zebra, and Figure \ref{zebra-net} shows a portion of one of its connected components.

The connected components of $\mathcal{L} \restriction \{\frac{p}{q}, -\frac{p}{q}\}$ will be very important to us in this section. We call them \emph{nets}, and we proceed to examine them more closely.

To begin with, every net is isomorphic, in the sense of graph theory, to the unit-distance graph on the integer lattice.

Two cells $a = (x_a, y_a)$ and $b = (x_b, y_b)$ belong to the same net if and only if all three of the following conditions hold:

(a) $x_a \equiv x_b \bmod q$,

(b) $y_a \equiv y_b \bmod p$, and

(c) $\frac{x_a - x_b}{q} \equiv \frac{y_a - y_b}{p} \bmod 2$.

Thus every grid of size $p \times 2q$ intersects each net precisely once, and so does every grid of size $2p \times q$ as well.

We are also going to need the following lemma.

\medbreak

\begin{lemma} Suppose that $2p \le h_A$ and $2q \le w_A$, and let $\mathcal{N}$ be any net. Then the restriction of $\mathcal{N}$ to $A$ is a connected subgraph of $\mathcal{N}$. \label{nets-connectedness} \end{lemma}

\medbreak

\begin{proof} Let $a$ and $b$ be two cells in the restriction of $\mathcal{N}$ to $A$.

Suppose first that $a$ and $b$ are neither in the same row nor in the same column. Then there exists some path in $\mathcal{N}$ from $a$ to $b$ contained entirely within the bounding box of $a$ and $b$. This path will also be in the restriction of $\mathcal{N}$ to $A$.

Suppose, then, that $a$ and $b$ are in the same row. The case when they are in the same column is analogous.

Let $a = (x_a, y_a)$ and $b = (x_b, y_b)$, with $y_a = y_b$. Suppose, without loss of generality, that $x_a < x_b$. Then, since $a$ and $b$ are in the same row and also both in $\mathcal{N}$, we obtain that $x_a \equiv x_b \bmod 2q$.

Since $2p \le h_A$, there are either at least $p$ rows of $A$ below the row of $a$ and $b$, or at least $p$ rows of $A$ above the row of $a$ and $b$. (Or perhaps both.) Without loss of generality, suppose the former. Starting from $a$, let $L$ move successively in directions $(q, -p)$, $(q, p)$, $(q, -p)$, $(q, p)$, \ldots, and so on and so forth. By following this route, $L$ will eventually reach $b$ without ever leaving $\mathcal{N}$ or $A$. \end{proof}

\medbreak

We are ready to add one more slope. That is, we move on to leaper graphs in which only one slope has been forbidden. Let that slope be $\frac{q}{p}$.

\begin{figure}[h] \centering \includegraphics[scale=\figuresScale]{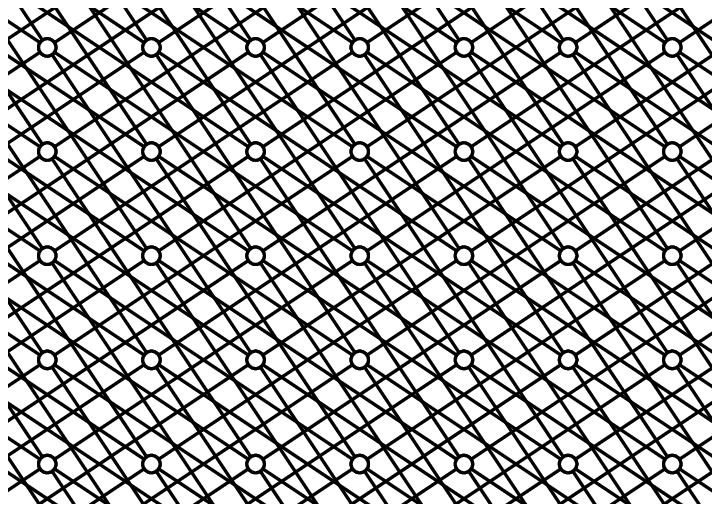} \caption{} \label{zebra-three-slopes} \end{figure}

Visually, the graph $\mathcal{L} \setminus \frac{q}{p}$ does not look significantly less complicated than $\mathcal{L}$. For example, Figure \ref{zebra-three-slopes} shows a portion of the resulting graph for the zebra.

We proceed to examine the structure of $\mathcal{L} \setminus \frac{q}{p}$. The nets of $\mathcal{L} \restriction \{\frac{p}{q}, -\frac{p}{q}\}$ (which is, clearly, a subgraph of $\mathcal{L} \setminus \frac{q}{p}$) will be exceptionally helpful to us in this task.

We say that two nets $\mathcal{N}'$ and $\mathcal{N}''$ are \emph{neighbours} if there is some edge of $L$ in $\mathcal{L} \setminus \frac{q}{p}$ which joins one cell of $\mathcal{N}'$ and one cell of $\mathcal{N}''$. Observe that this edge must necessarily be of slope $-\frac{q}{p}$. Moreover, all moves of $L$ from a given net in a given direction lead to the same net. Therefore, every net has exactly two neighbours, one in direction $(-p, q)$ and the other one in direction $(p, -q)$.

Thus we arrive at one crucial lemma.

\medbreak

\begin{lemma} All nets form a cycle with respect to neighbourhood. \label{nets-cycle} \end{lemma}

\medbreak

\begin{proof} Consider cells $a_i = (-pi, qi)$ for $0 \le i \le 2pq$. They form a path of $L$ of length $2pq$ edges such that all edges in it are of slope $-\frac{q}{p}$. This path starts and ends in the same net. We proceed to show that it also visits every other net precisely once. Clearly, this suffices.

To this end, let $b = (x_b, y_b)$ be any cell. Since $p$ and $q$ are relatively prime, by the Chinese Remainder Theorem there exists some $j$, unique modulo $pq$, such that $-pj \equiv x_b \bmod q$ and $qj \equiv y_b \bmod p$. This settles conditions (a) and (b) for two cells to be in the same net.

Also, for the coordinates of cells $a_j$ and $a_{j + pq}$, we obtain \[\left(\frac{-pj - x_b}{q} - \frac{qj - y_b}{p}\right) - \left(\frac{-p(j + pq) - x_b}{q} - \frac{q(j + pq) - y_b}{p}\right) = p^2 + q^2.\]

Since $p + q$ is odd, so is the right-hand side of the above identity. Therefore, precisely one of cells $a_j$ and $a_{j + pq}$ satisfies condition (c) as well, and so belongs to the net of cell $b$. \end{proof}

\medbreak

From Lemma \ref{nets-cycle}, we derive the following corollary.

\medbreak

\begin{corollary} The graph $\mathcal{L} \setminus \frac{q}{p}$ is connected. \label{three-slopes-lattice} \end{corollary}

\medbreak

Let $\mathcal{N}_1$, $\mathcal{N}_2$, \ldots, $\mathcal{N}_{2pq}$ be some cyclic enumeration of all nets such that two nets are neighbours if and only if they occupy neighbouring positions in it. In other words, we consider indices in the enumeration modulo $2pq$, so that, for example, $\mathcal{N}_1$ and $\mathcal{N}_{2pq + 1}$ are two different ways of denoting the same net; and two nets $\mathcal{N}_i$ and $\mathcal{N}_j$ are neighbours if and only if $j - i \equiv \pm 1 \bmod 2pq$.

Moreover, we may assume without loss of generality that the neighbour of $\mathcal{N}_i$ in direction $(-p, q)$ is $\mathcal{N}_{i + 1}$ for all $i$.

By this point, we have become properly equipped to tackle the following lemma.

\medbreak

\begin{lemma} Let $\Gamma$ be the rectangular grid of height $2p + q - 1$ and width $p + 2q - 1$. Then the graph $\mathcal{L}_\Gamma \setminus \frac{q}{p}$ is connected. \label{three-slopes-connectedness} \end{lemma}

\medbreak

\begin{proof} Let $\mathcal{N}'_i$ be the restriction of net $\mathcal{N}_i$ to $\Gamma$.

We say that two nets $\mathcal{N}'$ and $\mathcal{N}''$ are \emph{neighbours on $\Gamma$} if there is some edge of $L$ in $\mathcal{L}_\Gamma \setminus \frac{q}{p}$ which joins one cell of $\mathcal{N}'$ and one cell of $\mathcal{N}''$. Clearly, if two nets are neighbours on $\Gamma$, then they are also neighbours on the integer lattice. We proceed to examine the ``degree of truth'' of the converse.

We say that two nets form a \emph{stopper pair} if they are neighbours on the integer lattice but not on $\Gamma$. We are going to prove that there is precisely one stopper pair of nets.

Consider, to this end, two neighbouring nets $\mathcal{N}_i$ and $\mathcal{N}_{i + 1}$. They are also neighbours on $\Gamma$ if and only if some move of $L$ in direction $(-p, q)$ leads from a cell in $\mathcal{N}'_i$ to a cell in $\mathcal{N}'_{i + 1}$. Equivalently, if and only if some move of $L$ in direction $(-p, q)$ leads from a cell in $\mathcal{N}'_i$ to another cell in $\Gamma$. Let us see what we can learn about $\mathcal{N}_i$ and $\mathcal{N}_{i + 1}$ based on this.

Let $\Gamma'_\texttt{LR}$ and $\Gamma''_\texttt{LR}$ be the lower right subgrids of $\Gamma$ of sizes $p \times 2q$ and $2p \times q$, respectively. Then $\mathcal{N}_i$ intersects each one of them precisely once. Let $a'$ be the unique cell of $\mathcal{N}_i$ in $\Gamma'_\texttt{LR}$, and define $a''$ analogously. Of course, $a'$ and $a''$ are also the unique cells of $\mathcal{N}'_i$ in $\Gamma'_\texttt{LR}$ and $\Gamma''_\texttt{LR}$, respectively.

Every cell of $\mathcal{N}'_i$ is either nonstrictly above and to the left of $a'$, or nonstrictly above and to the left of $a''$. (Or perhaps both.) Therefore, some move of $L$ in direction $(-p, q)$ leads from a cell in $\mathcal{N}'_i$ to another cell in $\Gamma$ if and only if at least one of the moves of $L$ in direction $(-p, q)$ from $a'$ and $a''$ does so.

Consequently, nets $\mathcal{N}_i$ and $\mathcal{N}_{i + 1}$ form a stopper pair if and only if both cells $a' + (-p, q)$ and $a'' + (-p, q)$ are outside of $\Gamma$.

On the other hand, cell $a' + (-p, q)$ is outside of $\Gamma$ if and only if $a'$ is in the leftmost column of subgrid $\Gamma'_\texttt{LR}$, and cell $a'' + (-p, q)$ is outside of $\Gamma$ if and only if $a''$ is in the top row of subgrid $\Gamma''_\texttt{LR}$.

This information suffices to pin down $\mathcal{N}_i$ and $\mathcal{N}_{i + 1}$ precisely.

To this end, suppose first that $\mathcal{N}_i$ and $\mathcal{N}_{i + 1}$ do form a stopper pair, and so also that $a'$ is in the leftmost column of $\Gamma'_\texttt{LR}$ and $a''$ is in the top row of $\Gamma''_\texttt{LR}$.

Consider cell $a''' = a' + (q, p)$. Since $a'$ is in $\mathcal{N}_i$, so is $a'''$. Also, since $a'$ is in the leftmost column of $\Gamma'_\texttt{LR}$, we obtain that $a'''$ is in the leftmost column of $\Gamma''_\texttt{LR}$. Thus $a'''$ must be the unique cell where $\mathcal{N}_i$ intersects $\Gamma''_\texttt{LR}$. However, we already have a name for that cell: It is $a''$.

It follows that $a''$ and $a'''$ coincide and so cell $a'' = a'''$ lies both in the top row and in the leftmost column of $\Gamma''_\texttt{LR}$. Consequently, it must be the top left corner of $\Gamma''_\texttt{LR}$, and so also cell $a'$ must be the top left corner of $\Gamma'_\texttt{LR}$.

Conversely, the top left corners of $\Gamma'_\texttt{LR}$ and $\Gamma''_\texttt{LR}$ do, in fact, always belong to the same net.

Therefore, there is indeed precisely one stopper pair of nets. One of the nets in this pair is the net which contains the top left corners of $\Gamma'_\texttt{LR}$ and $\Gamma''_\texttt{LR}$, and the other one is its neighbour in direction $(-p, q)$.

We are almost done with the proof. By Lemma \ref{nets-connectedness}, the restriction of each net to $\Gamma$ is connected. And, since there is precisely one stopper pair of nets, by Lemma \ref{nets-cycle} all nets form a path with respect to neighbourhood on $\Gamma$. Therefore, the graph $\mathcal{L}_\Gamma \setminus \frac{q}{p}$ is connected, as needed. \end{proof}

\medbreak

We take a moment at this point to discuss, briefly, a couple of items related to Lemma \ref{three-slopes-connectedness}.

Curiously, when $p = 1$, the graph $\mathcal{L}_\Gamma \setminus \frac{q}{p}$ becomes a Hamiltonian path of $L$ on $\Gamma$. This Hamiltonian path was already known to fairy chess composer Thomas Dawson in the early twentieth century; see \cite{D}.

And, in the general case, from the proof of Lemma \ref{three-slopes-connectedness} we derive the following corollary. (Which is not connected to flexibility and rigidity in any way whatsoever.)

\medbreak

\begin{corollary} The diameter of $\mathcal{L}_\Gamma \setminus \frac{q}{p}$ is at least $2pq - 1$. \label{three-slopes-diameter} \end{corollary}

\medbreak

\begin{proof} Since we need to change nets at least this many times in order to reach from one net in the unique stopper pair to the other. \end{proof}

\medbreak

In essence, Corollary \ref{three-slopes-diameter} tells us that the diameter of $\mathcal{L}_\Gamma \setminus \frac{q}{p}$ is relatively large. When $p$ and $q$ are close together, it is at least about $\frac{1}{6}$ of the total number of edges in the graph. Intuitively, this means that $\mathcal{L}_\Gamma \setminus \frac{q}{p}$ is just barely connected.

With this, we return to our main topic of discussion.

Our technique of proving the connectedness of leaper graphs with one forbidden slope applies to other grids beyond $\Gamma$ as well. We discuss two such grids very briefly in the proof of Theorem \ref{rigidity-theta}.

Lemma \ref{three-slopes-connectedness} sums up almost all knowledge about the connectedness of forbidden-slope leaper graphs that we are going to need for Section \ref{rigid}, where we establish rigidity. In subsequent sections, however, where we study flexibility, we will require the following generalisation instead.

\medbreak

\begin{lemma} Let $h$ and $w$ be two positive integers such that $h \le p$ and $w \le p$, and let $A$ be the rectangular grid of height $2p + q - h$ and width $p + 2q - w$. Then the number of connected components of the graph $\mathcal{L}_A \setminus \frac{q}{p}$ is $hw$.

Moreover, let $A'_\texttt{LR}$ and $A''_\texttt{LR}$ be the lower right subgrids of $A$ of sizes $p \times 2q$ and $2p \times q$, respectively. Then the upper left subgrid of $A'_\texttt{LR}$ of size $h \times w$ contains precisely one representative out of each one of these connected components, and so does the upper left subgrid of $A''_\texttt{LR}$ of the same size. \label{three-slopes-components} \end{lemma}

\medbreak

\begin{proof} Fully analogous to the proof of Lemma \ref{three-slopes-connectedness}, except that this time around there are a total of $hw$ stopper pairs. \end{proof}

\medbreak

Of course, we obtain Lemma \ref{three-slopes-connectedness} by setting $h = w = 1$ in Lemma \ref{three-slopes-components}.

This completes our journey into the topic of the connectedness of leaper graphs with forbidden slopes.

\section{Rigidity} \label{rigid}

Our main goal in this section will be to prove the following theorem.

\medbreak

\begin{theorem} Let $\Phi$ be the square grid of side $2(p + q) - 1$. Then the leaper framework of $L$ on $\Phi$ is rigid. \label{rigidity-phi} \end{theorem}

\medbreak

By Lemma \ref{monotonicity}, Theorem \ref{rigidity-phi} implies that the leaper framework of $L$ is rigid on all larger grids as well. This confirms Solymosi and White's conjecture.

(In the form in which we stated it in the introduction. As we remarked in Section \ref{prelim}, Solymosi and White's original form of their conjecture was about infinitesimal rigidity rather than continuous rigidity. We handle this as well; see below.)

We are going to prove that $\mathcal{F}_\Phi$ is both globally rigid with distinctness, and also infinitesimally rigid. Of course, each one of these claims implies Theorem \ref{rigidity-phi} immediately. The two proofs have a lot in common, to the extent that they could reasonably be described as superficially distinct presentations of the same proof. We give both of them anyway since each one of these two types of rigidity is interesting in its own right.

We consider global rigidity first.

To begin with, we make some general observations about the leaper frameworks of $L$ on arbitrary grids. Then we apply these observations to grid $\Phi$ specifically.

We say that two edges of $L$ form a \emph{rhombus} if they are independent and they belong to some cycle of $L$ of length four. Equivalently, if they become the opposite sides of a nondegenerate rhombus, in the sense of Euclidean geometry, when we reinterpret them as segments; and, additionally, the other pair of opposite sides of that rhombus correspond to edges of $L$ as well.

We also say that two moves of $L$ form a rhombus if they point in the same direction and unorienting them results in two edges of $L$ which form a rhombus.

Let, then, $A$ be any grid. We say that two edges $e'$ and $e''$ of $L$ on $A$ are in the same \emph{rhombic class} of $L$ on $A$ if there exists some sequence $e_1 = e'$, $e_2$, \ldots, $e_n = e''$ of edges of $L$ on $A$ such that $e_i$ and $e_{i + 1}$ form a rhombus for all $i$. Of course, if two edges of $L$ are in the same rhombic class on any grid, then they are of the same slope.

We define the rhombic class of a move of $L$ on $A$ to be the same as the rhombic class of the edge which we obtain when we unorient that move. Thus if two moves of $L$ on $A$ are in the same rhombic class, then they point either in the same direction, or in opposite directions.

The notion of a rhombic class is the cornerstone on top of which we are going to build essentially all of our subsequent discussion of the flexibility and rigidity of leaper frameworks. Its full importance will become clear shortly.

Consider any placement $\mathbf{c}$ of $\mathcal{L}_A$ such that the points of $\mathbf{c}$ are pairwise distinct and $\mathbf{c}$ is equivalent to the canonical placement of $\mathcal{L}_A$. Then $\mathbf{c}$ determines an embedding of $\mathcal{L}_A$ into the plane. For simplicity, we write $\mathbf{c}(a)$ for the point which this embedding assigns to cell $a$ of $A$.

Let $e = a \to b$ be any move of $L$ on $A$. Then our embedding maps edge $ab$ of $L$ onto some segment, and it is natural to say that it also maps move $a \to b$ of $L$ onto the corresponding oriented segment. We say that the vector of that oriented segment \emph{realises} $e$ under $\mathbf{c}$. Formally, this vector is $\mathbf{c}(b) - \mathbf{c}(a)$.

\medbreak

\begin{lemma} Suppose that two moves $e'$ and $e''$ of $L$ on $A$ point in the same direction and belong to the same rhombic class. Then they are realised by the same vector under $\mathbf{c}$. \label{rhombic-vector} \end{lemma}

\medbreak

\begin{proof} Let $e' = a' \to b'$ and $e'' = a'' \to b''$.

Suppose first that $e'$ and $e''$ form a rhombus. Then, since all bars in a leaper framework are of the same length and points $\mathbf{c}(a')$, $\mathbf{c}(b')$, $\mathbf{c}(a'')$, and $\mathbf{c}(b'')$ are pairwise distinct, these four points must be the vertices of a rhombus (in the sense of Euclidean geometry) as well. Therefore, in this special case $e'$ and $e''$ are indeed realised by the same vector under $\mathbf{c}$.

Otherwise, if $e'$ and $e''$ do not form a rhombus directly, then let $e_1 = e'$, $e_2$, \ldots, $e_n = e''$ be some sequence of moves of $L$ on $A$ such as in the definition of a rhombic class, all pointing in the same direction. Then, by the first part of the proof, $e_i$ and $e_{i + 1}$ are realised by the same vector under $\mathbf{c}$ for all $i$. \end{proof}

\medbreak

We proceed to study the rhombic classes of $L$ on $A$. The following lemma shows exactly how our results of Section \ref{slope} relate to flexibility and rigidity.

\medbreak

\begin{lemma} Suppose that $p < w_A$ and $q < h_A$. Let $e' = a' \to b'$ and $e'' = a'' \to b''$ be two moves of $L$ on $A$, both of direction $(p, q)$, and let $B$ be the lower left subgrid of $A$ of size $(h_A - q) \times (w_A - p)$. Then $e'$ and $e''$ are in the same rhombic class of $L$ on $A$ if and only if $a'$ and $a''$ are in the same connected component of $\mathcal{L}_B \setminus \frac{q}{p}$. Analogous results hold for all other directions of the moves of $L$ as well. \label{rhombic-class} \end{lemma}

\medbreak

The condition that $p < w_A$ and $q < h_A$ is not actually particularly restrictive. If it is not satisfied, then there are no moves of $L$ on $A$ in direction $(p, q)$ at all.

\medbreak

\begin{proof} Suppose first that $a'$ and $a''$ are in the same connected component of $\mathcal{L}_B \setminus \frac{q}{p}$. Then there is some path from $a'$ to $a''$ in this graph, say $a_1a_2 \ldots a_n$, where $a_1 = a'$ and $a_n = a''$.

Consider the sequence of moves of $L$ on $A$ given by $e_i = a_i \to [a_i + (p, q)]$ for all $i$. Since no move in our path is of slope $\frac{q}{p}$, moves $e_i$ and $e_{i + 1}$ form a nondegenerate rhombus for all $i$. Therefore, this sequence is a witness for moves $e' = e_1$ and $e'' = e_n$ being in the same rhombic class of $L$ on $A$.

For a proof of the ``only if'' part, read the above argument backwards. \end{proof}

\medbreak

This completes our discussion of leaper frameworks in general, and we turn to grid $\Phi$ in particular.

The most important thing about $\Phi$, as far as flexibility and rigidity are concerned, is that there are really very few rhombic classes of $L$ on it.

\medbreak

\begin{lemma} Two edges of $L$ on $\Phi$ are in the same rhombic class if and only if they are of the same slope. \label{rhombic-phi} \end{lemma}

\medbreak

\begin{proof} The ``if'' part follows by Lemmas \ref{three-slopes-connectedness} and \ref{rhombic-class}, and the ``only if'' part is clear. \end{proof}

\medbreak

Consider any placement $\mathbf{c}$ of $\mathcal{L}_\Phi$ such that the points of $\mathbf{c}$ are pairwise distinct and $\mathbf{c}$ is equivalent to the canonical placement of $\mathcal{L}_\Phi$. To see that $\mathcal{F}_\Phi$ is globally rigid with distinctness, we must show that $\mathbf{c}$ is also congruent to the canonical placement of $\mathcal{L}_\Phi$.

By Lemmas \ref{rhombic-vector} and \ref{rhombic-phi}, all moves of $L$ on $\Phi$ which point in the same direction are realised by the same vector under $\mathbf{c}$. Intuitively, this means that $\mathbf{c}$ is more or less completely determined by some very small amount of information.

Of course, two moves of $L$ on $\Phi$ which point in opposite directions are always realised by opposite vectors. Thus we need to introduce notation for four vectors only. Suppose that all moves of $L$ on $\Phi$ in directions $(q, p)$, $(p, q)$, $(-p, q)$, and $(-q, p)$ are realised by vectors $v_1$, $v_2$, $v_3$, and $v_4$, respectively.

Let us see what we can learn about these vectors.

The main role in our argument from this point on will be taken up by one remarkably nice oriented cycle of $L$ on $\Phi$. We define it first as an oriented subgraph of $\mathcal{L}_\Phi$, and then we prove that this subgraph is in fact an oriented cycle.

(To be clear, to us a ``cycle'' is a simple closed walk, that is, one which repeats neither cells nor edges.)

Suppose that $\Phi$ is the grid $[1; 2(p + q) - 1] \times [1; 2(p + q) - 1]$. We define $C_\Phi$ as the subgraph of $\mathcal{L}_\Phi$ formed by the following moves of $L$ on $\Phi$:

(a) All moves in direction $(q, p)$ which start from a cell of the form $(2i - 1, 1)$, where $1 \le i \le p$.

(b) All moves in direction $(q, -p)$ which lead to a cell of the form $(2i - 1, 1)$, where $q + 1 \le i \le p + q$, .

(c) All moves in direction $(-p, q)$ which start from a cell of the form $(2i - 1, 1)$, where $p + 1 \le i \le p + q$.

(d) Lastly, all moves in direction $(-p, -q)$ which lead to a cell of the form $(2i - 1, 1)$, where $1 \le i \le q$.

\begin{figure}[t] \centering \includegraphics[scale=\figuresScale]{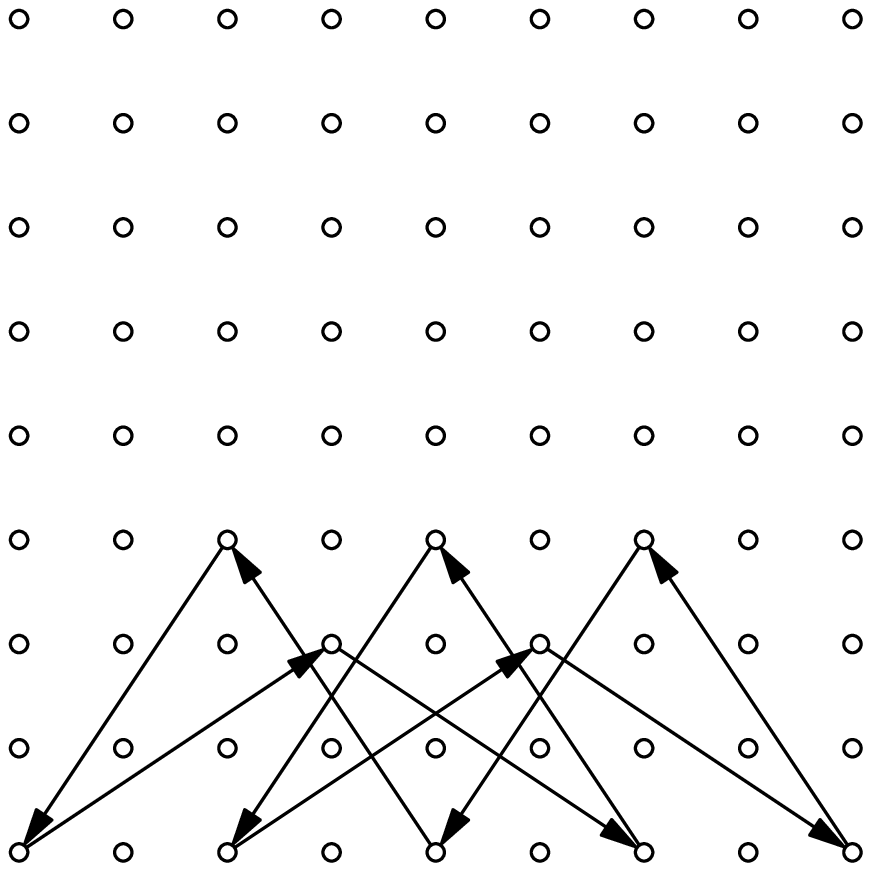} \caption{} \label{zebra-c-phi} \end{figure}

For example, Figure \ref{zebra-c-phi} shows $C_\Phi$ for the zebra.

It is straightforward to verify that every vertex of $C_\Phi$ is of in-degree one as well as of out-degree one. Therefore, $C_\Phi$ is the disjoint union of several oriented cycles.

(This is actually all we need to know about the structure of $C_\Phi$ for the rest of our argument to go through. But we are going to derive the complete description anyway.)

Let $C$ be the connected component of $C_\Phi$ which contains the lower left corner of $\Phi$. Suppose that $C$ contains a total of $\alpha$ moves of $L$ in direction $(q, p)$ and a total of $\beta$ moves of $L$ in direction $(-p, q)$.

Observe that, in $C_\Phi$, every endpoint of a move of $L$ in direction $(q, p)$ is also the starting point of another move of $L$ in direction $(q, -p)$. Therefore, this is true of $C$ as well, and so $C$ contains a total of $\alpha$ moves of $L$ in direction $(q, -p)$. Analogously, $C$ contains a total of $\beta$ moves of $L$ in direction $(-p, -q)$. Lastly, since $C$ contains the lower left corner of $\Phi$, both of $\alpha$ and $\beta$ are nonzero.

By summation along $C$, we obtain \[\alpha(q, p) + \alpha(q, -p) + \beta(-p, q) + \beta(-p, -q) = \mathbf{0}.\]

Therefore, $\alpha q = \beta p$. Since $p$ and $q$ are relatively prime, this means that $\alpha$ is a multiple of $p$ and $\beta$ is a multiple of $q$. Moreover, since both of $\alpha$ and $\beta$ are nonzero, it follows that $\alpha \ge p$ and $\beta \ge q$.

On the other hand, the entire $C_\Phi$ contains only $p$ moves of $L$ in each of directions $(q, p)$ and $(q, -p)$, as well as only $q$ moves of $L$ in each of directions $(-p, q)$ and $(-p, -q)$. Therefore, we must have that $\alpha = p$, $\beta = q$, and $C$ coincides with $C_\Phi$.

This completes our proof that $C_\Phi$ is an oriented cycle of $L$ on $\Phi$.

Our embedding of $\mathcal{L}_\Phi$ into the plane contains a subembedding of $C_\Phi$. Since $C_\Phi$ contains $p$ moves of $L$ in direction $(q, p)$, its embedding into the plane contains $p$ copies of vector $v_1$. Same goes for vectors $v_2$, $v_3$, and $v_4$. Then, by summation along the embedding of $C_\Phi$ into the plane, we obtain \[pv_1 - qv_2 + qv_3 - pv_4 = \mathbf{0}.\]

Let $C'_\Phi$ be the reflection of $C_\Phi$ in the line $x = y$. (Equivalently, in the up-and-to-the-right diagonal of $\Phi$. Reflection in the other diagonal of $\Phi$, or rotation by $\pm 90^\circ$ about the center of $\Phi$, would all have worked just as well.) Analogous reasoning applies to $C'_\Phi$, too, and in this way we obtain \[-qv_1 + pv_2 + pv_3 - qv_4 = \mathbf{0}.\]

This information suffices to pin down $v_1$, $v_2$, $v_3$, and $v_4$ almost precisely, up to rotation. Let us see how the calculations work out.

Define \begin{align*} u' &= \frac{v_1 - v_4}{2q} = \frac{v_2 - v_3}{2p} \text{ and }\\ u'' &= \frac{v_1 + v_4}{2p} = \frac{v_2 + v_3}{2q}. \end{align*}

Since all bars in a leaper framework are of the same length, so are vectors $v_1$, $v_2$, $v_3$, and $v_4$. Therefore, \begin{align*} u' \cdot u'' &= \frac{1}{4pq} \cdot 2qu' \cdot 2pu''\\ &= \frac{1}{4pq}(v_1 - v_4)(v_1 + v_4)\\ &= \frac{1}{4pq}(|v_1|^2 - |v_4|^2)\\ &= 0. \end{align*}

In other words, vectors $u'$ and $u''$ are orthogonal.

Consequently, we have that \begin{align*} p^2 + q^2 &= |v_1|^2\\ &= (qu' + pu'')^2\\ &= q^2|u'|^2 + p^2|u''|^2 \end{align*} and \begin{align*} p^2 + q^2 &= |v_2|^2\\ &= (pu' + qu'')^2\\ &= p^2|u'|^2 + q^2|u''|^2. \end{align*}

Therefore, both of $u'$ and $u''$ are unit vectors.

Without loss of generality, rotate $\mathbf{c}$ until vector $u'$ comes to point to the right and vector $u''$ comes to point upwards. Then vectors $v_1$, $v_2$, $v_3$, and $v_4$ become simply vectors $(q, p)$, $(p, q)$, $(-p, q)$, and $(-q, p)$, respectively.

Our proof of Theorem \ref{rigidity-phi} is almost complete. There is just one more detail which we need to work our way through: That vectors $v_1$, $v_2$, $v_3$, and $v_4$ determine $\mathbf{c}$ essentially uniquely, up to translation. To this end, it suffices to verify that graph $\mathcal{L}_\Phi$ is connected.

This follows immediately by Knuth's description of all grids on which the leaper graph of $L$ is connected. Alternatively, a more self-contained approach could be as follows. There is at least one cell of $\Phi$ which is incident with at least one edge of $L$ on $\Phi$ of each slope. Consequently, $\mathcal{L}_\Phi$ is connected by Lemma \ref{rhombic-phi}.

Thus $\mathbf{c}$ is indeed congruent to the canonical placement of the leaper graph of $L$ on $\Phi$. This completes our proof that the leaper framework of $L$ on $\Phi$ is globally rigid with distinctness.

We go on to infinitesimal rigidity.

Since the two proofs are very similar, we will do as follows. For the most part, we are simply going to stroll along the lines of our proof of global rigidity. However, whenever substantial changes are necessary, we will examine them in detail.

Again, we discuss the leaper frameworks of $L$ on arbitrary grids first, and then we apply our observations to grid $\Phi$ in particular.

Let $A$ be any grid and let $I$ be any infinitesimal motion of the leaper framework of $L$ on $A$.

Let $e = a \to b$ be any move of $L$ on $A$. By analogy with global rigidity, we say that vector $I(b) - I(a)$ \emph{realises} $e$ under $I$. By the definition of an infinitesimal motion, the direction of $e$ and the vector which realises $e$ are orthogonal for all $e$.

An analogue of Lemma \ref{rhombic-vector} holds for infinitesimal rigidity as well. The proof relies on the following useful lemma.

\medbreak

\begin{lemma} Suppose that cells $a_1$, $a_2$, $a_3$, and $a_4$ of $A$ form a rhombus of $L$ on $A$, in this order. Then $I(a_1) + I(a_3) = I(a_2) + I(a_4)$. \label{rhombic-inf-auxiliary} \end{lemma}

\medbreak

This is a special case of a more general result found in \cite{SW} and \cite{W}. For completeness, we include a retelling of the proof given there.

\medbreak

\begin{proof} We have that \begin{align*} &\phantom{= \mbox{ }} (a_1 - a_2)(I(a_1) - I(a_2) + I(a_3) - I(a_4))\\ &= (a_1 - a_2)(I(a_1) - I(a_2)) + (a_1 - a_2)(I(a_3) - I(a_4))\\ &= (a_1 - a_2)(I(a_1) - I(a_2)) - (a_3 - a_4)(I(a_3) - I(a_4))\\ &= \mathbf{0}. \end{align*}

Analogously, \[(a_1 - a_4)(I(a_1) - I(a_2) + I(a_3) - I(a_4)) = \mathbf{0}\] as well.

Thus vector $I(a_1) - I(a_2) + I(a_3) - I(a_4)$ is orthogonal to two noncollinear vectors. Therefore, it must be the zero vector. \end{proof}

\medbreak

We are ready to prove the analogue of Lemma \ref{rhombic-vector} for infinitesimal flexibility.

\medbreak

\begin{lemma} Suppose that two moves $e'$ and $e''$ of $L$ on $A$ point in the same direction and belong to the same rhombic class. Then they are realised by the same vector under $I$. \label{rhombic-inf-vector} \end{lemma}

\medbreak

\begin{proof} By Lemma \ref{rhombic-inf-auxiliary}, analogously to the proof of Lemma \ref{rhombic-vector}. \end{proof}

\medbreak

With this, we return to the concrete grid $\Phi$.

Define vectors $v_1$, $v_2$, $v_3$, and $v_4$ exactly as in the proof of global rigidity, only relative to $I$ instead of $\mathbf{c}$.

Consider $v_I = (v_1, v_2, v_3, v_4)$ as a vector in $\mathbb{R}^8$. Since all infinitesimal motions of $\mathcal{F}_\Phi$ form a vector space, so do all vectors $v_I$ as $I$ ranges over all infinitesimal motions of $\mathcal{F}_\Phi$. We denote this vector space by $\operatorname{Inf}'(\mathcal{F}_\Phi)$. Of course, it is a subspace of $\mathbb{R}^8$.

Just as in our proof of global rigidity, $v_I$ determines $I$ essentially uniquely, up to translation. Factoring out all translations shaves two degrees of freedom off of $\operatorname{Inf}(\mathcal{F}_\Phi)$. Therefore, the dimension of $\operatorname{Inf}'(\mathcal{F}_\Phi)$ is precisely two less than the dimension of $\operatorname{Inf}(\mathcal{F}_\Phi)$.

Thus, for infinitesimal rigidity in the setting of Theorem \ref{rigidity-phi}, it suffices to prove that $\operatorname{Inf}'(\mathcal{F}_\Phi)$ is of dimension at most one.

(We know already that its dimension is at least one, corresponding to the rotations of $\mathcal{F}_\Phi$. So it will follow from our proof that the dimension of $\operatorname{Inf}'(\mathcal{F}_\Phi)$ is in fact exactly one.)

Let us see what constraints $v_I$ must satisfy. For all $i$, let $v_i = (x_i, y_i)$, so that $v_I = (x_1, y_1, x_2, y_2, x_3, y_3, x_4, y_4)$.

By the definition of an infinitesimal motion, each $v_i$ must be orthogonal to its associated direction. For example, vectors $v_1$ and $(q, p)$ must be orthogonal. This gives us the four homogeneous linear constraints \begin{align*} qx_1 + py_1 &= 0,\\ px_2 + qy_2 &= 0,\\ -px_3 + qy_3 &= 0, \text{ and }\\ -qx_4 + py_4 &= 0. \end{align*}

On the other hand, the sum of the realisations of the moves of $L$ along any oriented cycle telescopes, and so must equal the zero vector. Applying this observation to $C_\Phi$ and $C'_\Phi$, we see that that the $v_i$ must satisfy the exact same constraints as with global rigidity. In terms of the $x_i$ and $y_i$, these constraints become \begin{align*} px_1 - qx_2 + qx_3 - px_4 &= 0,\\ py_1 - qy_2 + qy_3 - py_4 &= 0,\\ -qx_1 + px_2 + px_3 - qx_4 &= 0, \text{ and }\\ -qy_1 + py_2 + py_3 - qy_4 &= 0. \end{align*}

Thus we obtain four more homogeneous linear constraints for the $x_i$ and $y_i$.

The matrix of the coefficients of all eight constraints is \[M_\Phi = \begin{pmatrix} q & p & 0 & 0 & 0 & 0 & 0 & 0\\ 0 & 0 & p & q & 0 & 0 & 0 & 0\\ 0 & 0 & 0 & 0 & -p & q & 0 & 0\\ 0 & 0 & 0 & 0 & 0 & 0 & -q & p\\ p & 0 & -q & 0 & q & 0 & -p & 0\\ 0 & p & 0 & -q & 0 & q & 0 & -p\\ -q & 0 & p & 0 & p & 0 & -q & 0\\ 0 & -q & 0 & p & 0 & p & 0 & -q \end{pmatrix}.\]

Then \[\dim \operatorname{Inf}'(\mathcal{F}_\Phi) \le 8 - \operatorname{rank} M_\Phi.\]

Therefore, all that we are left to do is prove that the rank of $M_\Phi$ is at least seven.

(We know already that the rank of $M_\Phi$ is at most seven since the dimension of $\operatorname{Inf}'(\mathcal{F}_\Phi)$ is at least one. Indeed, the rows of $M_\Phi$ are linearly dependent with coefficients $-1$, $1$, $-1$, $1$, $0$, $1$, $-1$, and $0$, respectively. So it will follow from our proof that the rank of $M_\Phi$ is in fact exactly seven.)

To this end, consider any linear combination of the rows of $M_\Phi$ with coefficients $\alpha_1$, $\alpha_2$, $\alpha_3$, $\alpha_4$, $\beta_1$, $\beta_2$, $\beta_3$, and $\beta_4$, respectively, which equals the zero vector. Suppose that $\beta_2 = 0$; it suffices to show that, under this assumption, all other coefficients must equal zero as well. (We want to keep rows $1$--$4$ because this makes the calculations simpler, and omitting row $5$ decreases rank.)

By columns $2$, $4$, $6$, and $8$ of $M_\Phi$, we have that \[\alpha_1 : q^2 = \alpha_2 : (-p^2) = \alpha_3 : (-p^2) = \alpha_4 : q^2 = \beta_4 : pq.\]

Thus let $\alpha_1 = \alpha_4 = \gamma q^2$, $\alpha_2 = \alpha_3 = -\gamma p^2$, and $\beta_4 = \gamma pq$ for some real parameter $\gamma$.

Then columns $1$ and $3$ of $M_\Phi$ tell us that \begin{align*} \gamma q^3 + \beta_1 p - \beta_3 q &= 0 \text{ and }\\ -\gamma p^3 - \beta_1 q + \beta_3 p &= 0, \end{align*} respectively. Since $p \neq q$, this is enough to calculate $\beta_1$ and $\beta_3$; they work out to $\beta_1 = \gamma pq$ and $\beta_3 = \gamma (p^2 + q^2)$.

However, substituting these values for $\beta_1$ and $\beta_3$ in column $5$ of $M_\Phi$, and taking into account the fact that $p > 0$ and $q > 0$, yields $\gamma = 0$. Therefore, all seven of $\alpha_1$, $\alpha_2$, $\alpha_3$, $\alpha_4$, $\beta_1$, $\beta_3$, and $\beta_4$ must equal zero as well, and we are done.

This completes our proof that the framework of $L$ on $\Phi$ is infinitesimally rigid.

Both of our proofs of Theorem \ref{rigidity-phi} go through on all larger grids as well. The only technicality which is perhaps worth mentioning is that graph $\mathcal{L}_A \setminus \frac{q}{p}$ is connected on all grids $A$ larger than $\Gamma$ as well. This follows by induction on the sides of the grid, or alternatively simply by observing that on all larger grids there are no stopper pairs of nets at all. Thus we arrive at the following corollary.

\medbreak

\begin{corollary} The leaper framework of $L$ is both globally rigid with distinctness and also infinitesimally rigid on all grids larger than or equal to $\Phi$. \label{rigidity-phi-stronger} \end{corollary}

\medbreak

An alternative approach to Corollary \ref{rigidity-phi-stronger} could be to work out what the analogues of Lemma \ref{monotonicity} for global rigidity and infinitesimal rigidity should be.

Corollary \ref{rigidity-phi-stronger} completely settles Solymosi and White's conjecture in the form in which they stated it originally in \cite{SW} and \cite{W}.

One might naturally wonder if Corollary \ref{rigidity-phi-stronger} describes all grids on which the leaper framework of $L$ is rigid. To finish this section, we show that it does not.

\medbreak

\begin{theorem} Let $\Theta$ be the rectangular grid of height $p + 2q$ and width $2(p + q)$. Then the leaper framework of $L$ on $\Theta$ is rigid. \label{rigidity-theta} \end{theorem}

\medbreak

Before we go on to the proof, let us first have a quick look at the implications of this theorem. When $p \ge 2$, $\Theta$ is shorter than $\Phi$. Hence, the following corollary.

\medbreak

\begin{corollary} Suppose that $p \ge 2$. Then there is at least one rectangular grid $A$ such that $\Phi$ is not smaller than or equal to $A$, and yet the leaper framework of $L$ on $A$ is rigid. \label{incompleteness} \end{corollary}

\medbreak

In other words, there are infinitely many leapers $L$ such that Corollary \ref{rigidity-phi-stronger} does not describe all grids on which the leaper framework of $L$ is rigid.

With this, we move on to the proof.

Both of our proofs for Theorem \ref{rigidity-phi} mostly work for Theorem \ref{rigidity-theta}, too. In particular, an analogue of Corollary \ref{rigidity-phi-stronger} holds for grid $\Theta$ as well. Still, there are some important changes which we need to make, and we proceed to examine them in detail.

To begin with, we must adapt Lemma \ref{three-slopes-connectedness}. Let $\Lambda'$ be the rectangular grid of height $p + q$ and width $p + 2q$, and let $\Lambda''$ be the rectangular grid of height $2q$ and width $2p + q$. Then both graphs $\mathcal{L}_{\Lambda'} \setminus \frac{q}{p}$ and $\mathcal{L}_{\Lambda''} \setminus \frac{p}{q}$ are connected.

The proof is fully analogous to the proof of Lemma \ref{three-slopes-connectedness}, except that this time around there are no stopper pairs of nets at all, neither on grid $\Lambda'$ nor on grid $\Lambda''$.

Equipped with this knowledge, we establish the analogue of Lemma \ref{rhombic-phi} for $\Theta$ exactly as we established the original Lemma \ref{rhombic-phi}. Thus the rhombic class of an edge of $L$ on $\Theta$ is entirely determined by the edge's slope.

Oriented cycle $C_\Phi$ is a subgraph of $\mathcal{L}_\Theta$, too, and so this part of the proof remains the same. Oriented cycle $C'_\Phi$, though, is not a subgraph of $\mathcal{L}_\Theta$ when $p \ge 2$; it is too tall.

Luckily, this is not as much of an obstacle as it might seem. We used $C_\Phi$ and its reflection for $\Phi$ because it is a natural oriented cycle to come up with in this situation, and also because it is easy to describe. There are other subgraphs that we could have used, which are more complicated to describe but take up a lot less room. Here is one of them.

Suppose that $\Theta$ is the grid $[1; 2(p + q)] \times [1; p + 2q]$. When $p \ge 2$, we define $C_\Theta$ as the subgraph of $\mathcal{L}_\Theta$ formed by the following moves of $L$ on $\Theta$:

(a) All moves in direction $(q, p)$ which start from a cell of the form $(1, i)$, where $1 \le i \le 2q$.

(b) All moves in direction $(-q, p)$ which start from a cell of the form $(q + 1, i)$, where $1 \le i \le 2q$.

(c) All moves in direction $(p, -q)$ which start from a cell of the form $(1, i)$, where $2q + 1 \le i \le p + 2q$.

(d) All moves in direction $(-p, -q)$ which lead to a cell of the form $(1, i)$, where $1 \le i \le p$.

(e) All moves in direction $(-p, -q)$ which start from a cell of the form $(q + 1, i)$, where $2q + 1 \le i \le p + 2q$.

(f) Lastly, all moves in direction $(p, -q)$ which lead to a cell of the form $(q + 1, i)$, where $1 \le i \le p$.

\begin{figure}[t] \centering \includegraphics[scale=\figuresScale]{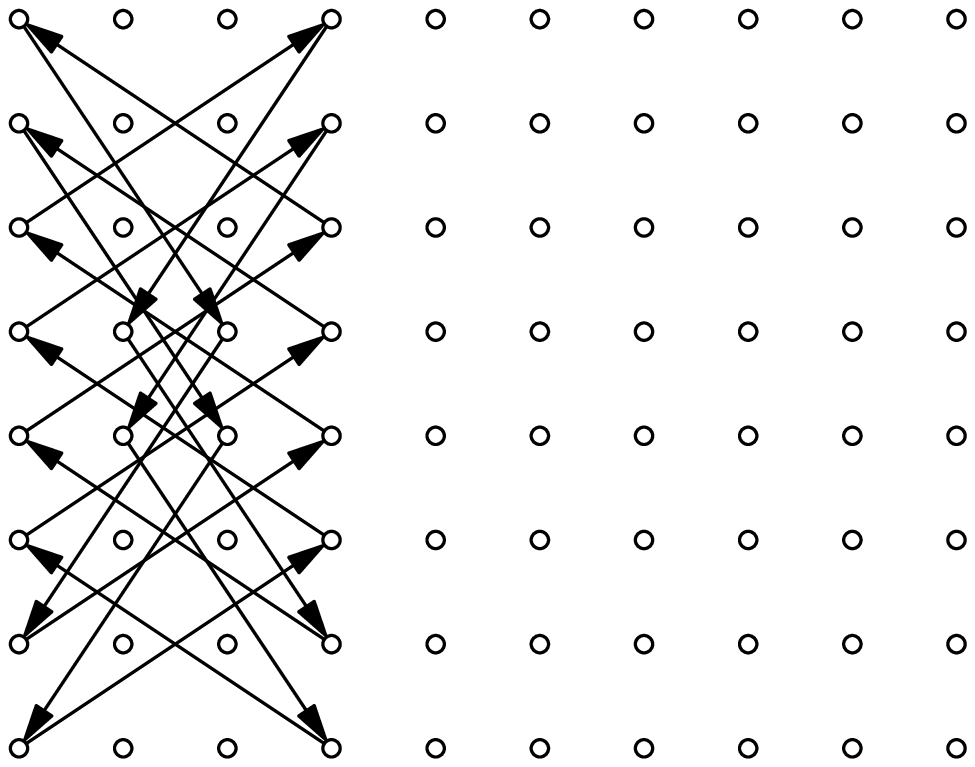} \caption{} \label{zebra-c-theta} \end{figure}

For example, Figure \ref{zebra-c-theta} shows $C_\Theta$ for the zebra.

We prove that $C_\Theta$ is the disjoint union of several oriented cycles precisely as we did the same thing for $C_\Phi$. Just as with $C_\Phi$, this is really all that we need for the rest of our argument to go through.

(Note that, this time around, we do not claim that $C_\Theta$ is a single oriented cycle. Indeed, it is no such thing. It is not too difficult to describe the structure of $C_\Theta$ completely; it always has two connected components, which are either reflections or translation copies of each other depending on the parities of $p$ and $q$. However, a proof of this fact would be too much of a distraction from our main focus.)

The rest of what we did for Theorem \ref{rigidity-phi} does not require any adaptations at all. This completes our proof of Theorem \ref{rigidity-theta}, as well as our discussion of rigidity.

\section{Infinitesimal Flexibility} \label{inf-flex}

Solymosi and White's conjecture is, as we said in the introduction, sharp with respect to the size of the grid: On all grids smaller than $\Phi$, the leaper framework of $L$ is flexible. Our goal for most of the rest of the present work will be to prove this.

Roughly speaking, the proof consists of three parts which belong to three different fields of mathematics: combinatorics, linear algebra, and mathematical analysis. We handle each part somewhat separately, with the obvious caveat that a perfectly neat division into fields is not really possible. Still, we do most of the combinatorics of the proof in the present section; most of the mathematical analysis in Section \ref{method}; and most of the linear algebra in Section \ref{flex-i}. Section \ref{flex-ii} then ties up some loose ends.

In this section, we focus on infinitesimal flexibility. By way of an intermediate step towards our long-term goal, we are going to prove the following theorem.

\medbreak

\begin{theorem} Let $\Psi$ be the rectangular grid of height $2(p + q) - 2$ and width $2(p + q) - 1$. Then the leaper framework of $L$ on $\Psi$ is infinitesimally flexible. \label{inf-flexibility-psi} \end{theorem}

\medbreak

Of course, Theorem \ref{inf-flexibility-psi} is also an immediate corollary of Theorem \ref{flexibility-psi}. We prove it separately anyway because this is a good way to break the argument down into manageable parts, and also because seeing how some of the main ideas of the proof of Theorem \ref{flexibility-psi} play out in this simpler setting will help throw them into starker relief.

Note that graph $\mathcal{L}_\Psi$ is connected by Knuth's description of all grids on which the leaper graph of $L$ is connected. Thus the question of the infinitesimal and continuous flexibility of $\mathcal{F}_\Psi$ is indeed nontrivial.

The case $p = 1$ of Theorem \ref{inf-flexibility-psi} is immediate by Lemma \ref{linear-bound}. (In fact, Lemma \ref{linear-bound} does not resolve any other cases of Theorem \ref{inf-flexibility-psi} at all.)

Suppose, then, from this point on throughout the rest of our discussion of $\Psi$ in this section, that $p \ge 2$.

Our first order of business, just as with rigidity, will be to examine the rhombic classes of $L$ on $\Psi$.

\medbreak

\begin{lemma} When $p \ge 2$, there are a total of eight rhombic classes of $L$ on $\Psi$, two per each slope. \label{rhombic-psi} \end{lemma}

\medbreak

\begin{proof} By Lemmas \ref{three-slopes-components} and \ref{rhombic-class}. \end{proof}

\medbreak

Note that our proof of Lemma \ref{rhombic-psi} does not work when $p = 1$, as then we cannot apply Lemma \ref{three-slopes-components}. Indeed, in the case of $p = 1$, the number of rhombic classes becomes larger.

We go on to introduce some notation for the rhombic classes of $L$ on $\Psi$.

Let $R_1$ and $R_2$ be the two rhombic classes of $L$ on $\Psi$ for slope $\frac{q}{p}$.

Since reflection with respect to the vertical axis of symmetry of $\Psi$ swaps slopes $\frac{q}{p}$ and $-\frac{q}{p}$, the rhombic classes of $L$ on $\Psi$ for slope $-\frac{q}{p}$ are reflections of the rhombic classes of $L$ on $\Psi$ for slope $\frac{q}{p}$. Let rhombic classes $R_3$ and $R_4$ be the reflections of rhombic classes $R_1$ and $R_2$, respectively.

We handle slopes $\frac{p}{q}$ and $-\frac{p}{q}$ analogously. Let $R_5$ and $R_6$ be the rhombic classes of $L$ on $\Psi$ for slope $\frac{p}{q}$, and let $R_7$ and $R_8$ be their reflections, respectively, in the vertical axis of symmetry of $\Psi$, so that $R_7$ and $R_8$ are the rhombic classes of $L$ on $\Psi$ for slope $-\frac{p}{q}$.

We take a moment here to discuss one related question which is not a required part of our argument, but will help us understand the structure of $\mathcal{L}_\Psi$ and its rhombic classes a little better.

Reflection with respect to the horizontal axis of symmetry of $\Psi$ also swaps slopes $\frac{q}{p}$ and $-\frac{q}{p}$, just as reflection with respect to the vertical axis of symmetry of $\Psi$ does. Therefore, horizontal reflection swaps the two pairs of rhombic classes $\{R_1, R_2\}$ and $\{R_3, R_4\}$ as well. Does it swap $R_1$ with $R_3$ and $R_2$ with $R_4$, or the other way around? Of course, the same question applies to rhombic classes $R_5$, $R_6$, $R_7$, and $R_8$.

Some thought shows that this is really a question about the connected components of a certain forbidden-slope leaper graph. Let $\Pi'$ be the rectangular grid of height $h_\Psi - q = 2p + q - 2$ and width $w_\Psi - p = p + 2q - 1$, and let $\Pi''$ be the rectangular grid of height $h_\Psi - p = p + 2q - 2$ and width $w_\Psi - q = 2p + q - 1$. We consider graph $\mathcal{L}_{\Pi'} \setminus \frac{q}{p}$ in detail, and the case of graph $\mathcal{L}_{\Pi''} \setminus \frac{p}{q}$ is analogous.

By Lemma \ref{three-slopes-components}, $\mathcal{L}_{\Pi'} \setminus \frac{q}{p}$ has two connected components. Consider the central symmetry $\sigma$ with respect to the center of $\Pi'$. Since $\sigma$ preserves $\mathcal{L}_{\Pi'} \setminus \frac{q}{p}$, either $\sigma$ preserves each connected component of $\mathcal{L}_{\Pi'} \setminus \frac{q}{p}$, or $\sigma$ swaps the two components.

If the former, then horizontal reflection on $\Psi$ swaps $R_1$ with $R_3$ and $R_2$ with $R_4$. Otherwise, if the latter, then horizontal reflection on $\Psi$ swaps $R_1$ with $R_4$ and $R_2$ with $R_3$.

The answer to our question, then, is given by the following lemma.

\medbreak

\begin{lemma} When $p \ge 2$, each connected component of $\mathcal{L}_{\Pi'} \setminus \frac{q}{p}$ is centrally symmetric with respect to the center of $\Pi'$. Same goes for $\mathcal{L}_{\Pi''} \setminus \frac{p}{q}$. \label{three-slopes-symmetry} \end{lemma}

\medbreak

\begin{proof} Suppose not, for the sake of contradiction.

Let the two connected components be $K'$ and $K''$. Then $\sigma$ swaps $K'$ and $K''$; consequently, $K'$ and $K''$ span the same number of nets.

On the other hand, the nets spanned by each one of $K'$ and $K''$ must form a continuous subsequence of $\mathcal{N}_1$, $\mathcal{N}_2$, \ldots, $\mathcal{N}_{2pq}$. Suppose, without loss of generality, that $K'$ spans nets $\mathcal{N}_1$, $\mathcal{N}_2$, \ldots, $\mathcal{N}_{pq}$ and $K''$ spans nets $\mathcal{N}_{pq + 1}$, $\mathcal{N}_{pq + 2}$, \ldots, $\mathcal{N}_{2pq}$.

Then the two stopper pairs of nets on $\Pi'$ must be $\{\mathcal{N}_1, \mathcal{N}_{2pq}\}$ and $\{\mathcal{N}_{pq}, \mathcal{N}_{pq + 1}\}$.

By the proof of Lemma \ref{nets-cycle}, the cells of the union of $\mathcal{N}_{pq}$ and $\mathcal{N}_{2pq}$ form a translation copy of the lattice generated by vectors $(q, 0)$ and $(0, p)$. Consequently, when $p \ge 2$, no two cells in the union of $\mathcal{N}_{pq}$ and $\mathcal{N}_{2pq}$ are one unit apart.

On the other hand, consider the representatives of $K'$ and $K''$ which Lemma \ref{three-slopes-components} gives us. They are one unit apart exactly, and by the proof of Lemma \ref{three-slopes-components} one of them belongs to $\mathcal{N}_{pq}$ and the other one to $\mathcal{N}_{2pq}$.

We have arrived at a contradiction. \end{proof}

\medbreak

By Lemma \ref{three-slopes-symmetry}, each rhombic class of $L$ on $\Psi$ is centrally symmetric with respect to the center of $\Psi$. In particular, horizontal reflection on $\Psi$ swaps $R_1$ with $R_3$, $R_2$ with $R_4$, $R_5$ with $R_7$, and $R_6$ with $R_8$, just as vertical reflection on $\Psi$ does.

This completes our digression into the structure of the rhombic classes of $L$ on $\Psi$.

Similarly to what we did in Section \ref{rigid} for infinitesimal rigidity, we proceed to discuss the infinitesimal flexibility of the leaper frameworks of $L$ on arbitrary grids. This discussion will occupy most of the remainder of this section, and we will reuse large portions of it in Section \ref{method}. Once we are done, we will apply our findings to the concrete grid $\Psi$.

Note that we drop our assumption that $p \ge 2$ when we consider arbitrary grids, and that we pick it up again as soon as we return back to $\Psi$.

In the beginning, our analysis of infinitesimal flexibility is going to follow closely along the lines of our earlier analysis of infinitesimal rigidity in Section \ref{rigid}. However, at some point the two paths will diverge.

Let, then, $A$ be any grid. Suppose that $\mathcal{L}_A$ is connected, and let $R_1$, $R_2$, \ldots, $R_D$ be the rhombic classes of $L$ on $A$.

Let $I$ be any infinitesimal motion of $\mathcal{F}_A$. By Lemma \ref{rhombic-inf-vector}, if two moves of $L$ on $A$ point in the same direction and belong to the same rhombic class, then they are realised by the same vector under $I$.

For each rhombic class $R_i$, fix one of the two directions of its moves; say, $d_i$. Then, for all $i$, let all moves of direction $d_i$ in rhombic class $R_i$ be realised by vector $v_i$ under $I$. Of course, this means that all moves of direction $-d_i$ in rhombic class $R_i$ are realised by vector $-v_i$ under $I$.

Just as in Section \ref{rigid}, let $v_i = (x_i, y_i)$ for all $i$, and define $v_I = (x_1, y_1, x_2, y_2, \ldots,\allowbreak x_D, y_D)$. Then we say that vector $v_I$ is the \emph{summary} of the infinitesimal motion $I$ of $\mathcal{F}_A$. Once again, since all infinitesimal motions of $\mathcal{F}_A$ form a vector space, so do their summaries. We denote this vector space by $\operatorname{Inf}'(\mathcal{F}_A)$.

Since $\mathcal{L}_A$ is connected by assumption, $v_I$ determines $I$ almost uniquely, up to translation. Just as in Section \ref{rigid}, it follows from this that the dimension of $\operatorname{Inf}'(\mathcal{F}_A)$ is precisely two less than the dimension of $\operatorname{Inf}(\mathcal{F}_A)$. Consequently, $\mathcal{F}_A$ is infinitesimally flexible if and only if the dimension of $\operatorname{Inf}'(\mathcal{F}_A)$ is at least two.

Thus we set about to investigate $\operatorname{Inf}'(\mathcal{F}_A)$.

Consider any vector $v$ in $\mathbb{R}^{2D}$. When is it the summary of some infinitesimal motion $I_v$ of $\mathcal{F}_A$?

In Section \ref{rigid}, we derived some necessary conditions which all summary vectors must satisfy, in the special case of grid $\Phi$. Let us see how these conditions generalise. Write, as above, $v = (x_1, y_1, x_2, y_2, \ldots, x_D, y_D)$, and set $v_i = (x_i, y_i)$ for all $i$.

To begin with, by the definition of an infinitesimal motion, $v_i$ must be orthogonal to $d_i$ for all $i$. This amounts to a total of $D$ homogeneous linear constraints on the $x_i$ and $y_i$. We call them the \emph{directional constraints} on $v$.

Similarly to what we did in Section \ref{rigid}, the oriented cycles of $L$ on $A$ impose some homogeneous linear constraints on the $x_i$ and $y_i$ as well.

In Section \ref{rigid}, our goal was to estimate the dimension of $\operatorname{Inf}'(\mathcal{F}_A)$ from above. To this end, it was sufficient to examine just a couple of oriented cycles, since from them we learned enough to push the dimension of $\operatorname{Inf}'(\mathcal{F}_A)$ as low as we needed it to go.

This time around, however, we want to estimate the dimension of $\operatorname{Inf}'(\mathcal{F}_A)$ from below. Equivalently, we want to show that there exist many linearly independent vectors $v$ in $\mathbb{R}^{2D}$ which satisfy all conditions that we could possibly impose on them. Therefore, this time around, we are going to have to keep track of all oriented cycles of $L$ on $A$.

Consider, then, any oriented cycle $C$ of $L$ on $A$.

We define the \emph{content} $\omega_i$ of $C$ in rhombic class $R_i$ to be the number of edges in $R_i$ which $C$ traverses in direction $d_i$ minus the number of edges in $R_i$ which $C$ traverses in direction $-d_i$. Equivalently, as we trace $C$, every edge of rhombic class $R_i$ which we traverse in direction $\varepsilon d_i$ contributes $\varepsilon$ to $\omega_i$. Then we define the \emph{weight} of $C$ to be vector $(\omega_1, \omega_2, \ldots, \omega_D)$, and we denote it by $\operatorname{Weight}(C)$. Thus the weight of $C$ is always in $\mathbb{Z}^D$.

Observe that our definition of weight extends in a natural way to the oriented paths, open walks, and closed walks of $L$ on $A$.

Just as in Section \ref{rigid}, the sum of the realisations of the moves of $L$ along $C$ telescopes, and so must equal the zero vector. In other words, we must have that \[\omega_1v_1 + \omega_2v_2 + \cdots + \omega_Dv_D = \mathbf{0}.\]

In terms of the $x_i$ and $y_i$, this works out to \begin{align*} \omega_1x_1 + \omega_2x_2 + \cdots + \omega_Dx_D &= 0 \text{ and }\\ \omega_1y_1 + \omega_2y_2 + \cdots + \omega_Dy_D &= 0. \end{align*}

Thus every oriented cycle of $L$ on $A$ contributes two homogeneous linear constraints which the $x_i$ and $y_i$ must satisfy. We call all such constraints over all oriented cycles of $L$ on $A$ the \emph{cyclic constraints} on $v$. Note that, since every finite graph contains only finitely many oriented cycles, there are only finitely many cyclic constraints on $v$ as well.

The following lemma shows that the constraints on $v$ which we have figured out thus far are in fact all that we need to know about $v$.

\medbreak

\begin{lemma} Let $v$ be any vector in $\mathbb{R}^{2D}$. Then $v$ is the summary of some infinitesimal motion $I_v$ of $\mathcal{F}_A$ if and only if $v$ satisfies all directional and all cyclic constraints. \label{inf-flexibility-summary} \end{lemma}

\medbreak

\begin{proof} We know necessity already, and so we turn directly to sufficiency. Our proof essentially follows along the lines of our earlier observation that, since $\mathcal{L}_A$ is connected by assumption, $v_I$ determines $I$ up to translation.

Suppose that $v$ does satisfy all directional and all cyclic constraints, and define the mapping $I_v$ from the cells of $A$ to $\mathbb{R}^2$ as follows.

Choose some cell $r$ of $A$. For each cell $a$ of $A$, choose some path $P_a$ of $L$ on $A$ which leads from $r$ to $a$, and let $\operatorname{Weight}(P_a) = (\omega_{a, 1}, \omega_{a, 2}, \ldots, \omega_{a, D})$. Then define \[I_v(a) = \omega_{a, 1}v_1 + \omega_{a, 2}v_2 + \cdots + \omega_{a, D}v_D.\]

Let $e = a \to b$ be any move of $L$ on $A$ of direction $d_i$ in rhombic class $R_i$. We claim that $I_v(b) - I_v(a) = v_i$.

To see this, consider the oriented closed walk $W$ of $L$ on $A$ which starts from $r$, traces path $P_a$ until it reaches cell $a$, traverses edge $ab$ of $L$ from $a$ to $b$, and then returns back to $r$ by tracing path $P_b$ in reverse.

Observe that the edge multiset of every oriented closed walk of $L$ on $A$ is the multiset union of the edge multisets of several oriented cycles of $L$ on $A$ together with several pairs of mutually annihilating moves of $L$ on $A$. Therefore, the weight of $W$, say, $(w_1, w_2, \ldots, w_D)$, is an integer linear combination of the weights of several oriented cycles of $L$ on $A$. Since $v$ satisfies all cyclic constraints, from this it follows that \[w_1v_1 + w_2v_2 + \cdots + w_Dv_D = \mathbf{0}.\]

On the other hand, \[\operatorname{Weight}(W) = \operatorname{Weight}(P_a) + (0, 0, \ldots, 0, 1, 0, 0, \ldots, 0) - \operatorname{Weight}(P_b).\]

Consequently, $w_j = \omega_{a, j} - \omega_{b, j}$ for all $j \neq i$, and $w_i = \omega_{a, i} - \omega_{b, i} + 1$. When we substitute these values for the $w_j$ into the identity $w_1v_1 + w_2v_2 + \cdots + w_Dv_D = \mathbf{0}$ and rearrange appropriately, we obtain $I_v(b) - I_v(a) = v_i$, as required.

Since $v$ satisfies all directional constraints, it follows that vectors $a - b$ and $I_v(a) - I_v(b)$ are orthogonal for all edges $ab$ of $L$ on $A$. Therefore, $I_v$ is an infinitesimal motion of $\mathcal{F}_A$; and it is clear that $v$ is the summary of $I_v$. \end{proof}

\medbreak

We proceed to ``package'' each set of constraints on $v$ into a more convenient form.

Given any constraint on $v$ of the form $\alpha_1 x_1 + \beta_1 y_1 + \alpha_2 x_2 + \beta_2 y_2 + \cdots + \alpha_D x_D + \beta_D y_D = 0$, we refer to the vector $(\alpha_1, \beta_1, \alpha_2, \beta_2, \ldots, \alpha_D, \beta_D)$ of its coefficients as a \emph{constraint vector}.

Let $\mathcal{H}^\text{Dir}_A$ be the vector space spanned by all directional constraint vectors and let $\mathcal{H}^\text{Cyc}_A$ be the vector space spanned by all cyclic constraint vectors. Then both of $\mathcal{H}^\text{Dir}_A$ and $\mathcal{H}^\text{Cyc}_A$ are subspaces of $\mathbb{R}^{2D}$. We call $\mathcal{H}^\text{Dir}_A$ the \emph{directional constraint space} of $L$ on $A$ and $\mathcal{H}^\text{Cyc}_A$ the \emph{cyclic constraint space} of $L$ on $A$.

Moreover, let $\mathcal{H}^\text{Const}_A$ be the sum of $\mathcal{H}^\text{Dir}_A$ and $\mathcal{H}^\text{Cyc}_A$. Clearly, $\mathcal{H}^\text{Const}_A$ is also a subspace of $\mathbb{R}^{2D}$. Naturally enough, we call $\mathcal{H}^\text{Const}_A$ the \emph{constraint space} of $L$ on $A$.

Then $v$ satisfies all directional constraints if and only if it is orthogonal to $\mathcal{H}^\text{Dir}_A$, and $v$ satisfies all cyclic constraints if and only if it is orthogonal to $\mathcal{H}^\text{Cyc}_A$. By Lemma \ref{inf-flexibility-summary}, it follows that $v$ is the summary of some infinitesimal motion $I_v$ of $\mathcal{F}_A$ if and only if $v$ is orthogonal to $\mathcal{H}^\text{Const}_A$.

Consequently, \[\dim \operatorname{Inf}'(\mathcal{F}_A) = 2D - \dim \mathcal{H}^\text{Const}_A\] and \[\dim \operatorname{Inf}(\mathcal{F}_A) = 2D + 2 - \dim \mathcal{H}^\text{Const}_A.\]

Hence, the following criteria for infinitesimal flexibility.

\medbreak

\begin{theorem} The leaper framework of $L$ on $A$ is infinitesimally flexible if and only if the dimension of $\mathcal{H}^\text{Const}_A$ is at most $2D - 2$. \label{inf-flexibility-method} \end{theorem}

\medbreak

\begin{corollary} If the dimension of $\mathcal{H}^\text{Cyc}_A$ does not exceed $D - 2$, then the leaper framework of $L$ on $A$ is infinitesimally flexible. \label{inf-flexibility-sufficient} \end{corollary}

\medbreak

\begin{proof} By $\dim \mathcal{H}^\text{Dir}_A = D$, $\dim \mathcal{H}^\text{Const}_A \le \dim \mathcal{H}^\text{Dir}_A + \dim \mathcal{H}^\text{Cyc}_A$, and Theorem \ref{inf-flexibility-method}. \end{proof}

\medbreak

Corollary \ref{inf-flexibility-sufficient} is less general than Theorem \ref{inf-flexibility-method}, but its conditions are easier to verify in the cases where it does apply. Soon we will see that the dimension of $\mathcal{H}^\text{Cyc}_A$ can be bounded from above without too much trouble by combinatorial means; while in Sections \ref{flex-i} and \ref{flex-ii} we are going to convince ourselves that calculating the dimension of $\mathcal{H}^\text{Const}_A$ precisely is somewhat trickier.

We bound the dimension of $\mathcal{H}^\text{Cyc}_A$ from above as follows.

Let $G$ be any subgraph of $\mathcal{L}_A$. Then we write $\mathcal{H}^\text{Cyc}_A \restriction G$ for the vector space spanned by all cyclic constraint vectors obtained from oriented cycles in $G$. Of course, $\mathcal{H}^\text{Cyc}_A \restriction G$ is a subspace of $\mathcal{H}^\text{Cyc}_A$ for all $G$.

Our plan will be reduce the problem of estimating the dimension of $\mathcal{H}^\text{Cyc}_A$ to the problem of estimating the dimension of $\mathcal{H}^\text{Cyc}_A \restriction G$ for successively simpler and simpler subgraphs $G$. We simplify subgraphs by deleting edges from them, and the following lemma allows us to keep track of dimension as we carry out the deletions.

\medbreak

\begin{lemma} Let $e$ be any edge of $G$. If $e$ is not part of any cycle in $G$ or if $e$ is part of some oriented cycle in $G$ of zero weight, then \[\dim \mathcal{H}^\text{Cyc}_A \restriction G = \dim \mathcal{H}^\text{Cyc}_A \restriction (G \setminus e).\]

Moreover, \[\dim \mathcal{H}^\text{Cyc}_A \restriction G \le 2 + \dim \mathcal{H}^\text{Cyc}_A \restriction (G \setminus e)\] for all edges $e$ of $G$. \label{edge-deletion} \end{lemma}

\medbreak

\begin{proof} The case when $e$ is not part of any cycle in $G$ is clear. Suppose, then, that $e$ is part of the oriented cycle $C$ in $G$.

Let $C'$ be any oriented cycle in $G$. It suffices to show that $\operatorname{Weight}(C')$ is in the linear hull of $\operatorname{Weight}(C)$ and the weights of all oriented cycles in $G \setminus e$. Then it will follow that $\mathcal{H}^\text{Cyc}_A \restriction G$ is a subspace of the sum of the vector space spanned by the two cyclic constraint vectors obtained from $C$ and the vector space $\mathcal{H}^\text{Cyc}_A \restriction (G \setminus e)$.

(In particular, when $C$ is of zero weight, we obtain that $\mathcal{H}^\text{Cyc}_A \restriction G$ and $\mathcal{H}^\text{Cyc}_A \restriction (G \setminus e)$ actually coincide.)

If $e$ is not part of $C'$, then there is nothing to prove. Otherwise, suppose, without loss of generality, that $e = ab$ and that both of $C$ and $C'$ traverse $e$ from $a$ to $b$. Consider the closed walk $W$ which starts from $b$, traces $C$ until it reaches $a$, and then traces $C'$ in reverse until it returns back to $b$.

Then $\operatorname{Weight}(W) = \operatorname{Weight}(C) - \operatorname{Weight}(C')$. On the other hand, just as in the proof of Lemma \ref{inf-flexibility-summary}, the weight of $W$ is an integer linear combination of the weights of several oriented cycles in $G \setminus e$. This completes the proof. \end{proof}

\medbreak

Lemma \ref{edge-deletion} suggests the following definition. We say that edge $e$ of subgraph $G$ of $\mathcal{L}_A$ is \emph{superfluous} in $G$ if $e$ is part of some oriented cycle in $G$ of zero weight. By Lemma \ref{edge-deletion}, removing superfluous edges does not affect the dimension of the cyclic constraint space corresponding to the subgraph in any way.

Thus we proceed as follows. Starting from $\mathcal{L}_A$, we remove superfluous edges, one by one, for as long as we can manage. Then, in the end, we apply the second part of Lemma \ref{edge-deletion} to the resulting subgraph in order to obtain our upper bound for the dimension of $\mathcal{H}^\text{Cyc}_A$. The following lemma encapsulates this procedure.

\medbreak

\begin{lemma} Suppose that some sequence of superfluous edge deletions starts from $\mathcal{L}_A$ and leaves a total of $|A| - 1 + n$ edges on $A$. Then $\dim \mathcal{H}^\text{Cyc}_A \le 2n$. \label{inf-flexibility-procedure} \end{lemma}

\medbreak

\begin{proof} Let $H$ be the subgraph of $\mathcal{L}_A$ which we obtain in the end. By Lemma \ref{edge-deletion}, $\dim \mathcal{H}^\text{Cyc}_A = \dim \mathcal{H}^\text{Cyc}_A \restriction H$.

Since every superfluous edge is part of some cycle before we delete it, the deletion of superfluous edges does not affect connectedness. Therefore, $H$ is connected.

Let $T$ be any spanning tree of $H$. Since the vertices of $T$ are exactly the cells of $A$, $T$ contains a total of $|A| - 1$ edges.

Delete all edges of $H$ outside of $T$. By Lemma \ref{edge-deletion}, each such deletion decreases the dimension of the cyclic constraint space corresponding to the subgraph by at most two. On the other hand, since $T$ does not contain any cycles, the dimension of the cyclic constraint space corresponding to $T$ is zero. \end{proof}

\medbreak

In effect, Lemma \ref{inf-flexibility-procedure} reduces the linear-algebraic problem of bounding the dimension of $\mathcal{H}^\text{Cyc}_A$ from above to the purely combinatorial problem of finding a very thorough sequence of superfluous edge deletions.

We take a moment here for a brief discussion of one theoretical question related to our superfluous edge deletion procedure, and then we go on to grid $\Psi$.

One might naturally wonder if our procedure is ``stable'', in the following sense. Can we always delete several superfluous edges so as to obtain a sharp upper bound for the dimension of $\mathcal{H}^\text{Cyc}_A$? Moreover, suppose that we make the first few superfluous edge deletions ``at random''. Can we still always delete a few more superfluous edges so that we obtain a sharp upper bound in the end, or could it happen that we have painted ourselves in a corner?

Note that Lemma \ref{inf-flexibility-procedure} always yields even upper bounds. Thus let us, first of all, verify that the dimension of $\mathcal{H}^\text{Cyc}_A$ is even.

Let $\mathcal{H}^\text{Weight}_A$ be the vector space spanned by the weights of all oriented cycles in $\mathcal{L}_A$, and let $E$ be the dimension of $\mathcal{H}^\text{Weight}_A$. We claim that the dimension of $\mathcal{H}^\text{Cyc}_A$ is then $2E$. Indeed, $\mathcal{H}^\text{Cyc}_A$ is the direct sum of the vector space spanned by all cyclic constraint vectors which only address the $x_i$ and the vector space spanned by all cyclic constraint vectors which only address the $y_i$. Clearly, each one of these two vector spaces is isomorphic to $\mathcal{H}^\text{Weight}_A$.

With that taken care of, there is one generalisation of Lemma \ref{edge-deletion} which makes our procedure stable in the sense outlined above.

Let $e$ be any edge of $L$ on $A$, and let $W$ be any oriented closed walk of $L$ on $A$. Then we say that $e$ is \emph{essential} for $W$ if $W$ traverses $e$ a different number of times in both directions. (Intuitively, if $e$ contributes nontrivially towards the weight of $W$.) Moreover, we say that edge $e$ of subgraph $G$ is \emph{essentially superfluous} in $G$ if $e$ is essential for some oriented closed walk in $G$ of zero weight.

Clearly, the first part of Lemma \ref{edge-deletion} generalises to essentially superfluous edges. Moreover, there is the following lemma.

\medbreak

\begin{lemma} Suppose that $G$ is connected and that it contains at least $|A| + E$ edges. Then $G$ contains at least one essentially superfluous edge. \label{stable} \end{lemma}

\medbreak

Note that, similarly to superfluous edge deletion, the deletion of essentially superfluous edges does not affect connectedness.

\medbreak

\begin{proof} Let $T$ be any spanning tree of $G$. Choose $E + 1$ distinct edges of $G$ outside of $T$, say, $e_1$, $e_2$, \ldots, $e_{E + 1}$. For all $i$, let $C_i$ be the unique cycle in subgraph $T \cup e_i$ of $\mathcal{L}_A$.

Since $\dim \mathcal{H}^\text{Weight}_A < E + 1$, vectors $\operatorname{Weight}(C_i)$ are linearly dependent. Since all of them are in $\mathbb{Z}^D$, without loss of generality the coefficients of their linear dependence are integers. Let, then, \[\alpha_1\operatorname{Weight}(C_1) + \alpha_2\operatorname{Weight}(C_2) + \cdots + \alpha_{E + 1}\operatorname{Weight}(C_{E + 1}) = \mathbf{0},\] where $\alpha_i$ is an integer for all $i$ and, without loss of generality, $\alpha_1$ is nonzero.

For each $i$, take $|\alpha_i|$ copies of $C_i$. If $\alpha_i > 0$, then leave each copy as is; otherwise, if $\alpha_i < 0$, reverse the directions of all moves within each copy.

We splice all such copies over all $i$ together into a single oriented closed walk $W$ in $G$ by iterating the following operation. Given two oriented closed walks $W'$ and $W''$ in $G$ such that $W'$ visits some cell $a'$ and $W''$ visits some cell $a''$, we splice $W'$ and $W''$ together into a single oriented closed walk in $G$ by tracing $W'$ in its entirety from $a'$ to $a'$, travelling along any path $P$ in $G$ from $a'$ to $a''$, tracing $W''$ in its entirety from $a''$ to $a''$, and, finally, travelling along $P$ in reverse from $a''$ until we return back to $a'$.

Then $e_1$ is essential for $W$ and $W$ is a witness for $e_1$ being essentially superfluous in $G$. \end{proof}

\medbreak

Note that Lemma \ref{stable} has nothing at all to say about the difficulty of proving that any edge in $G$ is essentially superfluous when we do not know $E$ in advance.

We are not going to need the more general version of our procedure for grid $\Psi$. In fact, as we are about to see, this generalisation is way overpowered for our purposes; but it is also a reassuring indication that we are on the right track.

At long last, we proceed to grid $\Psi$.

By Lemma \ref{rhombic-psi}, Lemma \ref{inf-flexibility-procedure}, and Corollary \ref{inf-flexibility-sufficient}, it suffices to delete several superfluous edges, starting from $\mathcal{L}_\Psi$, so that at most $|\Psi| + 2$ edges remain on $\Psi$.

Throughout the proof, we are only ever going to apply Lemma \ref{edge-deletion} to rhombuses; that is, to oriented cycles of length four. Thus even the original Lemma \ref{edge-deletion} is much more general than grid $\Psi$ requires.

There is one manoeuvre which we will apply multiple times over the course of the deletion process, as follows.

Let $G$ be any subgraph of $\mathcal{L}_\Psi$. We say that two edges $e'$ and $e''$ of $G$ are in the same \emph{rhombic class relative to $G$} if there exists some sequence of edges $e_1 = e'$, $e_2$, \ldots, $e_n = e''$ of $L$ on $\Psi$ such that, for all $i$, $e_i$ and $e_{i + 1}$ form a rhombus all four of whose sides are edges of $G$. This is a natural extension of the definition of a rhombic class on a grid.

Let $R$ be any rhombic class relative to $G$. We claim that there exists some sequence of superfluous edge deletions, starting from $G$ and with rhombuses as witnesses, which eliminates all but one edges of $R$.

To see this, define graph $\mathcal{R}$ as follows: The vertices of $\mathcal{R}$ are all edges of $R$, and two edges of $R$ are joined by an edge in $\mathcal{R}$ if and only if they form a rhombus all four of whose sides are edges of $G$. Since $R$ is a rhombic class relative to $G$, $\mathcal{R}$ is connected. Choose any rooted spanning tree $T$ of $\mathcal{R}$, and then delete the edges of $R$ in descending order of depth within $T$.

This completes our description of that useful manoeuvre, and we get started deleting superfluous edges.

We proceed slope by slope, and we begin with slope $\frac{q}{p}$.

We leave one edge of rhombic class $R_1$ and one edge of rhombic class $R_2$, and we delete all other edges of these rhombic classes. Let $G'$ be the subgraph of $\mathcal{L}_\Psi$ which we obtain in this way.

We go on to slope $-\frac{q}{p}$.

Let $\Psi_\texttt{LR}$ be the lower right subgrid of $\Psi$ of height $h_\Psi - q$ and width $w_\Psi - p$. Observe that the two edges of slope $\frac{q}{p}$ still on the grid do not participate in any rhombuses. On the other hand, $G'$ still contains all edges of the other three slopes. Therefore, analogously to the proof of Lemma \ref{rhombic-class}, two moves $a' \to b'$ and $a'' \to b''$ of direction $(-p, q)$ in $G'$ are in the same rhombic class relative to $G'$ if and only if $a'$ and $a''$ are in the same connected component of graph $\mathcal{L}_{\Psi_\texttt{LR}} \restriction \{\frac{p}{q}, -\frac{p}{q}\}$.

By Lemma \ref{nets-connectedness}, it follows that there are exactly $2pq$ rhombic classes of slope $-\frac{q}{p}$ relative to $G'$. We leave one edge out of each such class, and we delete all other edges of this slope. Let $G''$ be the subgraph of $\mathcal{L}_\Psi$ which we obtain in this way.

We continue to slope $\frac{p}{q}$, and we handle it analogously to slope $-\frac{q}{p}$.

Let $\Psi_\texttt{LL}$ be the lower left subgrid of $\Psi$ of height $h_\Psi - p$ and width $w_\Psi - q$. Observe that there are no edges of slopes $\frac{q}{p}$ and $-\frac{q}{p}$ still on the grid which participate in any rhombuses at all. On the other hand, $G''$ still contains all edges of the other two slopes, $\frac{p}{q}$ and $-\frac{p}{q}$. Therefore, analogously to the proof of Lemma \ref{rhombic-class}, two moves $a' \to b'$ and $a'' \to b''$ of direction $(q, p)$ in $G''$ are in the same rhombic class relative to $G''$ if and only if $a'$ and $a''$ are in the same connected component of graph $\mathcal{L}_{\Psi_\texttt{LL}} \restriction -\frac{p}{q}$.

By Lemma \ref{single-slope}, it follows that there are exactly $p(w_\Psi - q) + q(h_\Psi - p) - pq$ rhombic classes of slope $\frac{p}{q}$ relative to $G''$. As before, we leave one edge out of each such class, and we delete all other edges of this slope.

There are no more rhombuses left on the grid, and so we do not delete any edges of the fourth and final slope, $-\frac{p}{q}$. There are a total of $(h_\Psi - p)(w_\Psi - q)$ edges of $L$ on $\Psi$ of this slope.

Let us take stock. The total number of edges of $L$ which remain on grid $\Psi$ by this point is \begin{align*} &\phantom{= \mbox{ }} 2 + 2pq + [p(w_\Psi - q) + q(h_\Psi - p) - pq] + (h_\Psi - p)(w_\Psi - q)\\ &= 2 + 2pq + [(h_\Psi - p) + p][(w_{\Psi} - q) + q] - 2pq\\ &= 2 + h_\Psi w_\Psi\\ &= |\Psi| + 2. \end{align*}

This completes our proof of Theorem \ref{inf-flexibility-psi}, and we go on to continuous flexibility.

\section{A General Method} \label{method}

In this section, we develop one general method for establishing the flexibility of leaper frameworks. Theorem \ref{flexibility-method} below does most of the heavy lifting, and its proof occupies the bulk of our considerations. Then we give a step-by-step outline of our method, together with some comments on how we apply it in practice.

Let $A$ be any grid such that $\mathcal{L}_A$ is connected. All definitions of Section \ref{inf-flex} carry over in this section; in particular, so do the definition of the rhombic classes of $L$ on $A$, namely $R_1$, $R_2$, \ldots, $R_D$, and the definitions of the three constraint spaces of $L$ on $A$, namely $\mathcal{H}^\text{Dir}_A$, $\mathcal{H}^\text{Cyc}_A$, and $\mathcal{H}^\text{Const}_A$.

To begin with, we review a couple of basic notions from linear algebra.

Let $V$ be some vector space and let $V'$ and $V''$ be two subspaces of $V$. We say that $V'$ and $V''$ are in \emph{general position} if \[\dim (V' + V'') = \min \{\dim V' + \dim V'', \dim V\}.\]

This is a special case of the definition of general position for affine subspaces in Euclidean geometry. In both settings, we aim to capture the intuitive notion of two objects whose positioning relative to each other is as arbitrary as possible, so that they are not ``aligned'' in any special way.

Moreover, we say that $V'$ and $V''$ are \emph{essentially disjoint} if their intersection consists only of the zero vector. When $\dim V' + \dim V'' \le \dim V$, $V'$ and $V''$ are essentially disjoint if and only if they are in general position.

In this language, our criterion for flexibility sounds as follows.

\medbreak

\begin{theorem} Suppose that the dimension of $\mathcal{H}^\text{Const}_A$ is at most $2D - 2$. Furthermore, suppose that $\mathcal{H}^\text{Dir}_A$ and $\mathcal{H}^\text{Cyc}_A$ are in general position; or, equivalently, that they are essentially disjoint. Then the leaper framework of $L$ on $A$ is flexible. \label{flexibility-method} \end{theorem}

\medbreak

By Theorem \ref{inf-flexibility-method}, the condition that $\dim \mathcal{H}^\text{Const}_A \le 2D - 2$ is equivalent to $\mathcal{F}_A$ being infinitesimally flexible, and is therefore necessary.

Our proof of Theorem \ref{flexibility-method} follows more or less the same overall plan as our proof of Theorem \ref{inf-flexibility-method}. However, in the setting of continuous flexibility, this plan becomes a lot more complicated to implement.

Let $\mathbf{c}$ be some placement of $\mathcal{L}_A$. For simplicity, this time around we do not require that the points of $\mathbf{c}$ are pairwise distinct. We say that $\mathbf{c}$ is \emph{proper} if it satisfies the following conditions:

(a) If two moves of $L$ on $A$ form a rhombus, then they are realised by equal vectors under $\mathbf{c}$.

(b) Placement $\mathbf{c}$ of $\mathcal{L}_A$ is equivalent to the canonical placement of $\mathcal{L}_A$.

(Intuitively, condition (a) says that the embedding of $\mathcal{L}_A$ into the plane defined by $\mathbf{c}$ preserves rhombuses. Though note that, while we forbid degeneracy for combinatorial rhombuses, in this instance we permit it for geometrical rhombuses.)

In particular, if all points of $\mathbf{c}$ are pairwise distinct and $\mathbf{c}$ is equivalent to the canonical placement of $\mathcal{L}_A$, then $\mathbf{c}$ is proper. The converse is false, and so propriety is a strictly weaker requirement on the placements of $\mathcal{L}_A$ than what we had in Section \ref{rigid}.

Echoing Section \ref{inf-flex}, let $\mathbf{c}$ be some proper placement of $\mathcal{L}_A$ and, for all $i$, let all moves of direction $d_i$ in rhombic class $R_i$ be realised by vector $v_i = (x_i, y_i)$ under $\mathbf{c}$. Of course, then all moves of direction $-d_i$ in rhombic class $R_i$ are realised by vector $-v_i$ under $\mathbf{c}$. We define vector $v_\mathbf{c}$ by $v_\mathbf{c} = (x_1, y_1, x_2, y_2, \ldots,\allowbreak x_D, y_D)$, and we call this vector the \emph{summary} of $\mathbf{c}$. Since $\mathcal{L}_A$ is connected by assumption, $v_\mathbf{c}$ determines $\mathbf{c}$ essentially uniquely, up to translation.

Still following along the lines of Section \ref{inf-flex}, let $v$ be any vector in $\mathbb{R}^{2D}$. When is $v$ the summary of some proper placement $\mathbf{c}_v$ of $\mathcal{L}_A$?

Write, as above, $v = (x_1, y_1, x_2, y_2, \ldots, x_D, y_D)$, and set $v_i = (x_i, y_i)$ for all $i$. In every summary of a proper placement of $\mathcal{L}_A$, the constraints $|v_i| = \sqrt{p^2 + q^2}$ must hold for all $i$. We call these the \emph{bar-length constraints} on $v$. The \emph{cyclic constraints} on $v$ we carry over from Section \ref{inf-flex} unchanged. Every summary of a proper placement of $\mathcal{L}_A$ must satisfy all of them as well, by summation along each subembedding of an oriented cycle of $\mathcal{L}_A$ into the plane. (As in Section \ref{rigid}, where we applied this argument to $C_\Phi$, $C'_\Phi$, and $C_\Theta$.)

Here follows the analogue of Lemma \ref{inf-flexibility-summary} for global and continuous flexibility.

\medbreak

\begin{lemma} Let $v$ be any vector in $\mathbb{R}^{2D}$. Then $v$ is the summary of some proper placement $\mathbf{c}_v$ of $\mathcal{L}_A$ if and only if $v$ satisfies all bar-length and all cyclic constraints. \nolinebreak \label{flexibility-summary} \end{lemma}

\medbreak

\begin{proof} Fully analogous to the proof of Lemma \ref{inf-flexibility-summary}. \end{proof}

\medbreak

Note that the ``if'' part of Lemma \ref{flexibility-summary} does not hold anymore when condition (a) of the definition of a proper placement is replaced (as it was in Section \ref{rigid}) with the condition that all points of $\mathbf{c}$ are pairwise distinct.

We denote the summary of the canonical placement of $\mathcal{L}_A$, corresponding to framework $\mathcal{F}_A$, by $v_\text{Canon}$. It is given by $v_i = d_i$ for all $i$. Equivalently, let $d_i = (d'_i, d''_i)$ for all $i$; then $v_\text{Canon} = (d'_1, d''_1, d'_2, d''_2, \ldots, d'_D, d''_D)$. Of course, $v_\text{Canon}$ satisfies all bar-length and all cyclic constraints.

Similarly to what we did in Section \ref{inf-flex}, we are going to work with summaries rather than with placements directly. Our plan for the proof of Theorem \ref{flexibility-method} will be to show that we can deform $v_\text{Canon}$ continuously so that it continues to satisfy all bar-length and all cyclic constraints. Once we have done that, obtaining a flexion of $\mathcal{F}_A$ will require very little additional effort.

We begin by encoding all bar-length and all cyclic constraints into a single vector-valued function.

For $1 \le i \le D$, define the function $f_i : \mathbb{R}^{2D} \to \mathbb{R}$ by \[f_i(x_1, y_1, x_2, y_2, \ldots, x_D, y_D) = x_i^2 + y_i^2 - p^2 - q^2.\]

Let vectors $\lambda_1$, $\lambda_2$, \ldots, $\lambda_{2E}$ form a basis of $\mathcal{H}^\text{Cyc}_A$. For $1 \le j \le 2E$, let \[\lambda_j = (\alpha_{j, 1}, \beta_{j, 1}, \alpha_{j, 2}, \beta_{j, 2}, \ldots, \alpha_{j, D}, \beta_{j, D}),\] and then define the function $g_j : \mathbb{R}^{2D} \to \mathbb{R}$ by \[g_j(x_1, y_1, x_2, y_2, \ldots, x_D, y_D) = \alpha_{j, 1}x_1 + \beta_{j, 1}y_1 + \alpha_{j, 2}x_2 + \beta_{j, 2}y_2 + \cdots + \alpha_{j, D}x_D + \beta_{j, D}y_D.\]

Lastly, define the vector-valued function $\mathbf{f} : \mathbb{R}^{2D} \to \mathbb{R}^{D + 2E}$ by \[\mathbf{f} = (f_1, f_2, \ldots, f_D, g_1, g_2, \ldots, g_{2E}).\]

Then, by Lemma \ref{flexibility-summary}, $v$ is the summary of some proper placement $\mathbf{c}_v$ of $\mathcal{L}_A$ if and only if $v$ is a root of $\mathbf{f}$. In particular, $v_\text{Canon}$ is a root of $\mathbf{f}$.

Since each component function of $\mathbf{f}$ is a polynomial, $\mathbf{f}$ is everywhere continuously differentiable. Let $J$ be the Jacobian matrix of $\mathbf{f}$ at point $v_\text{Canon}$. We proceed to calculate $J$ explicitly.

For $1 \le i \le D$, row $i$ of $J$ has $2d'_i$ at position $2i - 1$, $2d''_i$ at position $2i$, and zeros everywhere else. In other words, row $i$ of $J$ equals twice the vector of the coefficients of the directional constraint on $v$ which addresses $v_i$ and $d_i$. Consequently, the first $D$ rows of $J$ form a basis of $\mathcal{H}^\text{Dir}_A$.

For $1 \le j \le 2E$, row $D + j$ of $J$ is simply $\lambda_j$. Consequently, the remaining $2E$ rows of $J$ form a basis of $\mathcal{H}^\text{Cyc}_A$.

Since $\mathcal{H}^\text{Dir}_A$ and $\mathcal{H}^\text{Cyc}_A$ are essentially disjoint by assumption, the rows of $J$ are linearly independent. Therefore, $J$ has full rank and its rows form a basis of $\mathcal{H}^\text{Const}_A$.

Suppose that $J$ has $\mu = D - 2E$ fewer rows than columns. Since the dimension of $\mathcal{H}^\text{Const}_A$ is at most $2D - 2$ by assumption, $\mu \ge 2$.

Since $J$ has full rank, it is possible to delete some $\mu$ columns from $J$ so as to obtain a square matrix $\widehat{J}$ which has full rank as well. Out of all component variables of $v$, we call the $\mu$ variables which correspond to these $\mu$ columns of $J$ \emph{free}, and the remaining $2D - \mu = D + 2E$ component variables of $v$ we call \emph{bound}.

Note that $x_i$ and $y_i$ cannot both be free for any $i$, since then row $i$ of $\widehat{J}$ would equal the zero vector, and so it would be impossible for $\widehat{J}$ to have full rank. Thus suppose, without loss of generality, that the free variables are precisely $x_1$, $x_2$, \ldots, $x_{\mu'}$ and $y_{\mu'' + 1}$, $y_{\mu'' + 2}$, \ldots, $y_D$, where $\mu' + (D - \mu'') = \mu$. For convenience, we write $\mathbf{x}$ for $x_1$, $x_2$, \ldots, $x_{\mu'}$ and $\mathbf{y}$ for $y_{\mu'' + 1}$, $y_{\mu'' + 2}$, \ldots, $y_D$. We also write $\mathbf{x}_\text{Canon}$ for $d'_1$, $d'_2$, \ldots, $d'_{\mu'}$ and $\mathbf{y}_\text{Canon}$ for $d''_{\mu'' + 1}$, $d''_{\mu'' + 2}$, \ldots, $d''_D$, the values of $\mathbf{x}$ and $\mathbf{y}$ in $v_\text{Canon}$.

The point of this classification of all component variables of $v$ into free and bound ones is as follows. Roughly speaking, we are going to show that we can choose the values of all free variables somewhat arbitrarily, and that for each such choice there is a unique way to fill in the values of all bound variables so as to obtain a valid summary of a proper placement of $\mathcal{L}_A$. Then it will follow that the summaries of the proper placements of $\mathcal{L}_A$ vary with $\mu$ degrees of freedom, one per free component variable of $v$.

We achieve this feat by the Implicit Function Theorem. This step is the core of the proof.

The Implicit Function Theorem applies to $\mathbf{f}$ since $v_\text{Canon}$ is a root of $\mathbf{f}$, $\mathbf{f}$ is everywhere continuously differentiable, and $\widehat{J}$ is invertible. Thus we obtain that there is some open neighbourhood $U$ of $(\mathbf{x}_\text{Canon}, \mathbf{y}_\text{Canon})$ in $\mathbb{R}^\mu$ such that there exists a unique vector-valued function $\mathbf{F} : U \to \mathbb{R}^{D + 2E}$ which satisfies certain conditions, to be given shortly.

Let $\mathbf{F} = (F_{\mu' + 1}, F_{\mu' + 2}, \ldots, F_D, G_1, G_2, \ldots, G_{\mu''})$. For convenience, we write $v(\mathbf{x}, \mathbf{y})$ for the vector $v$ in $\mathbb{R}^{2D}$ whose bound variables have been set to $x_i = F_i(\mathbf{x}, \mathbf{y})$ and $y_j = G_j(\mathbf{x}, \mathbf{y})$ for all admissible $i$ and $j$. Then $\mathbf{F}$ is the unique vector-valued function from $U$ to $\mathbb{R}^{D + 2E}$ such that:

(a) $v(\mathbf{x}_\text{Canon}, \mathbf{y}_\text{Canon}) = v_\text{Canon}$, and

(b) For all $(\mathbf{x}, \mathbf{y})$ in $U$, vector $v(\mathbf{x}, \mathbf{y})$ is a root of $\mathbf{f}$, and thus also the summary of some proper placement $\mathbf{c}_{v(\mathbf{x}, \mathbf{y})}$ of $\mathcal{L}_A$.

Moreover, by the Implicit Function Theorem, $\mathbf{F}$ is continuous and continuously differentiable in $U$.

Thus $\mathbf{F}$ does precisely what we needed to be done: When we hand over to it the values of all free variables, chosen arbitrarily within some reasonable open neighbourhood of their canonical values, $\mathbf{F}$ sets all bound variables so as to give us the summary of some proper placement of $\mathcal{L}_A$.

In broad strokes, the rest of the proof goes as follows. We have established that, in some open neighbourhood of $v_\text{Canon}$, the summaries of the proper placements of $\mathcal{L}_A$ vary with $\mu$ degrees of freedom. We must factor out all rotations; this shaves off one degree of freedom. However, since $\mu \ge 2$, we still have at least one degree of freedom left, and this is enough for $\mathcal{F}_A$ to flex continuously.

We make this sketch rigorous as follows.

Suppose, without loss of generality, that $\mu' \ge 1$, and so also that $x_1$ is a free variable.

Consider any placement $\mathbf{c}$ of $\mathcal{L}_A$. Observe that rotation does not affect propriety in any way. Moreover, if $\mathbf{c}$ is proper, then up to translation there is a unique rotation $\rho$ such that all moves of direction $d_1$ in rhombic class $R_1$ are realised by vector $d_1$ under placement $\rho(\mathbf{c})$ of $\mathcal{L}_A$.

Our approach to factoring out all rotations is rather unsophisticated. We simply fix the value of $x_1$.

Fix, then, $x_1$ to its canonical value of $d'_1$. Because of the bar-length constraint on $v_1$, this forces $y_1$ to equal either $d''_1$ or $-d''_1$. However, since $\mathbf{F}$ is continuous, $y_1$ cannot jump discontinuously between values. Therefore, this fixes $y_1$ to its canonical value of $d''_1$, too, and so it fixes $v_1$ to its canonical value of $d_1$ as well.

We write $\widehat{\mathbf{x}}$ for $x_2$, $x_3$, \ldots, $x_{\mu'}$ and $\widehat{\mathbf{x}}_\text{Canon}$ for $d'_2, d'_3, \ldots, d'_{\mu'}$. Moreover, let $\widehat{U}$ be some sufficiently small open neighbourhood of $(\widehat{\mathbf{x}}_\text{Canon}, \mathbf{y}_\text{Canon})$ in $\mathbb{R}^{\mu - 1}$ such that $(d'_1, \widehat{\mathbf{x}}, \mathbf{y})$ is in $U$ for all $(\widehat{\mathbf{x}}, \mathbf{y})$ in $\widehat{U}$.

Choose and fix some move $a \to b$ of $L$ on $A$ which points in direction $d_1$ and belongs to rhombic class $R_1$. For all $(\widehat{\mathbf{x}}, \mathbf{y})$ in $\widehat{U}$, we define $F(\widehat{\mathbf{x}}, \mathbf{y})$ to be the unique placement of $\mathcal{L}_A$ which maps cells $a$ and $b$ onto integer points $a$ and $b$, respectively, as in the canonical placement of $\mathcal{L}_A$, and whose summary is $\mathbf{F}(d'_1, \widehat{\mathbf{x}}, \mathbf{y})$. Then $F$ is a mapping from $\widehat{U}$ to $(\mathbb{R}^2)^{|A|}$.

We proceed to verify that $F$ is a flexion of $\mathcal{F}_A$. Note that $F$ does not quite fit the standard definition of a flexion which we gave in Section \ref{prelim}. However, it is clear how conditions (a), (b), and (c) of that definition generalise to flexions which accommodate multiple degrees of freedom, and it is just as clear how to derive a flexion of $\mathcal{F}_A$ in the narrow sense of Section \ref{prelim} from $F$.

To begin with, since $\mathbf{F}$ is continuous, so is $F$.

Of course, $F(\widehat{\mathbf{x}}_\text{Canon}, \mathbf{y}_\text{Canon})$ yields the canonical placement of $\mathcal{L}_A$.

Then, since placement $F(\widehat{\mathbf{x}}, \mathbf{y})$ of $\mathcal{L}_A$ is proper for all $(\widehat{\mathbf{x}}, \mathbf{y})$ in $\widehat{U}$, in particular it is equivalent to the canonical placement of $\mathcal{L}_A$.

Lastly, the realisations under $F(\widehat{\mathbf{x}}, \mathbf{y})$ of all moves in rhombic class $R_1$ are canonical. On the other hand, since $\mu \ge 2$, when $(\widehat{\mathbf{x}}, \mathbf{y}) \neq (\widehat{\mathbf{x}}_\text{Canon}, \mathbf{y}_\text{Canon})$ there is at least one $i$ such that free variable $x_i$ or $y_i$ is not set to its canonical value, and so the realisations under $F(\widehat{\mathbf{x}}, \mathbf{y})$ of all moves in rhombic class $R_i$ are noncanonical. Therefore, when $(\widehat{\mathbf{x}}, \mathbf{y}) \neq (\widehat{\mathbf{x}}_\text{Canon}, \mathbf{y}_\text{Canon})$, placement $F(\widehat{\mathbf{x}}, \mathbf{y})$ of $\mathcal{L}_A$ cannot be a translation or a rotation of the canonical placement of $\mathcal{L}_A$. And, provided that $\widehat{U}$ is sufficiently small, $F(\widehat{\mathbf{x}}, \mathbf{y})$ cannot be a reflection or a glide reflection of the canonical placement of $\mathcal{L}_A$, either.

This settles all three conditions (a), (b), and (c) of the definition of a flexion. Our proof of Theorem \ref{flexibility-method} is complete.

Concerning distinctness, note the following. All points in the canonical placement of $\mathcal{L}_A$ are pairwise distinct. Therefore, there is some open neighbourhood of $(\widehat{\mathbf{x}}_\text{Canon}, \mathbf{y}_\text{Canon})$ in $\mathbb{R}^{\mu - 1}$ such that, for all $(\widehat{\mathbf{x}}, \mathbf{y})$ in it, all points of $F(\widehat{\mathbf{x}}, \mathbf{y})$ are pairwise distinct as well.

Conversely, consider any placement $\mathbf{c}$ of $\mathcal{L}_A$ such that $\mathbf{c}$ is sufficiently close to the canonical placement of $\mathcal{L}_A$ as well as equivalent to it. Since $\mathbf{c}$ is sufficiently close to the canonical placement of $\mathcal{L}_A$, all of its points are pairwise distinct. Then, since $\mathbf{c}$ is equivalent to the canonical placement of $\mathcal{L}_A$, $\mathbf{c}$ is proper. By Lemma \ref{flexibility-summary}, it follows that the summary of $\mathbf{c}$ must necessarily be a root of $\mathbf{f}$. In other words, in our proof of Theorem \ref{flexibility-method} we have in fact described all such placements $\mathbf{c}$ of $\mathcal{L}_A$, even though we only set out to prove their existence. Hence, the following corollary.

\medbreak

\begin{corollary} Suppose that the dimension of $\mathcal{H}^\text{Const}_A$ is $2D - 1 - n$. Furthermore, suppose that $\mathcal{H}^\text{Dir}_A$ and $\mathcal{H}^\text{Cyc}_A$ are in general position; or, equivalently, that they are essentially disjoint. Then the leaper framework of $L$ on $A$ flexes with precisely $n$ degrees of freedom. \label{degrees-of-freedom} \end{corollary}

\medbreak

Observe that $n$ is always a nonnegative integer: Since $v_\text{Canon}$ is orthogonal to $\mathcal{H}^\text{Const}_A$ and $\mathcal{H}^\text{Const}_A$ is a subspace of $\mathbb{R}^{2D}$, the dimension of $\mathcal{H}^\text{Const}_A$ cannot exceed $2D - 1$. The special case when $n = 0$ also follows by Theorem \ref{inf-flexibility-method}.

The special case when $n = 1$ is particularly interesting. In it, $\mathcal{F}_A$ flexes with precisely one degree of freedom. Intuitively, this means that $\mathcal{F}_A$ flexes ``neatly''. Equivalently, it means that $\mathcal{F}_A$ admits an essentially unique flexion in the narrow sense of Section \ref{prelim}. (``Essentially unique'' here means that the flexion is unique up to translation and rotation, and in an appropriate open neighbourhood of the canonical placement of $\mathcal{L}_A$.) We encounter leaper frameworks of this kind in Sections \ref{flex-i} and \ref{flex-ii}.

There is one simple but important generalisation of our ideas, as follows.

Everything that we did in Sections \ref{rigid}, \ref{inf-flex}, and \ref{method} for rectangular grids continues to hold when $A$ is any set of cells such that $\mathcal{L}_A$ is connected. Indeed, we hardly ever really made much use of our assumption that $A$ is a rectangular grid.

The only exception is Lemma \ref{rhombic-class}, where we talk about some subgrids of $A$. In the more general setting where $A$ is an arbitrary set of cells, we need to amend this lemma accordingly.

Let $d = (d', d'')$ be some direction of the moves of $L$. We define $A \llbracket d \rrbracket$ to be the set of the initial cells of all moves of $L$ on $A$ which point in direction $d$; or, equivalently, $A \cap \{a - d \mid a \in A\}$. Then two moves $a' \to b'$ and $a'' \to b''$ of $L$ on $A$ in direction $d$ are in the same rhombic class of $L$ on $A$ if and only if cells $a'$ and $a''$ are in the same connected component of graph $\mathcal{L}_{A \llbracket d \rrbracket} \setminus \frac{d''}{d'}$.

Naturally, in the case when $A$ is a rectangular grid, we obtain the original Lemma \ref{rhombic-class}.

We are in fact going to need our techniques in this generality in Section \ref{flex-ii}, where we encounter leaper frameworks on somewhat stranger sets of cells.

This completes our brief detour into arbitrary sets of cells, and we go on to the step-by-step outline which we promised to give at the beginning of this section.

Theorems \ref{inf-flexibility-method} and \ref{flexibility-method}, taken together, form the engine which powers our method. Note that Theorem \ref{inf-flexibility-method} supplies a necessary and sufficient condition, whereas Theorem \ref{flexibility-method} supplies a sufficient condition only.

Let $A$ be any set of cells. We wish to determine if the leaper framework of $L$ on $A$ is flexible or rigid, and to this end we proceed as follows.

If $\mathcal{L}_A$ is not connected, then, as we remarked in the introduction, $\mathcal{F}_A$ is flexible, albeit trivially. Suppose, from this point on, that $\mathcal{L}_A$ is connected.

To begin with, we investigate the structure of the rhombic classes of $L$ on $A$. Lemma \ref{rhombic-class} and its generalisation above reduce this question to a study of the connected components of some forbidden-slope leaper graphs. In Sections \ref{rigid} and \ref{inf-flex} we did this for grids $\Phi$, $\Psi$, and $\Theta$.

Then we estimate the dimension of $\mathcal{H}^\text{Const}_A$. This mostly boils down to a study of $\mathcal{H}^\text{Weight}_A$, since we already know one nice basis of $\mathcal{H}^\text{Dir}_A$ and $\mathcal{H}^\text{Weight}_A$ uniquely determines $\mathcal{H}^\text{Cyc}_A$. In particular, every basis of $\mathcal{H}^\text{Weight}_A$ immediately yields a corresponding basis of $\mathcal{H}^\text{Cyc}_A$.

We estimate the dimensions of $\mathcal{H}^\text{Weight}_A$ and $\mathcal{H}^\text{Const}_A$ from below by constructing suitable oriented cycles of $L$ on $A$. This is exactly what we did for grids $\Phi$ and $\Theta$ in Section \ref{rigid}. On the other hand, we estimate the dimensions of $\mathcal{H}^\text{Weight}_A$ and $\mathcal{H}^\text{Const}_A$ from above by constructing suitable sequences of superfluous edge deletions. This is exactly what we did for grid $\Psi$ in Section \ref{inf-flex}.

If the dimension of $\mathcal{H}^\text{Const}_A$ is at least $2D - 1$, and so equals $2D - 1$, then by Theorem \ref{inf-flexibility-method} we obtain that $\mathcal{F}_A$ is rigid. This is the essence of our proofs of rigidity for grids $\Phi$ and $\Theta$ in Section \ref{rigid}. Otherwise, if the dimension of $\mathcal{H}^\text{Const}_A$ is at most $2D - 2$, then by Theorem \ref{inf-flexibility-method} we obtain that $\mathcal{F}_A$ is infinitesimally flexible, and the stage is set for an application of Theorem \ref{flexibility-method}.

To apply Theorem \ref{flexibility-method}, first we figure out sufficiently nice bases of $\mathcal{H}^\text{Weight}_A$ and $\mathcal{H}^\text{Cyc}_A$, and then we do some linear algebra as called for by the task at hand. Once again, we obtain our bases of $\mathcal{H}^\text{Weight}_A$ and $\mathcal{H}^\text{Cyc}_A$ by constructing suitable oriented cycles of $L$ on $A$. This is exactly what we are going to do for grid $\Psi$ in Sections \ref{flex-i} and \ref{flex-ii}.

This completes the outline of our method.

Granted, this approach is not guaranteed to always work. It might happen that one step of it or another is far too difficult to implement, or that the conditions of Theorem \ref{flexibility-method} are not met even though $\mathcal{F}_A$ is, in fact, flexible. Note, however, that it is far from obvious whether the latter issue ever really occurs or not. We revisit this question in Section \ref{further}.

Lastly, let us address some practical concerns.

Suppose that we wish to prove that $\mathcal{F}_A$ is flexible, but the number of rhombic classes of $L$ on $A$ is large. This could make it difficult to apply our method directly. In such cases, it might be best to look for some larger set of cells $B$ such that $A$ is a subset of $B$ (and so $\mathcal{F}_A$ is a subframework of $\mathcal{F}_B$) and $\mathcal{F}_B$ is still flexible. Then we could apply our method to $\mathcal{F}_B$ instead. Provided that we manage to fine-tune $B$ so that the number of rhombic classes of $L$ on $B$ is reasonably small, in this way our task could become vastly easier.

For example, it appears to be rather difficult to approach the flexibility part of Theorem \ref{square} directly. Instead, in Sections \ref{flex-i} and \ref{flex-ii} we prove Theorem \ref{flexibility-psi}, about the larger but more tractable grid $\Psi$, and then we derive Theorem \ref{square} from it.

Another example is the case $p = 1$ of Theorem \ref{flexibility-psi}. In it, the advantages of studying a larger set of cells instead are so great that they outweigh all troubles caused by that larger set of cells being oddly shaped, rather than a rectangular grid.

We go on to apply our method to grid $\Psi$ in Sections \ref{flex-i} and \ref{flex-ii}.

\section{Flexibility I} \label{flex-i}

Our main goal in this as well as the next section will be to prove the following theorem.

\medbreak

\begin{theorem} Let $\Psi$ be the rectangular grid of height $2(p + q) - 2$ and width $2(p + q) - 1$. Then the leaper framework of $L$ on $\Psi$ is flexible. \label{flexibility-psi} \end{theorem}

\medbreak

Before we embark, though, let us first say a couple of words about the significance of Theorem \ref{flexibility-psi}. Theorems \ref{rigidity-phi} and \ref{flexibility-psi}, taken together, completely resolve the question of the flexibility and rigidity of leaper frameworks on square grids. The full classification, disregarding the trivial case of the $1 \times 1$ grid, is as follows.

\medbreak

\begin{theorem} The leaper framework of $L$ on the square grid of side $n \ge 2$ is flexible if and only if $n \le 2(p + q) - 2$, and it is rigid if and only if $n \ge 2(p + q) - 1$. \nolinebreak \label{square} \end{theorem}

\medbreak

\begin{proof} By Lemma \ref{monotonicity} and Theorems \ref{rigidity-phi} and \ref{flexibility-psi}. \end{proof}

\medbreak

We move on to the proof of Theorem \ref{flexibility-psi}. We split the proof in two parts, as follows. In this section, we are going to take care of all free leapers $L$ such that $p \ge 2$. The case when $p = 1$ is, in a sense, degenerate, and we handle it separately in Section \ref{flex-ii}. For the most part, the two cases are amenable to the same techniques; however, in the case when $p = 1$ there are some additional obstacles to overcome and some steps of the argument need to be done slightly differently.

(Recall that we split the proof of Theorem \ref{inf-flexibility-psi} in two parts in the exact same way.)

Suppose, then, from this point on throughout the rest of this section, that $p \ge 2$.

Our plan for the proof of Theorem \ref{flexibility-psi} will be simply to apply Theorem \ref{flexibility-method}. This is, of course, easier said than done. The main difficulty is to find sufficiently nice bases of $\mathcal{H}^\text{Weight}_\Psi$ and $\mathcal{H}^\text{Cyc}_\Psi$. Solving this problem will take up the bulk of the proof.

We have already carried out most of the requisite preparatory work in Section \ref{inf-flex}. Let us review it briefly.

We denote the rhombic classes of $L$ on $\Psi$ by $R_1$, $R_2$, \ldots, $R_8$. Rhombic classes $R_1$ and $R_2$ together contain all moves of $L$ on $\Psi$ of slope $\frac{q}{p}$. Rhombic classes $R_3$ and $R_4$ are the reflections of $R_1$ and $R_2$, respectively, across the vertical axis of symmetry of $\Psi$, and together they contain all moves of $L$ on $\Psi$ of slope $-\frac{q}{p}$. Rhombic classes $R_5$ and $R_6$ together contain all moves of $L$ on $\Psi$ of slope $\frac{p}{q}$. Lastly, rhombic classes $R_7$ and $R_8$ are the reflections of $R_5$ and $R_6$, respectively, across the vertical axis of symmetry of $\Psi$, and together they contain all moves of $L$ on $\Psi$ of slope $-\frac{p}{q}$.

By Lemma \ref{three-slopes-symmetry}, each rhombic class of $L$ on $\Psi$ is centrally symmetric with respect to the center of $\Psi$, and so horizontal reflection on $\Psi$ swaps the same pairs of rhombic classes as vertical reflection on $\Psi$ does.

By the proof of Theorem \ref{inf-flexibility-psi}, we know that the dimension of $\mathcal{H}^\text{Weight}_\Psi$ is at most three, and so also that the dimension of the cyclic constraint space $\mathcal{H}^\text{Cyc}_\Psi$ of $L$ on $\Psi$ is at most six.

This completes our overview of prerequisites.

In the notation of Sections \ref{inf-flex} and \ref{method}, we assign directions to rhombic classes as follows: $d_1 = d_2 = (p, q)$, $d_3 = d_4 = (-p, q)$, $d_5 = d_6 = (q, p)$, and $d_7 = d_8 = (-q, p)$. (In Section \ref{inf-flex}, we never had to do this explicitly; however, the proof of Theorem \ref{flexibility-psi} requires greater sensitivity to detail.)

We set out to find a nice basis of $\mathcal{H}^\text{Weight}_\Psi$. As per Sections \ref{inf-flex} and \ref{method}, then we will be able to derive a nice basis of $\mathcal{H}^\text{Cyc}_\Psi$ from it immediately.

First we make one general observation about the structure of $\mathcal{H}^\text{Weight}_\Psi$ which makes very little use of anything specific to $\Psi$. Given any vector $\omega = (\omega_1, \omega_2, \omega_3, \omega_4,\allowbreak \omega_5, \omega_6, \omega_7, \omega_8)$ in $\mathcal{H}^\text{Weight}_\Psi$, we write $\widehat{\omega}$ for the vector $(\omega_3, \omega_4, \omega_1, \omega_2, \omega_7, \omega_8, \omega_5, \omega_6)$. The significance of this permutation of the components of $\omega$ is given by the following lemma.

\medbreak

\begin{lemma} Let $\omega$ be any vector in $\mathcal{H}^\text{Weight}_\Psi$. Then $\widehat{\omega}$ belongs to $\mathcal{H}^\text{Weight}_\Psi$ as well. \nolinebreak \label{weight-psi-reflection} \end{lemma}

\medbreak

\begin{proof} It suffices to consider the case when $\omega$ is the weight of some oriented cycle $C$ of $L$ on $\Psi$. Then $\widehat{\omega}$ is the weight of the reflection of $C$ across the vertical axis of symmetry of $\Psi$. \end{proof}

\medbreak

Similarly to what we did for grids $\Phi$ and $\Theta$ in Section \ref{rigid}, our construction of a nice basis of $\mathcal{H}^\text{Weight}_\Psi$ will revolve around one exceptionally nice oriented cycle of $L$ on $\Psi$. The oriented cycles $C_\Phi$ and $C_\Theta$ which we constructed back then are not well-suited to this task because it is difficult to say anything sufficiently specific about their weights.

Once again, first we define our oriented cycle of $L$ on $\Psi$ as an oriented subgraph of $\mathcal{L}_\Psi$, and then we prove that it is in fact an oriented cycle. For this construction, recall that $p < q$ by convention.

Suppose that $\Psi$ is the grid $[1; 2(p + q) - 1] \times [1; 2(p + q) - 2]$. We define $C_\Psi$ as the subgraph of $\mathcal{L}_\Psi$ formed by the following moves of $L$ on $\Psi$:

(a) All moves in direction $(q, p)$ which start from a cell in the subgrid of $\Psi$ given by $[1; q] \times [1; q - p]$ .

(b) All moves in direction $(q, p)$ which start from a cell in the subgrid of $\Psi$ given by $[1; p + q] \times [q - p + 1, q]$.

(c) All moves in direction $(-q, p)$ which start from a cell in the subgrid of $\Psi$ given by $[q + 1; 2q] \times [1; q - p]$.

(d) All moves in direction $(-q, p)$ which start from a cell in the subgrid of $\Psi$ given by $[p + q + 1, 2q] \times [q - p + 1, q]$.

(e) Lastly, all moves in direction $(-p, -q)$ which lead to a cell in the subgrid of $\Psi$ given by $[1; 2q] \times [1; p]$.

\begin{figure}[t] \centering \includegraphics[scale=\figuresScale]{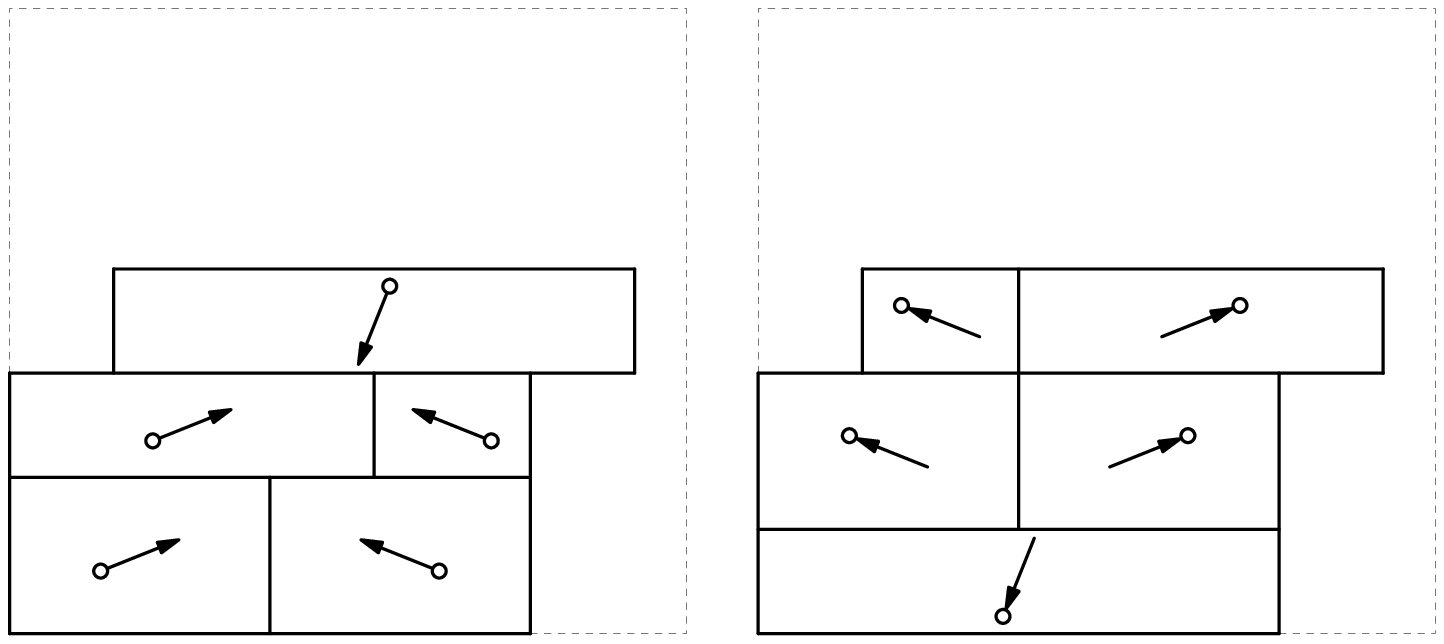} \caption{} \label{blueprint} \end{figure}

\begin{figure}[t] \centering \includegraphics[scale=\figuresScale]{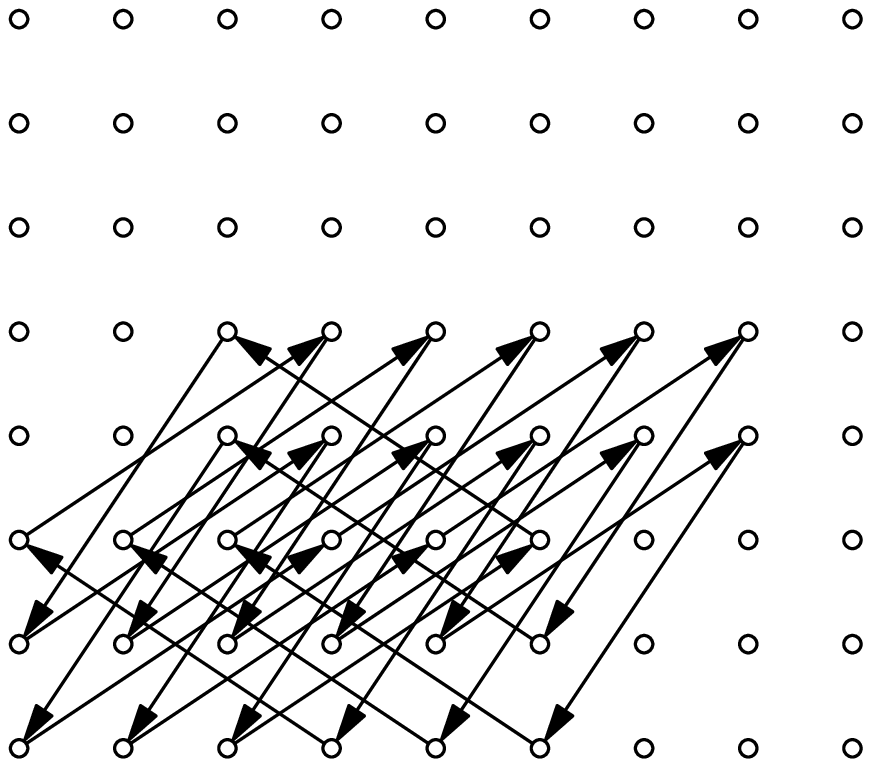} \caption{} \label{zebra-c-psi} \end{figure}

Of course, the endpoints of the moves of $L$ in each one of pencils (a), (b), (c), and (d), as well as the starting points of the moves of $L$ in pencil (e), form subgrids of $\Psi$ as well. Figure \ref{blueprint} shows how all ten subgrids fit together. The label of each subgrid indicates what pencil of moves of $L$ it is associated with. For example, Figure \ref{zebra-c-psi} shows $C_\Psi$ for the zebra.

It is straightforward to verify that every vertex of $C_\Psi$ is of in-degree one as
well as of out-degree one. Therefore, $C_\Psi$ is the disjoint union of several oriented
cycles. Just as in Section \ref{rigid}, this is really all that we need to know about the structure of $C_\Psi$ for the rest of our argument to go through; however, we are going to prove that $C_\Psi$ is an oriented cycle anyway. The proof is largely similar to the analogous proof for $C_\Phi$ in Section \ref{rigid}.

Observe that all moves of $L$ in $C_\Psi$ point in one of the three directions $(q, p)$, $(-q, p)$, and $(-p, -q)$.

Consider any connected component $C$ of $C_\Psi$. Suppose that $C$ contains a total of $\alpha$ moves of $L$ in direction $(q, p)$, a total of $\beta$ moves of $L$ in direction $(-q, p)$, and a total of $\gamma$ moves of $L$ in direction $(-p, -q)$. By summation along $C$, we obtain \[\alpha(q, p) + \beta(-q, p) + \gamma(-p, -q) = \mathbf{0}.\]

Since $p$ and $q$ are positive and distinct, it follows that all three of $\alpha$, $\beta$, and $\gamma$ are nonzero and that \[\alpha : \beta : \gamma = (p^2 + q^2) : (q^2 - p^2) : 2pq.\]

Since $p$ and $q$ are relatively prime positive integers of opposite parity, $p^2 + q^2$, $q^2 - p^2$, and $2pq$ are pairwise relatively prime. Thus $p^2 + q^2$, $q^2 - p^2$, and $2pq$ divide $\alpha$, $\beta$, and $\gamma$, respectively. Consequently, since $\alpha$, $\beta$, and $\gamma$ are positive integers, \[\alpha + \beta + \gamma \ge (p^2 + q^2) + (q^2 - p^2) + 2pq.\]

However, the right-hand side of this inequality equals the total number of moves of $L$ in $C_\Psi$. Therefore, we must have that $\alpha = p^2 + q^2$, $\beta = q^2 - p^2$, $\gamma = 2pq$, and $C$ coincides with $C_\Psi$.

This completes our proof that $C_\Psi$ is an oriented cycle of $L$ on $\Psi$, and we proceed to study the weight of $C_\Psi$.

\medbreak

\begin{lemma} Let the weight of $C_\Psi$ be $\mathbf{r} = (r_1, r_2, r_3, r_4, r_5, r_6, r_7, r_8)$. Then $r_1 + r_2 = -2pq$, $r_3 = r_4 = 0$, $r_5 + r_6 = p^2 + q^2$, and $r_7 + r_8 = q^2 - p^2$. Furthermore, $r_1$ and $r_2$ are negative, whereas $r_5$ and $r_6$ are positive. \label{weight-psi-nonzero-i} \end{lemma}

\medbreak

\begin{proof} The first part of the lemma is immediate by our earlier observation that all moves of $L$ in $C_\Psi$ point in one of the three directions $(q, p)$, $(-q, p)$, and $(-p, -q)$. For the same reason, $r_1$ and $r_2$ are nonpositive, whereas $r_5$ and $r_6$ are nonnegative. Thus we are left to show that all four of $r_1$, $r_2$, $r_5$, and $r_6$ are nonzero.

Let $\Pi'$ be the lower left subgrid of $\Psi$ of height $h_\Psi - q = 2p + q - 2$ and width $w_\Psi - p = p + 2q - 1$, and let $\Pi''$ be the lower left subgrid of $\Psi$ of height $h_\Psi - p = p + 2q - 2$ and width $w_\Psi - q = 2p + q - 1$. We first introduced grids $\Pi'$ and $\Pi''$ in the context of Lemma \ref{three-slopes-symmetry}; however, back then we did not anchor either one of them to any specific location on the integer lattice.

Let $\Omega_\texttt{LR}$ be the lower right subgrid of $\Pi'$ of height $p$ and width $2q$. Moreover, let cells $a_\texttt{I}$ and $a_\texttt{II}$ form the upper left subgrid of $\Omega_\texttt{LR}$ of height two and width one. Explicitly, $a_\texttt{I}$ and $a_\texttt{II}$ are cells $(p, p - 1)$ and $(p, p)$ of $\Psi$.

By Lemma \ref{three-slopes-components}, cells $a_\texttt{I}$ and $a_\texttt{II}$ are representatives of the two connected components of graph $\mathcal{L}_{\Pi'} \setminus \frac{q}{p}$.

Let $e_\texttt{I}$ and $e_\texttt{II}$ be the two moves of $L$ on $\Psi$ which point in direction $(-p, -q)$ and lead to cells $a_\texttt{I}$ and $a_\texttt{II}$, respectively. Then, by Lemma \ref{rhombic-class}, $e_\texttt{I}$ and $e_\texttt{II}$ are representatives of the two rhombic classes of $L$ on $\Psi$ for slope $\frac{q}{p}$. Since both of $e_\texttt{I}$ and $e_\texttt{II}$ belong to $C_\Psi$, we conclude that $r_1$ and $r_2$ are nonzero.

We handle components $r_5$ and $r_6$ of $\mathbf{r}$ analogously, as follows.

Let $\Omega_\texttt{UL}$ be the upper left subgrid of $\Pi''$ of height $2q$ and width $p$. Consider the lower right subgrid of $\Omega_\texttt{UL}$ of height two and width one. As luck would have it, this subgrid of $\Omega_\texttt{UL}$ consists exactly of cells $a_\texttt{I}$ and $a_\texttt{II}$, as defined above.

By Lemma \ref{three-slopes-components}, then, cells $a_\texttt{I}$ and $a_\texttt{II}$ are representatives of the two connected components of graph $\mathcal{L}_{\Pi''} \setminus \frac{p}{q}$.

Let $e_\texttt{V}$ and $e_\texttt{VI}$ be the two moves of $L$ on $\Psi$ which point in direction $(q, p)$ and start from cells $a_\texttt{I}$ and $a_\texttt{II}$, respectively. Then, by Lemma \ref{rhombic-class}, $e_\texttt{V}$ and $e_\texttt{VI}$ are representatives of the two rhombic classes of $L$ on $\Psi$ for slope $\frac{p}{q}$. Since both of $e_\texttt{V}$ and $e_\texttt{VI}$ belong to $C_\Psi$, we conclude that $r_5$ and $r_6$ are nonzero as well. \end{proof}

\medbreak

Note that Lemma \ref{weight-psi-nonzero-i} keeps silent regarding the chances of components $r_7$ and $r_8$ of $\mathbf{r}$ to equal zero. Indeed, Lemma \ref{three-slopes-components} does not provide us with suitable representatives in $C_\Psi$ of the rhombic classes $R_7$ and $R_8$ of $L$ on $\Psi$. Fortunately, what Lemma \ref{weight-psi-nonzero-i} does tell us will prove to be sufficient for our purposes.

Let $C'_\Psi$ be the reflection of $C_\Psi$ in the line $x = y$. Here follows the analogue of Lemma \ref{weight-psi-nonzero-i} for $C'_\Psi$.

\medbreak

\begin{lemma} Let the weight of $C'_\Psi$ be $\mathbf{s} = (s_1, s_2, s_3, s_4, s_5, s_6, s_7, s_8)$. Then $s_1 + s_2 = p^2 + q^2$, $s_3 + s_4 = p^2 - q^2$, $s_5 + s_6 = -2pq$, and $s_7 = s_8 = 0$. Furthermore, $s_1$ and $s_2$ are positive, whereas $s_5$ and $s_6$ are negative. \label{weight-psi-nonzero-ii} \end{lemma}

\medbreak

\begin{proof} Analogous to the proof of Lemma \ref{weight-psi-nonzero-i}. In particular, when we reverse the directions of moves $e_\texttt{I}$, $e_\texttt{II}$, $e_\texttt{V}$, and $e_\texttt{VI}$ as defined in the proof of Lemma \ref{weight-psi-nonzero-i}, we obtain representatives in $C'_\Psi$ of the rhombic classes $R_1$, $R_2$, $R_5$, and $R_6$ of $L$ on $\Psi$. \end{proof}

\medbreak

The weights of $C_\Psi$ and $C'_\Psi$ are special in the sense that their components in certain slopes of $L$ equal zero. The main point of the following lemma is that, apart from that, these weights are ``sufficiently arbitrary''. Note that each vector in Lemma \ref{weight-psi-independence}'s statement is associated with one slope of $L$.

\medbreak

\begin{lemma} Vectors $(r_1, r_2)$, $(s_1, s_2)$, and $(s_3, s_4)$ are pairwise linearly independent. So are vectors $(r_5, r_6)$, $(r_7, r_8)$, and $(s_5, s_6)$. \label{weight-psi-independence} \end{lemma}

\medbreak

\begin{proof} We consider the case of vectors $(r_1, r_2)$ and $(s_3, s_4)$ in detail, and all other cases are analogous.

Suppose, for the sake of contradiction, that \[\alpha(r_1, r_2) + \beta(s_3, s_4) = \mathbf{0},\] where $\alpha$ and $\beta$ are not both zero.

Then \[\alpha(r_1 + r_2) + \beta(s_3 + s_4) = 0.\]

By Lemmas \ref{weight-psi-nonzero-i} and \ref{weight-psi-nonzero-ii}, this amounts to \[\alpha : \beta = (p^2 - q^2) : 2pq.\]

Consequently, \begin{align*} (p^2 - q^2)r_1 + 2pqs_3 &= 0 \text{ and }\\ (p^2 - q^2)r_2 + 2pqs_4 &= 0. \end{align*}

Since $p$ and $q$ are relatively prime positive integers of opposite parity, $p^2 - q^2$ and $2pq$ are relatively prime. Thus $2pq$ divides both of $r_1$ and $r_2$. However, by Lemma \ref{weight-psi-nonzero-i}, both of $r_1$ and $r_2$ are negative. Therefore, \[r_1 + r_2 \le (-2pq) + (-2pq),\] and we arrive at a contradiction with Lemma \ref{weight-psi-nonzero-i}. \end{proof}

\medbreak

By Lemma \ref{weight-psi-reflection}, both vectors $\widehat{\mathbf{r}}$ and $\widehat{\mathbf{s}}$ are in $\mathcal{H}^\text{Weight}_\Psi$ as well. Building upon Lemma \ref{weight-psi-independence}, we proceed to construct our long sought for basis.

\medbreak

\begin{lemma} Any three vectors out of $\mathbf{r}$, $\widehat{\mathbf{r}}$, $\mathbf{s}$, and $\widehat{\mathbf{s}}$ form a basis of $\mathcal{H}^\text{Weight}_\Psi$. \label{weight-psi-basis} \end{lemma}

\medbreak

\begin{proof} We know already that the dimension of $\mathcal{H}^\text{Weight}_\Psi$ is at most three. Thus it suffices to show that any three vectors out of $\mathbf{r}$, $\widehat{\mathbf{r}}$, $\mathbf{s}$, and $\widehat{\mathbf{s}}$ are linearly independent. We consider the case of vectors $\mathbf{r}$, $\widehat{\mathbf{r}}$, and $\mathbf{s}$ in detail, and all other cases are analogous.

Suppose that \[\alpha\mathbf{r} + \beta\widehat{\mathbf{r}} + \gamma\mathbf{s} = \mathbf{0}\] for some real numbers $\alpha$, $\beta$, and $\gamma$.

Then, looking just at the first two components of all three vectors, we obtain \[\alpha(r_1, r_2) + \gamma(s_1, s_2) = \mathbf{0}.\]

By Lemma \ref{weight-psi-independence}, it follows that $\alpha = \gamma = 0$. Since vector $\widehat{\mathbf{r}}$ is nonzero, then $\beta = 0$ as well. \end{proof}

\medbreak

Given any vector $\omega$ in $\mathcal{H}^\text{Weight}_\Psi$, we write $\omega_\texttt{X}$ for the cyclic constraint vector associated with $\omega$ which only addresses the $x_i$, namely $(\omega_1, 0, \omega_2, 0, \ldots, \omega_8, 0)$. Analogously, we write $\omega_\texttt{Y}$ for the cyclic constraint vector associated with $\omega$ which only addresses the $y_i$, namely $(0, \omega_1, 0, \omega_2, \ldots, 0, \omega_8)$.

By Lemma \ref{weight-psi-basis}, vectors $\mathbf{r}_\texttt{X}$, $\mathbf{r}_\texttt{Y}$, $\widehat{\mathbf{r}}_\texttt{X}$, $\widehat{\mathbf{r}}_\texttt{Y}$, $\mathbf{s}_\texttt{X}$, and $\mathbf{s}_\texttt{Y}$ form a basis of $\mathcal{H}^\text{Cyc}_\Psi$.

We are ready to verify the conditions of Theorem \ref{flexibility-method}. We use the vectors of the coefficients of all directional constraints as a basis of $\mathcal{H}^\text{Dir}_\Psi$, as well as vectors $\mathbf{r}_\texttt{X}$, $\mathbf{r}_\texttt{Y}$, $\widehat{\mathbf{r}}_\texttt{X}$, $\widehat{\mathbf{r}}_\texttt{Y}$, $\mathbf{s}_\texttt{X}$, and $\mathbf{s}_\texttt{Y}$ as a basis of $\mathcal{H}^\text{Cyc}_\Psi$. Collected together, all fourteen of these vectors form the matrix \[\setcounter{MaxMatrixCols}{16} M_\Psi = \begin{pmatrix} p & q & 0 & 0 & 0 & 0 & 0 & 0 & 0 & 0 & 0 & 0 & 0 & 0 & 0 & 0\\ 0 & 0 & p & q & 0 & 0 & 0 & 0 & 0 & 0 & 0 & 0 & 0 & 0 & 0 & 0\\ 0 & 0 & 0 & 0 & -p & q & 0 & 0 & 0 & 0 & 0 & 0 & 0 & 0 & 0 & 0\\ 0 & 0 & 0 & 0 & 0 & 0 & -p & q & 0 & 0 & 0 & 0 & 0 & 0 & 0 & 0\\ 0 & 0 & 0 & 0 & 0 & 0 & 0 & 0 & q & p & 0 & 0 & 0 & 0 & 0 & 0\\ 0 & 0 & 0 & 0 & 0 & 0 & 0 & 0 & 0 & 0 & q & p & 0 & 0 & 0 & 0\\ 0 & 0 & 0 & 0 & 0 & 0 & 0 & 0 & 0 & 0 & 0 & 0 & -q & p & 0 & 0\\ 0 & 0 & 0 & 0 & 0 & 0 & 0 & 0 & 0 & 0 & 0 & 0 & 0 & 0 & -q & p\\ r_1 & 0 & r_2 & 0 & 0 & 0 & 0 & 0 & r_5 & 0 & r_6 & 0 & r_7 & 0 & r_8 & 0\\ 0 & r_1 & 0 & r_2 & 0 & 0 & 0 & 0 & 0 & r_5 & 0 & r_6 & 0 & r_7 & 0 & r_8\\ 0 & 0 & 0 & 0 & r_1 & 0 & r_2 & 0 & r_7 & 0 & r_8 & 0 & r_5 & 0 & r_6 & 0\\ 0 & 0 & 0 & 0 & 0 & r_1 & 0 & r_2 & 0 & r_7 & 0 & r_8 & 0 & r_5 & 0 & r_6\\ s_1 & 0 & s_2 & 0 & s_3 & 0 & s_4 & 0 & s_5 & 0 & s_6 & 0 & 0 & 0 & 0 & 0\\ 0 & s_1 & 0 & s_2 & 0 & s_3 & 0 & s_4 & 0 & s_5 & 0 & s_6 & 0 & 0 & 0 & 0 \end{pmatrix}.\]

Our task, then, is to show that the rows of $M_\Psi$ are linearly independent.

Consider any linear combination of the rows of $M_\Psi$ which equals the zero vector. Let the coefficients of the first eight rows of $M_\Psi$ (equivalently, of the basis of $\mathcal{H}^\text{Dir}_\Psi$) be $\alpha_1$, $\alpha_2$, \ldots, $\alpha_8$, respectively, and let the coefficients of the last six rows of $M_\Psi$ (equivalently, of the basis of $\mathcal{H}^\text{Cyc}_\Psi$) be $\beta_1$, $\beta_2$, \ldots, $\beta_6$, respectively.

By columns $1$ and $3$ of $M_\Psi$, we have that \[p(\alpha_1, \alpha_2) + \beta_1(r_1, r_2) + \beta_5(s_1, s_2) = \mathbf{0}.\]

Analogously, by columns $2$ and $4$ of $M_\Psi$, we have that \[q(\alpha_1, \alpha_2) + \beta_2(r_1, r_2) + \beta_6(s_1, s_2) = \mathbf{0}.\]

Eliminating $(\alpha_1, \alpha_2)$, we obtain \[(\beta_1 q - \beta_2 p)(r_1, r_2) + (\beta_5 q - \beta_6 p)(s_1, s_2) = \mathbf{0}.\]

By Lemma \ref{weight-psi-independence}, it follows that \begin{align*} \beta_1 q - \beta_2 p &= 0 \text{ and } \tag{\textbf{A}}\\ \beta_5 q - \beta_6 p &= 0. \tag{\textbf{B}} \end{align*}

Analogously, by columns $5$, $6$, $7$, and $8$ of $M_\Psi$, and taking into account Lemma \ref{weight-psi-independence}, we obtain \begin{align*} \beta_3 q + \beta_4 p &= 0 \text{ and } \tag{\textbf{C}}\\ \beta_5 q + \beta_6 p &= 0. \tag{\textbf{D}} \end{align*}

Lastly, by columns $13$, $14$, $15$, and $16$ of $M_\Psi$, and taking into account Lemma \ref{weight-psi-independence}, we obtain \begin{align*} \beta_1 p + \beta_2 q &= 0 \text{ and } \tag{\textbf{E}}\\ \beta_3 p + \beta_4 q &= 0. \tag{\textbf{F}} \end{align*}

Since $p$ and $q$ are positive and distinct, by (\textbf{A}) and (\textbf{E}) we have that $\beta_1 = \beta_2 = 0$; by (\textbf{C}) and (\textbf{F}) we have that $\beta_3 = \beta_4 = 0$; and, finally, by (\textbf{B}) and (\textbf{D}) we have that $\beta_5 = \beta_6 = 0$.

Thus all six of $\beta_1$, $\beta_2$, \ldots, $\beta_6$ must equal zero. On the other hand, since the first eight rows of $M_\Psi$ are linearly independent, all eight of $\alpha_1$, $\alpha_2$, \ldots, $\alpha_8$ must then equal zero as well. Therefore, the rows of $M_\Psi$ are linearly independent, as required.

This completes our proof of Theorem \ref{flexibility-psi} in the case when $p \ge 2$.

As an aside, we obtain the following corollary.

\medbreak

\begin{corollary} Suppose that $p \ge 2$. Then the leaper framework of $L$ on $\Psi$ flexes with precisely one degree of freedom. \label{one-degree-psi} \end{corollary}

\medbreak

\begin{proof} By Corollary \ref{degrees-of-freedom} and the proof of Theorem \ref{flexibility-psi}. \end{proof}

\medbreak

\begin{figure}[p] \centering \begin{tabular}{>{\centering}m{0.30\textwidth} >{\centering}m{0.30\textwidth} >{\centering}m{0.30\textwidth}} \includegraphics[scale=\stillFramesScale]{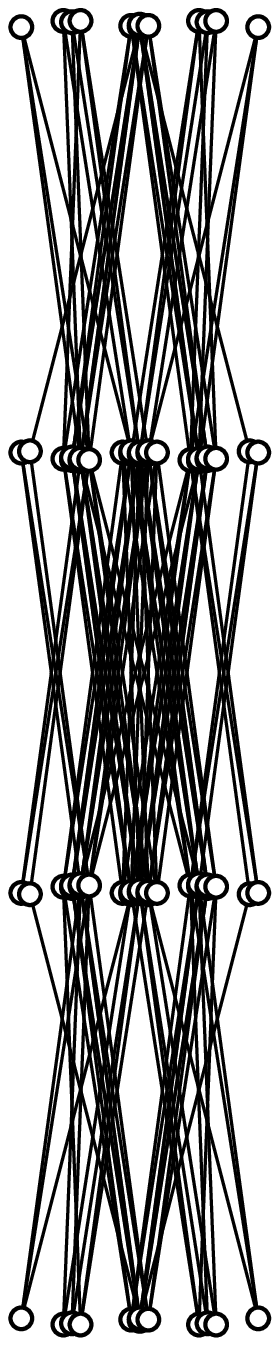} & \includegraphics[scale=\stillFramesScale]{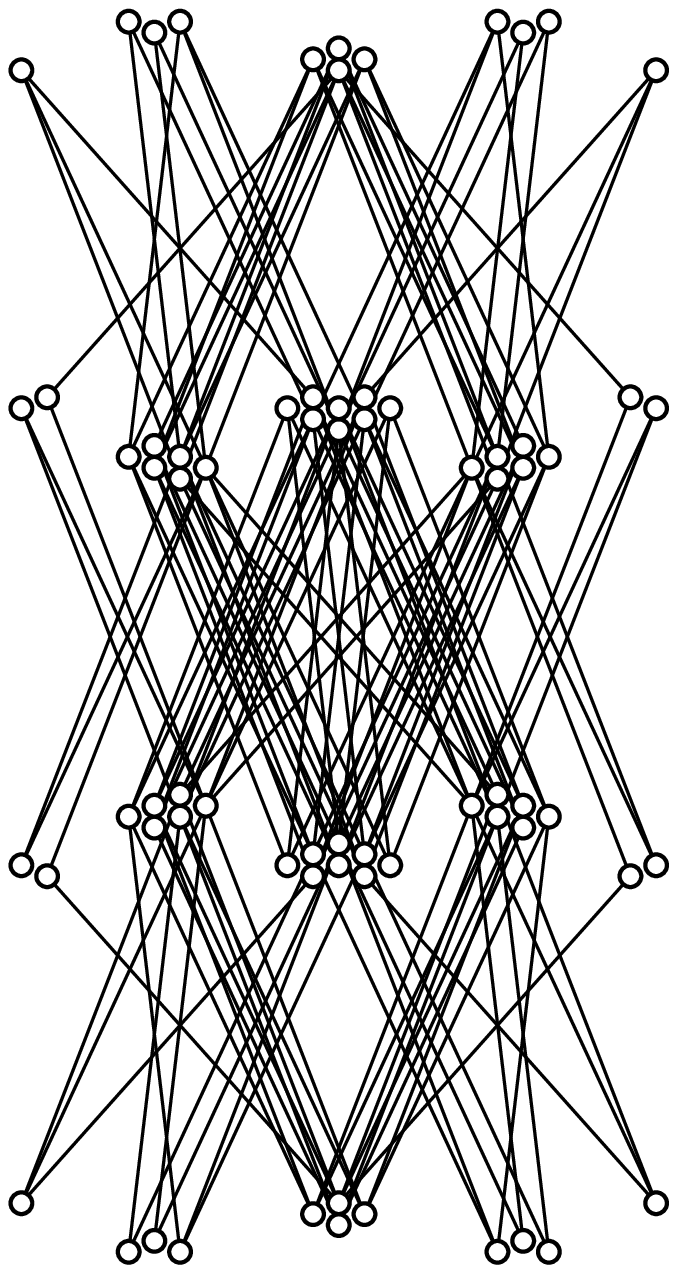} & \includegraphics[scale=\stillFramesScale]{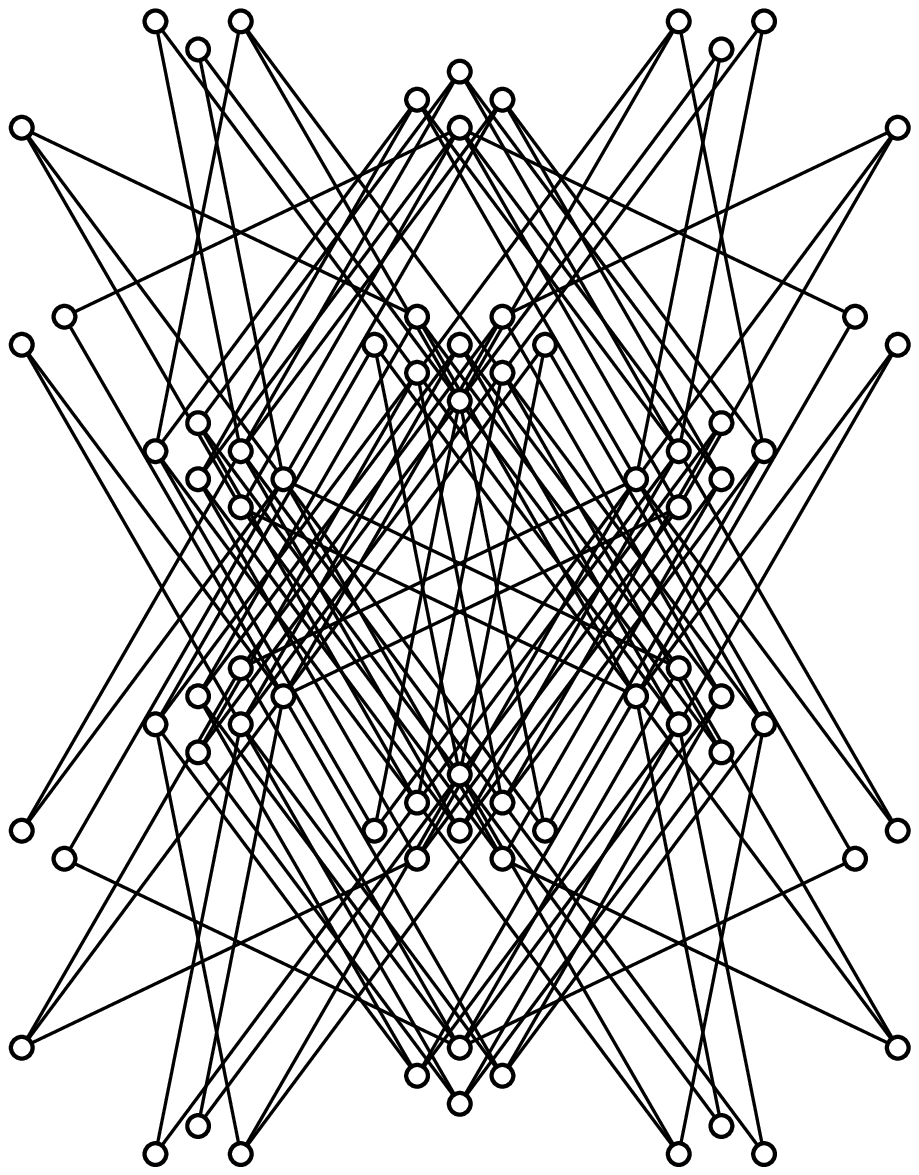} \tabularnewline \smallbreak \tabularnewline \includegraphics[scale=\stillFramesScale]{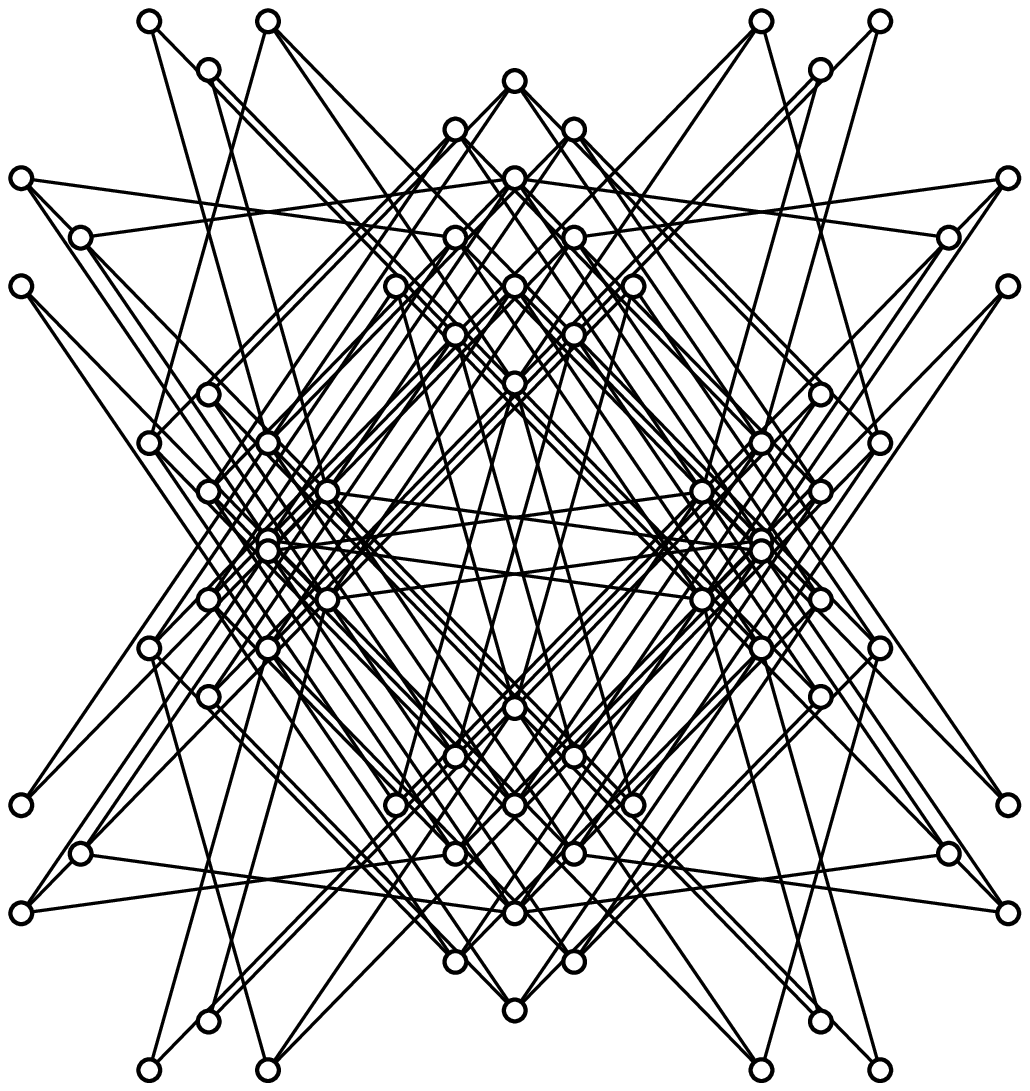} & \includegraphics[scale=\stillFramesScale]{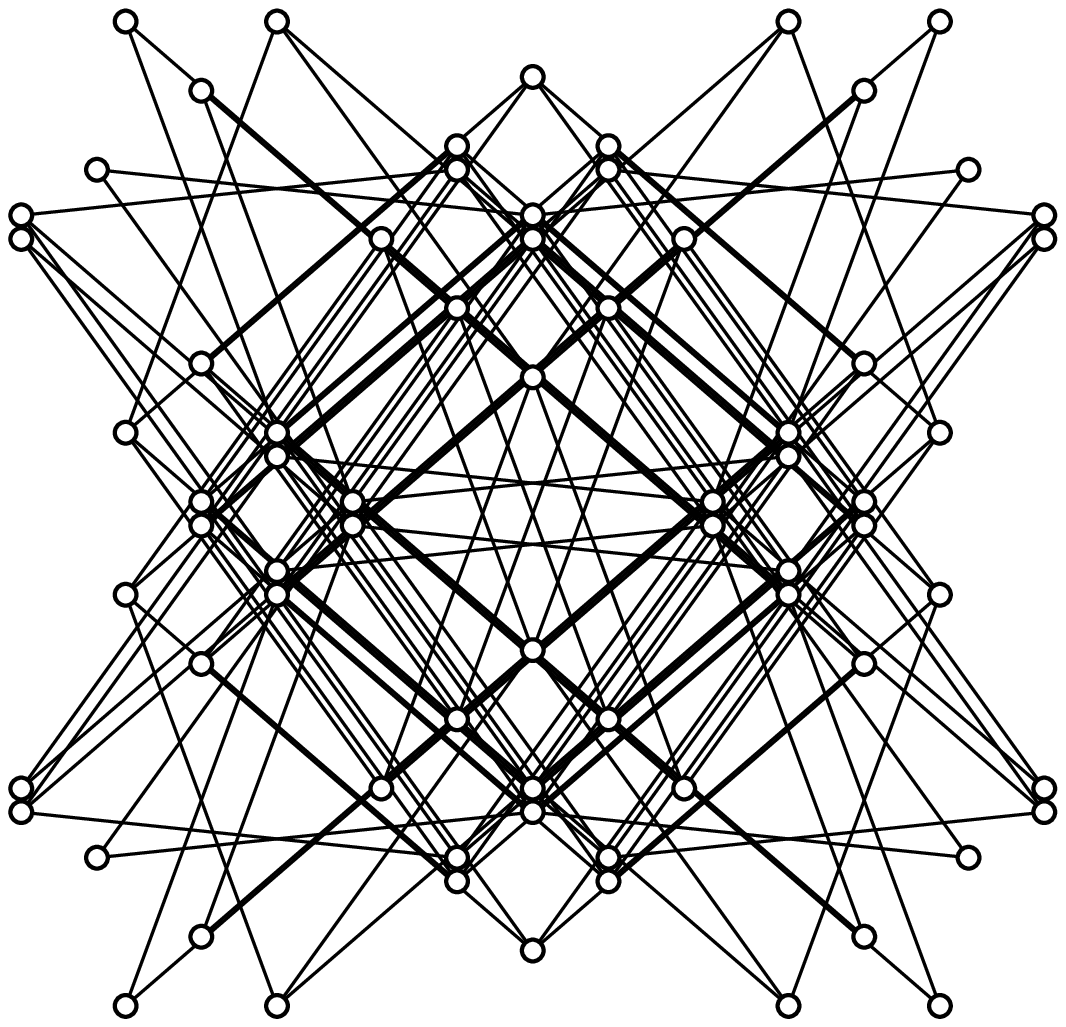} & \includegraphics[scale=\stillFramesScale]{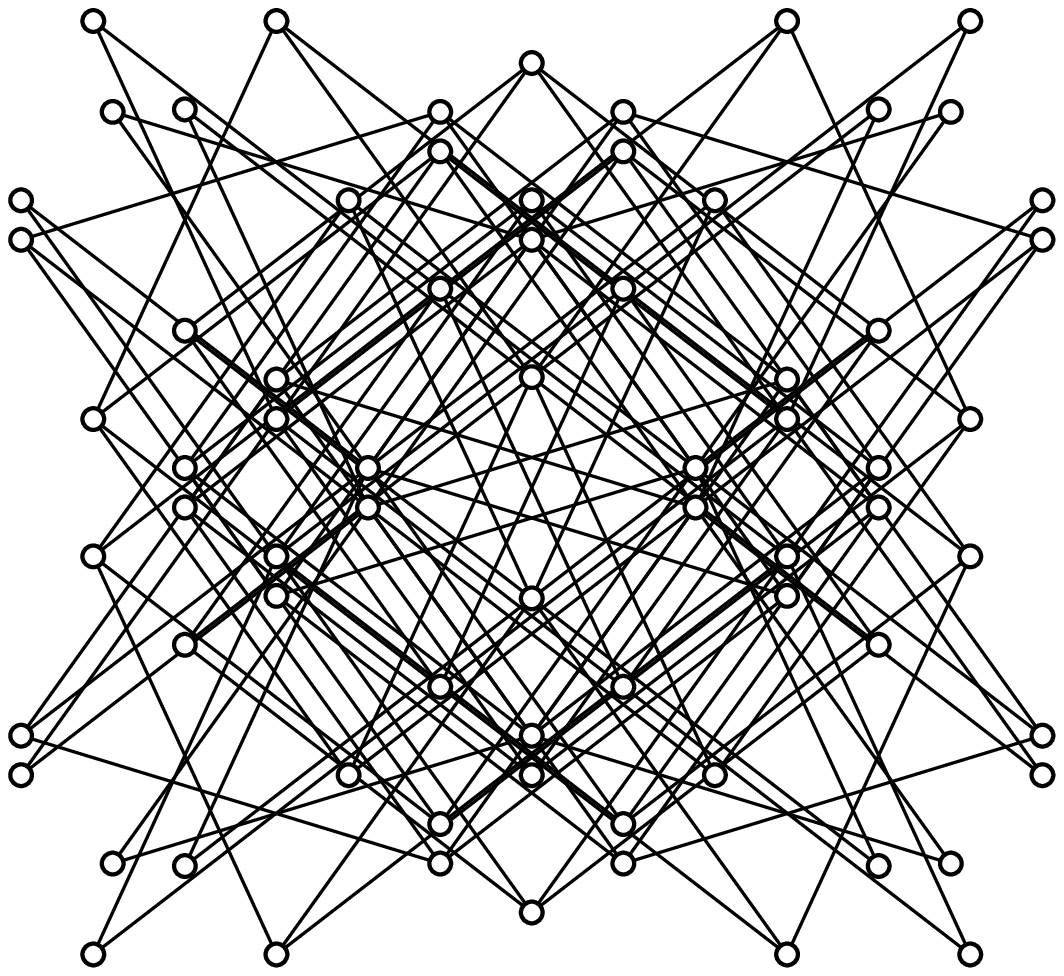} \tabularnewline \smallbreak \tabularnewline \includegraphics[scale=\stillFramesScale]{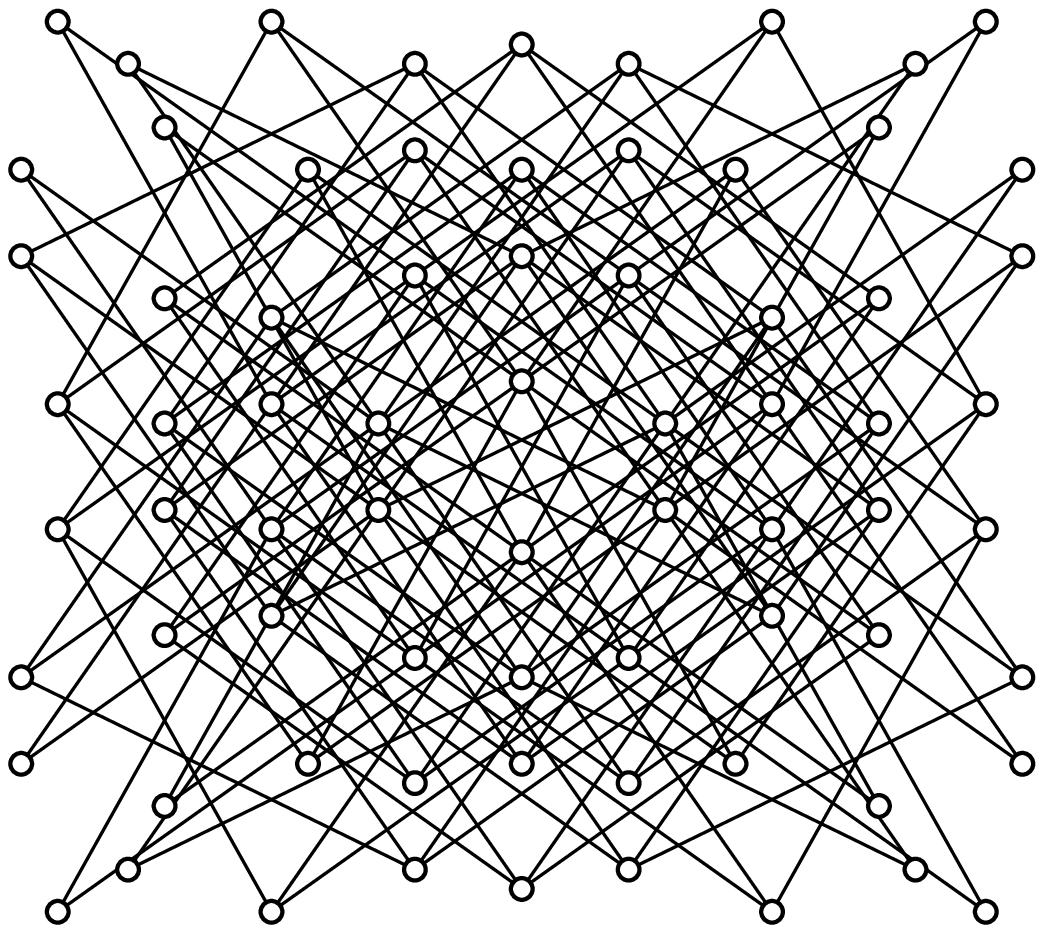} & \includegraphics[scale=\stillFramesScale]{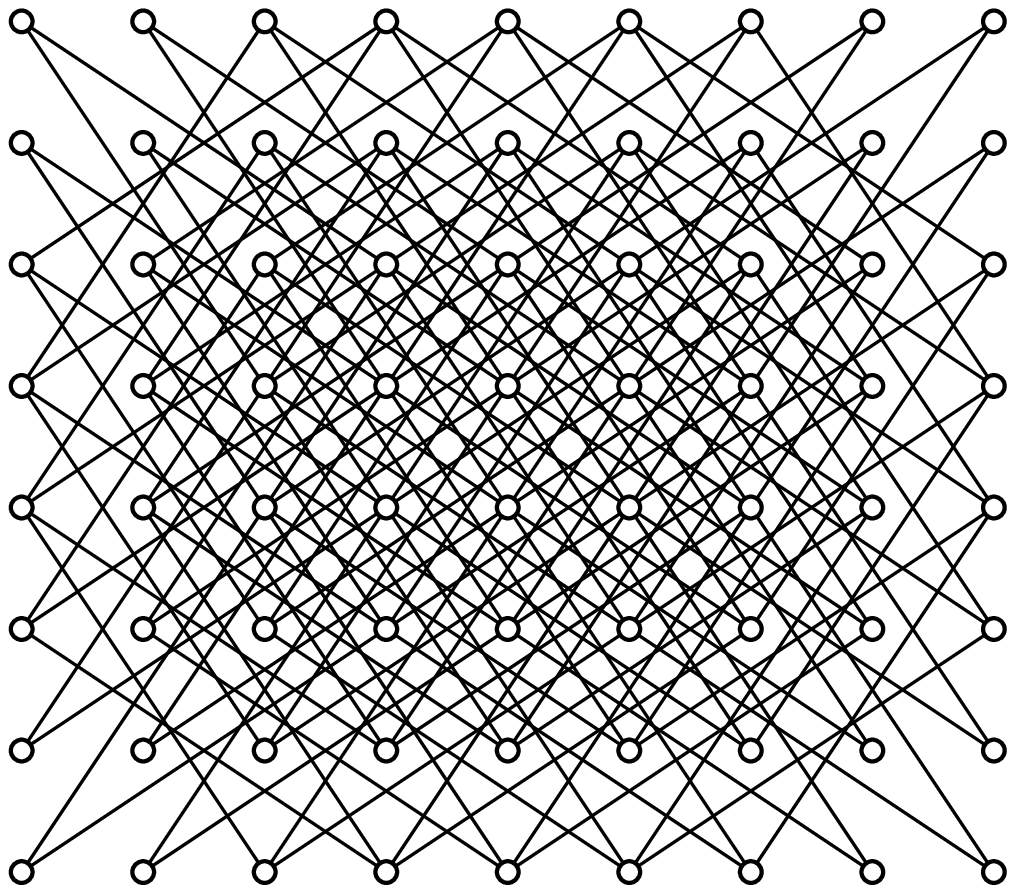} & \includegraphics[scale=\stillFramesScale]{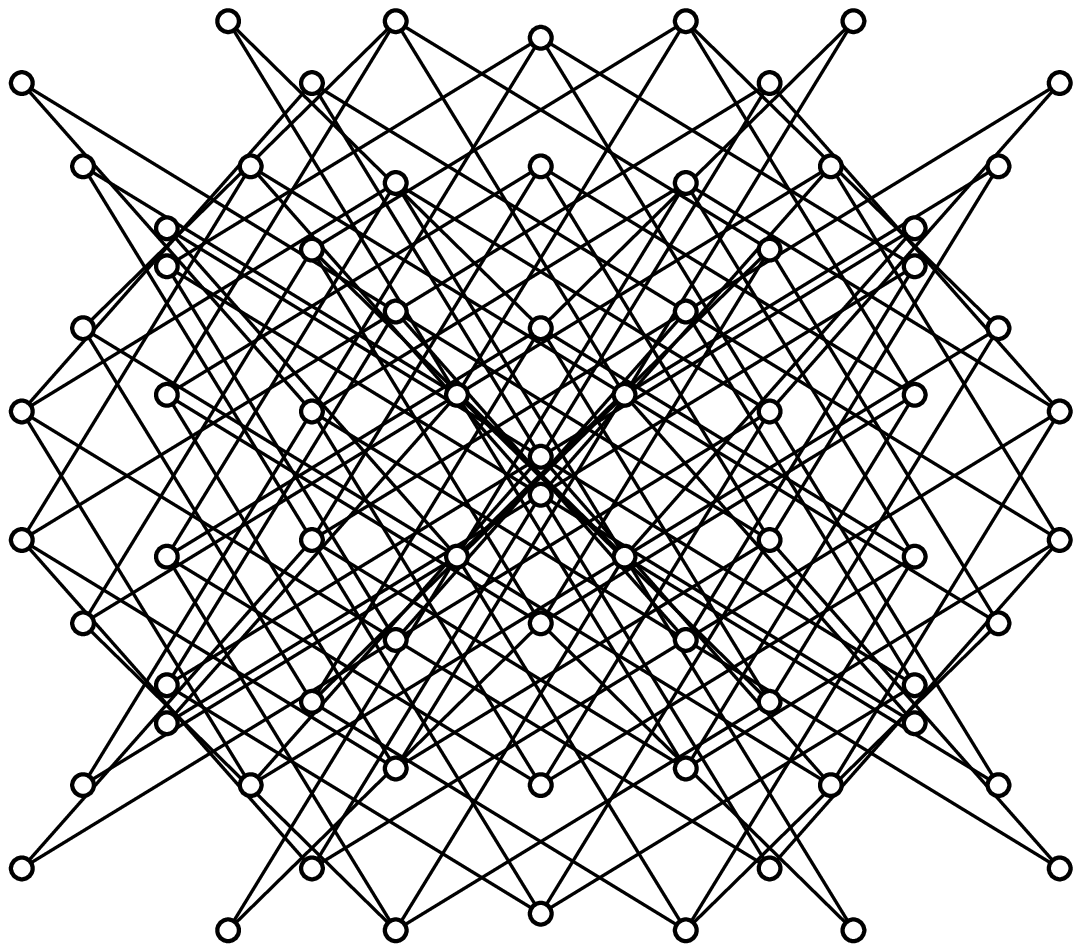} \tabularnewline \smallbreak \tabularnewline \includegraphics[scale=\stillFramesScale]{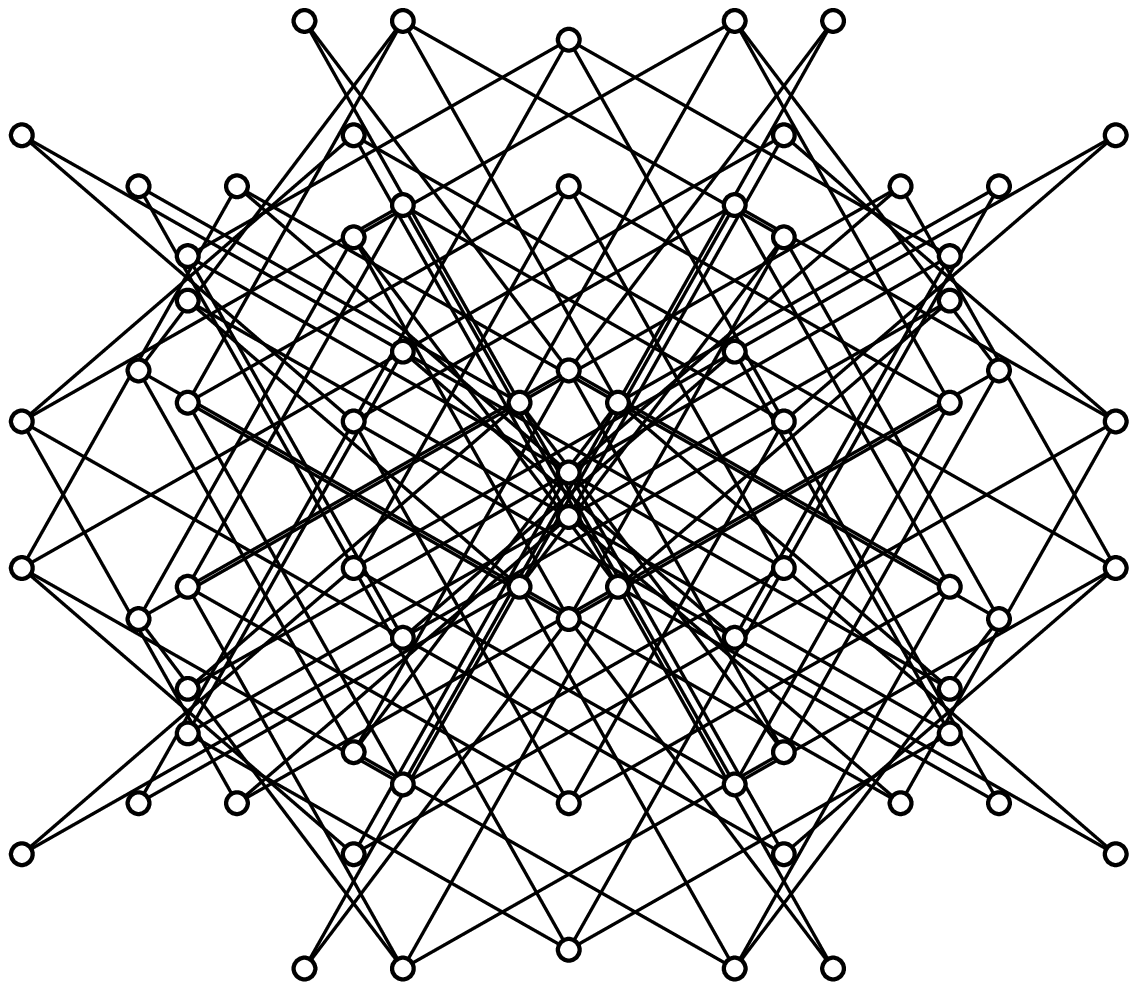} & \includegraphics[scale=\stillFramesScale]{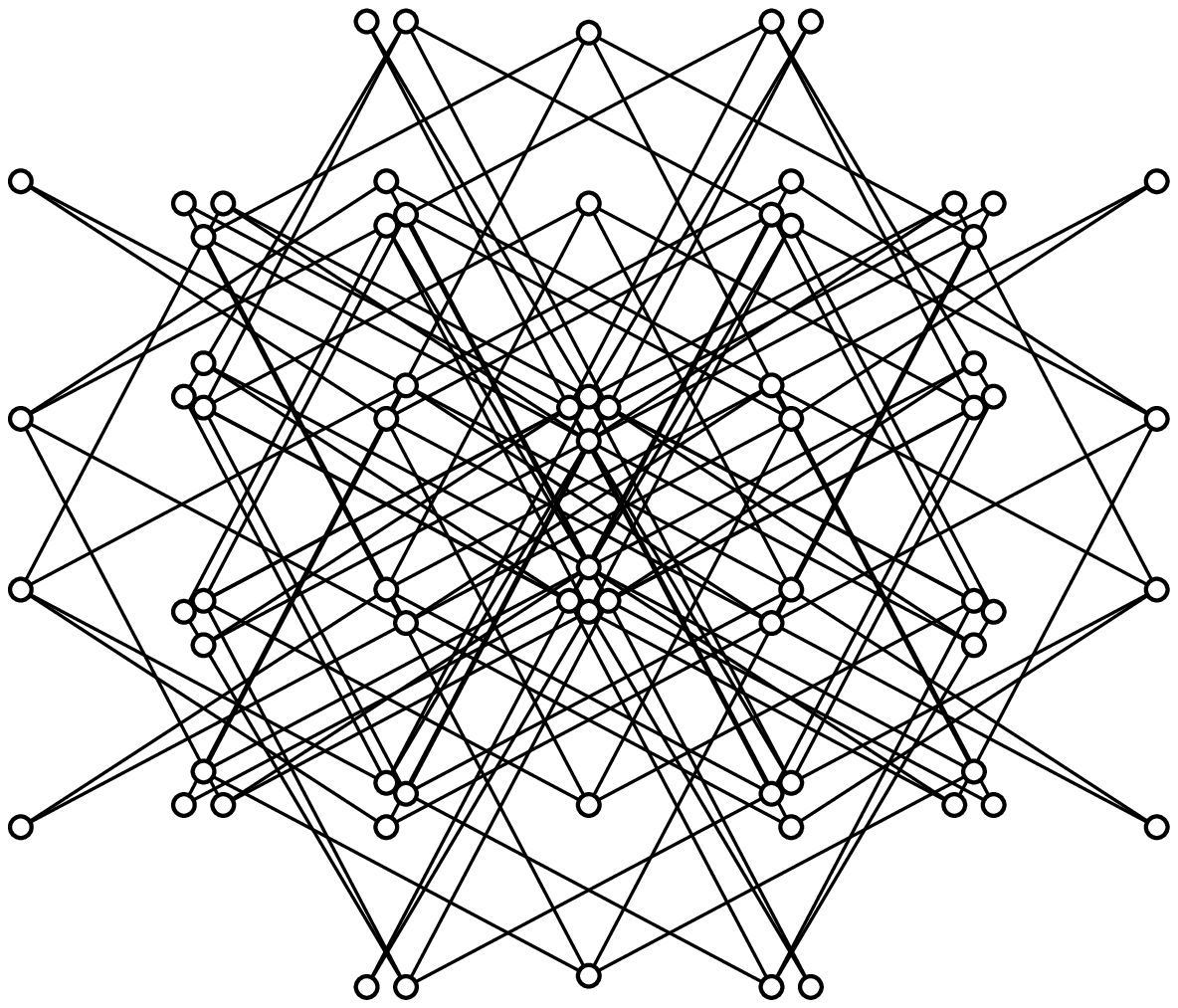} & \includegraphics[scale=\stillFramesScale]{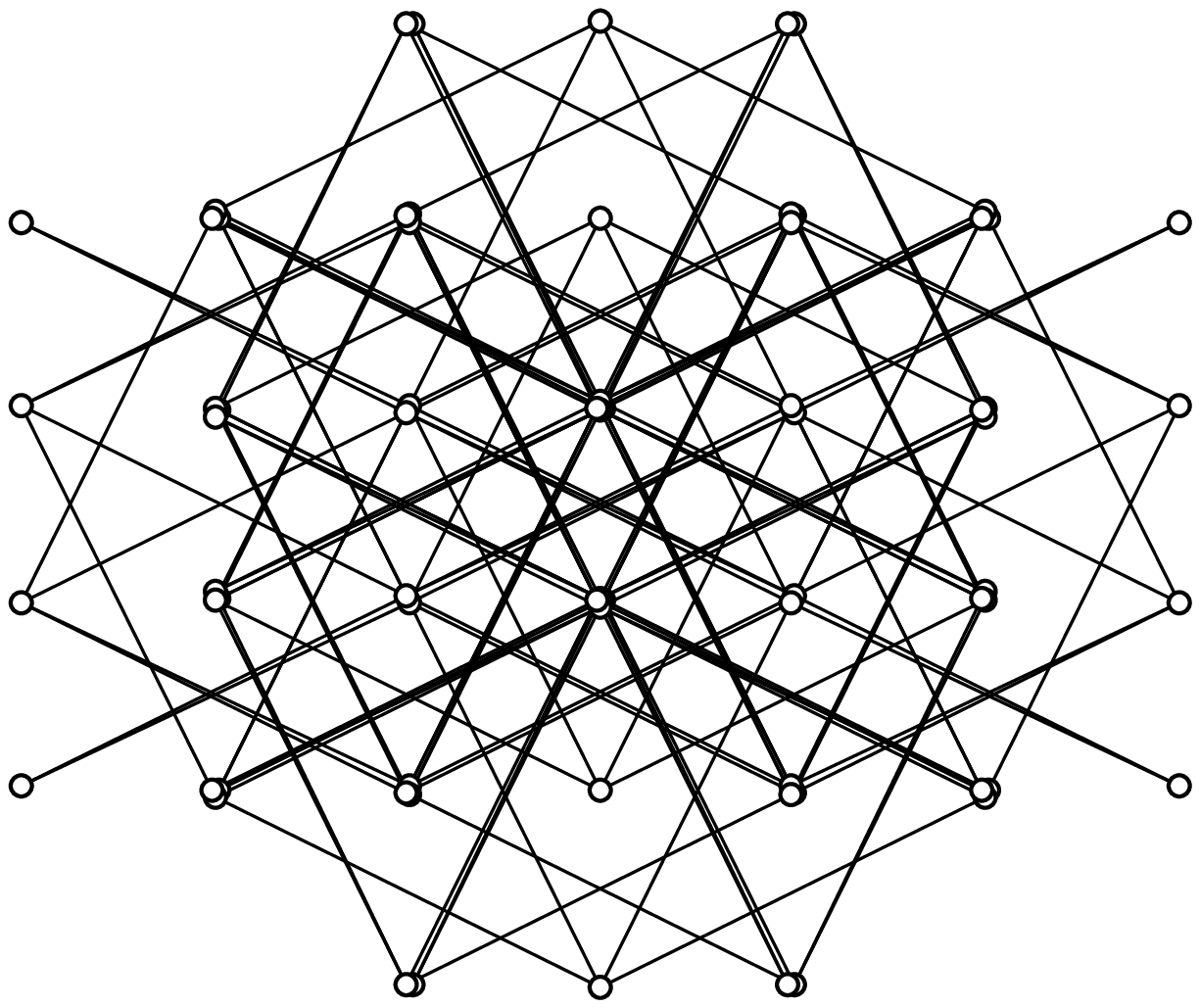} \tabularnewline \smallbreak \tabularnewline \includegraphics[scale=\stillFramesScale]{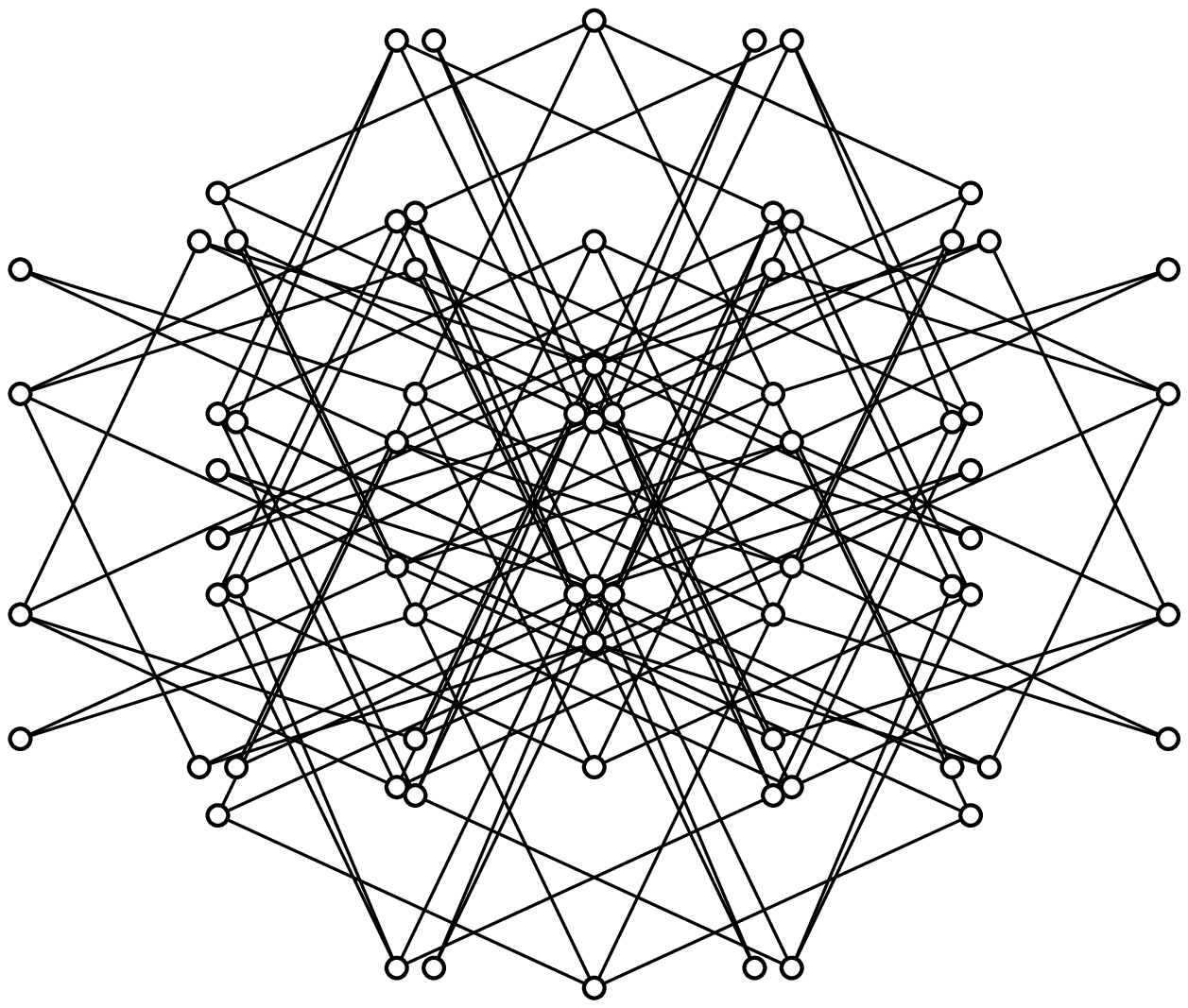} & \includegraphics[scale=\stillFramesScale]{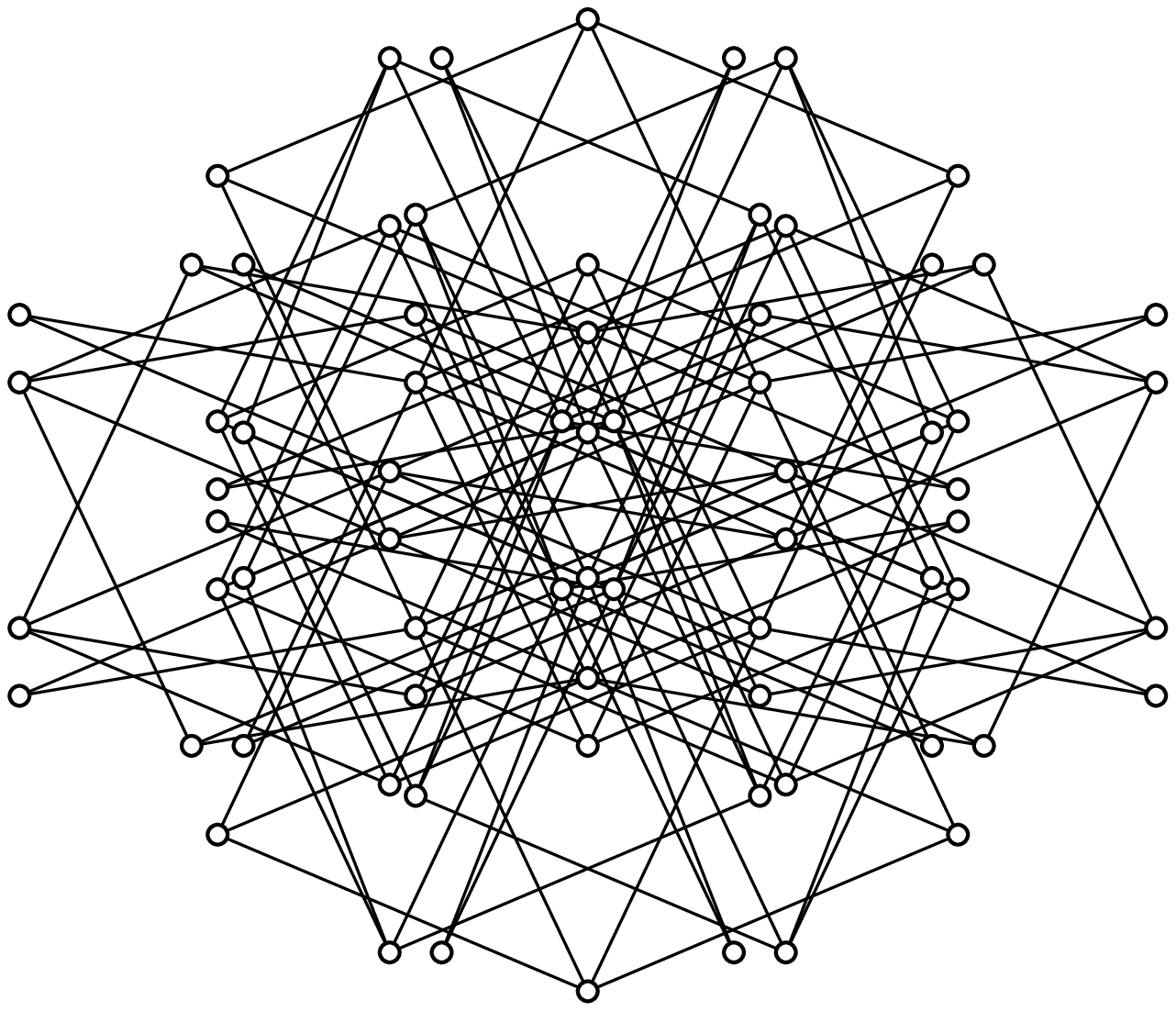} & \includegraphics[scale=\stillFramesScale]{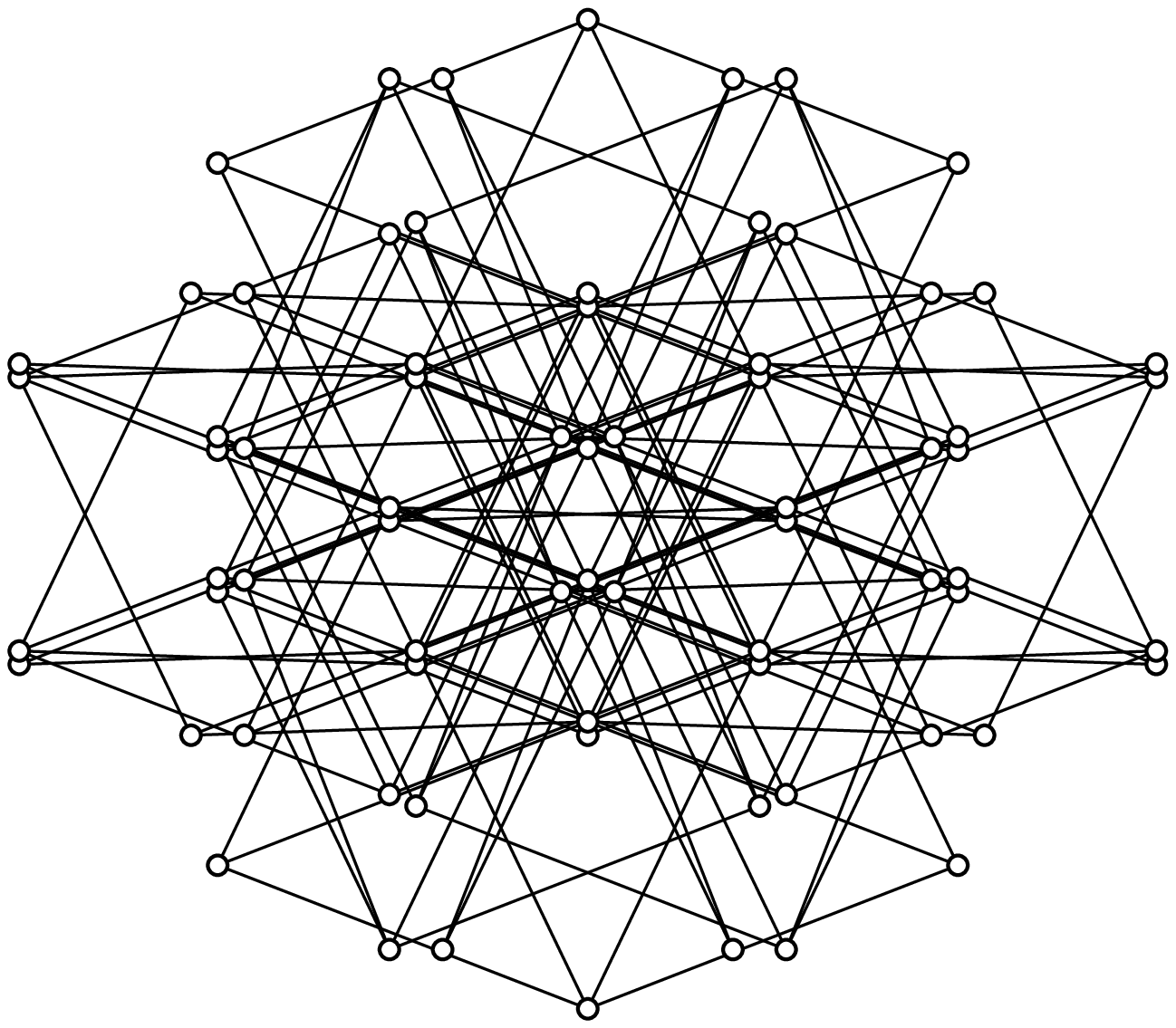} \end{tabular} \caption{} \label{zebra-still-frames} \end{figure}

For example, Figure \ref{zebra-still-frames} shows a sequence of ``still frames'' from the unique flexion of $\mathcal{F}_\Psi$ for the zebra. (Note that, ``off-screen'', it occasionally happens that multiple joints come to momentarily occupy the same point in the plane.)

Looking at Figure \ref{zebra-still-frames}, it is immediately striking how $\mathcal{F}_\Psi$ appears to retain its symmetries as it flexes.

The preservation of symmetry is not a general property of flexible frameworks. Counterexamples are readily available in the literature; a couple of them we review here for completeness.

Let $\mathcal{F}_\talloblong$ be any framework such that its joints $c_1$, $c_2$, $c_3$, and $c_4$ form a rectangle $c_1c_2c_3c_4$ which is not a square and its bars are the sides of that rectangle. Then $\mathcal{F}_\talloblong$ is axially symmetric with respect to two perpendicular axes and flexible with precisely one degree of freedom. However, flexing $\mathcal{F}_\talloblong$ destroys both of its axial symmetries.

On the other hand, let $\mathcal{F}_\boxslash$ be any framework such that its joints $c_1$, $c_2$, $c_3$, $c_4$, and $c_5$ form a square $c_1c_2c_3c_4$ together with its center $c_5$ and its bars are the sides of that square together with segments $c_1c_5$ and $c_3c_5$. Then $\mathcal{F}_\boxslash$ is centrally symmetric and flexible with precisely one degree of freedom. However, flexing $\mathcal{F}_\boxslash$ destroys its central symmetry.

We proceed to show that $\mathcal{F}_\Psi$ does indeed preserve all of its symmetries as it flexes. To this end, first we prove one lemma which sheds a bit more light on vectors $\mathbf{r}$, $\widehat{\mathbf{r}}$, $\mathbf{s}$, and $\widehat{\mathbf{s}}$.

\medbreak

\begin{lemma} $\dfrac{r_1}{s_1 + s_3} = \dfrac{r_2}{s_2 + s_4} = \dfrac{r_5 + r_7}{s_5} = \dfrac{r_6 + r_8}{s_6}$. \label{proportions} \end{lemma}

\medbreak

\begin{proof} Since $\dim \mathcal{H}^\text{Weight}_\Psi = 3$, vectors $\mathbf{r}$, $\widehat{\mathbf{r}}$, $\mathbf{s}$, and $\widehat{\mathbf{s}}$ are linearly dependent. Let \[\alpha_1\mathbf{r} + \alpha_2\widehat{\mathbf{r}} + \alpha_3\mathbf{s} + \alpha_4\widehat{\mathbf{s}} = \mathbf{0},\] where $\alpha_1$, $\alpha_2$, $\alpha_3$, and $\alpha_4$ are not all zero.

Then \begin{align*} (\alpha_1 - \alpha_2)(r_1, r_2, r_5 - r_7, r_6 - r_8) + (\alpha_3 - \alpha_4)(s_1 - s_3, s_2 - s_4, s_5, s_6) &= \mathbf{0} \text{ and }\\ (\alpha_1 + \alpha_2)(r_1, r_2, r_5 + r_7, r_6 + r_8) + (\alpha_3 + \alpha_4)(s_1 + s_3, s_2 + s_4, s_5, s_6) &= \mathbf{0}. \end{align*}

Since $\alpha_1$, $\alpha_2$, $\alpha_3$, and $\alpha_4$ are not all zero, at least one of pairs $\{\alpha_1 - \alpha_2, \alpha_3 - \alpha_4\}$ and $\{\alpha_1 + \alpha_2, \alpha_3 + \alpha_4\}$ is nonzero. Let $\varepsilon = \pm 1$ be such that pair $\{\alpha_1 + \varepsilon\alpha_2, \alpha_3 + \varepsilon\alpha_4\}$ is nonzero. Then \[\frac{r_1}{s_1 + \varepsilon s_3} = \frac{r_2}{s_2 + \varepsilon s_4} = \frac{r_5 + \varepsilon r_7}{s_5} = \frac{r_6 + \varepsilon r_8}{s_6}.\]

Observe that cases $\varepsilon = -1$ and $\varepsilon = 1$ cannot occur at the same time because of Lemma \ref{weight-psi-independence}. We are left to rule out the case of $\varepsilon = -1$ and to prove that $\varepsilon = 1$ is the case which does occur.

(Strictly speaking, by this point we already know enough about vectors $\mathbf{r}$, $\widehat{\mathbf{r}}$, $\mathbf{s}$, and $\widehat{\mathbf{s}}$ in order to prove Theorem \ref{flexibility-psi-symmetry}. But we are going to work out the exact value of $\varepsilon$ anyway.)

By Lemmas \ref{weight-psi-nonzero-i} and \ref{weight-psi-nonzero-ii}, \begin{align*} \frac{r_1}{s_1 + \varepsilon s_3} = \frac{r_2}{s_2 + \varepsilon s_4} &\Rightarrow \frac{r_1}{r_2} = \frac{s_1 + \varepsilon s_3}{s_2 + \varepsilon s_4}\\ &\Rightarrow \frac{r_1}{r_1 + r_2} = \frac{s_1 + \varepsilon s_3}{(s_1 + s_2) + \varepsilon(s_3 + s_4)}\\ &\Rightarrow \frac{r_1}{-2pq} = \frac{s_1 + \varepsilon s_3}{(p^2 + q^2) + \varepsilon(p^2 - q^2)}. \end{align*}

When $\varepsilon = -1$, this simplifies to $qr_1 + p(s_1 - s_3) = 0$. Since $p$ and $q$ are relatively prime, in this case we obtain that $p$ divides $r_1$.

Otherwise, when $\varepsilon = 1$, this simplifies to $pr_1 + q(s_1 + s_3) = 0$. Analogously, in this case we obtain that $q$ divides $r_1$.

This is how we are going to tell apart the two cases $\varepsilon = -1$ and $\varepsilon = 1$.

Recall cells $a_\texttt{I}$ and $a_\texttt{II}$ as well as subgrid $\Pi'$ of $\Psi$, as defined in the proof of Lemma \ref{weight-psi-nonzero-i}. Observe that the initial cells of all moves of $L$ in $C_\Psi$ of slope $\frac{q}{p}$ form exactly the lower left subgrid $\Pi'_\texttt{LL}$ of $\Pi'$ of height $p$ and width $2q$.

Since $\Pi'_\texttt{LL}$ intersects each net precisely once, by Lemma \ref{rhombic-class} we have that $r_1$ is the number of nets spanned by one connected component of graph $\mathcal{L}_{\Pi'} \setminus \frac{q}{p}$ and $r_2$ is the number of nets spanned by the other connected component of the same graph.

Denote these two connected components by $K'$ and $K''$, as in the proof of Lemma \ref{three-slopes-symmetry}. Suppose, without loss of generality, that $a_\texttt{I}$ belongs to $K'$, $r_1$ is the number of nets spanned by $K'$, $a_\texttt{II}$ belongs to $K''$, and $r_2$ is the number of nets spanned by $K''$.

By the proof of Lemma \ref{nets-cycle}, cell $a_\texttt{I} + r_1(p, -q)$ belongs to the same net as cell $a_\texttt{II}$. On the other hand, because $a_\texttt{II} = a_\texttt{I} + (0, 1)$, from this it follows that $pr_1 \equiv 0 \bmod q$ and $-qr_1 \equiv 1 \bmod p$.

The latter implies that $p$ does not divide $r_1$, and so it rules out the case of $\varepsilon = -1$. (Incidentally, by the former we also obtain that $q$ divides $r_1$, which is consistent with the case of $\varepsilon = 1$.) \end{proof}

\medbreak

We are ready to prove that $\mathcal{F}_\Psi$ preserves its symmetries as it flexes.

\medbreak

\begin{theorem} Suppose that $p \ge 2$. Then, as it flexes, the leaper framework of $L$ on $\Psi$ preserves all of its symmetries; namely, central symmetry as well as axial symmetry with respect to two perpendicular axes. \label{flexibility-psi-symmetry} \end{theorem}

\medbreak

\begin{proof} That $\mathcal{F}_\Psi$ preserves its central symmetry is clear, as each rhombic class of $L$ on $\Psi$ is symmetric with respect to the center of $\Psi$ by Lemma \ref{three-slopes-symmetry}. The preservation of axial symmetry, however, is somewhat trickier to establish.

By Corollary \ref{one-degree-psi}, it is enough to prove that $\mathcal{F}_\Psi$ admits some sufficiently symmetric flexion. This is exactly what we are going to do, in a manner very similar to our proof of Theorem \ref{flexibility-method}.

To factor out rotations, we decree that one axis of symmetry is vertical and the other one is horizontal. This has the added benefit of making all calculations much simpler.

In the notation of Sections \ref{inf-flex} and \ref{method}, set $x_3 = -x_1$, $y_3 = y_1$, $x_4 = -x_2$, $y_4 = y_2$, $x_7 = -x_5$, $y_7 = y_5$, $x_8 = -x_6$, and $y_8 = y_6$. Since each rhombic class of $L$ on $\Psi$ is symmetric with respect to the center of $\Psi$ by Lemma \ref{three-slopes-symmetry}, this is necessary and sufficient for the placement of $\mathcal{L}_\Psi$ whose summary is vector $(x_1, y_1, x_2, y_2, \ldots, x_8, y_8)$ to be axially symmetric with respect to one vertical and one horizontal axis.

We proceed to rewrite all bar-length constraints as well as our basis of all cyclic constraints in terms of the eight remaining variables $x_1$, $y_1$, $x_2$, $y_2$, $x_5$, $y_5$, $x_6$, and $y_6$.

The bar-length constraints which address $v_1$ and $v_3$ become one and the same, as do the ones which address $v_2$ and $v_4$, $v_5$ and $v_7$, and $v_6$ and $v_8$. Thus in this setting we only get four bar-length constraints.

The cyclic constraints associated with $\mathbf{r}_\texttt{X}$ and $\widehat{\mathbf{r}}_\texttt{X}$ become identical up to sign, and the ones associated with $\mathbf{r}_\texttt{Y}$ and $\widehat{\mathbf{r}}_\texttt{Y}$ become one and the same altogether. Thus in this setting we only get a spanning set of four cyclic constraints, namely \begin{align*} r_1x_1 + r_2x_2 + (r_5 - r_7)x_5 + (r_6 - r_8)x_6 &= 0,\\ r_1y_1 + r_2y_2 + (r_5 + r_7)y_5 + (r_6 + r_8)y_6 &= 0,\\ (s_1 - s_3)x_1 + (s_2 - s_4)x_2 + s_5x_5 + s_6x_6 &= 0, \text{ and }\\ (s_1 + s_3)y_1 + (s_2 + s_4)y_2 + s_5y_5 + s_6y_6 &= 0. \end{align*}

However, by Lemma \ref{proportions}, the second and fourth of these are scalar multiples of each other, and so we may omit either one of them.

In short, we obtain a system of four bar-length and three cyclic constraints such that vector $v = (x_1, y_1, x_2, y_2, -x_1, y_1, \ldots, -x_6, y_6)$ is a solution to this system if and only if it is the summary of some sufficiently symmetric placement of $\mathcal{L}_\Psi$.

From here on out, our plan for the proof of Theorem \ref{flexibility-psi-symmetry} will be exactly the same as our plan for the proof of Theorem \ref{flexibility-method}. That is, first we encode all constraints on $v$ into a single vector-valued function. Then we calculate the Jacobian matrix $M^\text{Sym}_\Psi$ of this vector-valued function at the point determined by the canonical values of all components of $v$. We prove that $M^\text{Sym}_\Psi$ has full rank. From this, we derive that it is possible to deform the canonical value of $v$ continuously so that it remains the summary of some sufficiently symmetric placement of $\mathcal{L}_\Psi$. Lastly, in conclusion, we obtain that $\mathcal{F}_\Psi$ admits some sufficiently symmetric flexion.

The only non-routine step of this plan is the calculation of the rank of $M^\text{Sym}_\Psi$, and so this is the only step which we are going to spell out in full.

Explicitly, \[M^\text{Sym}_\Psi = \begin{pmatrix} 2p & 2q & 0 & 0 & 0 & 0 & 0 & 0\\ 0 & 0 & 2p & 2q & 0 & 0 & 0 & 0\\ 0 & 0 & 0 & 0 & 2q & 2p & 0 & 0\\ 0 & 0 & 0 & 0 & 0 & 0 & 2q & 2p\\ r_1 & 0 & r_2 & 0 & r_5 - r_7 & 0 & r_6 - r_8 & 0\\ 0 & r_1 & 0 & r_2 & 0 & r_5 + r_7 & 0 & r_6 + r_8\\ s_1 - s_3 & 0 & s_2 - s_4 & 0 & s_5 & 0 & s_6 & 0 \end{pmatrix},\] and we aim to show that the rows of $M^\text{Sym}_\Psi$ are linearly independent.

Consider any linear combination of the rows of $M^\text{Sym}_\Psi$ whose coefficients are $\alpha_1$, $\alpha_2$, $\alpha_3$, $\alpha_4$, $\beta_1$, $\beta_2$, and $\beta_3$, respectively, and which equals the zero vector.

By columns $1$ and $3$ of $M^\text{Sym}_\Psi$, we have that \[2p(\alpha_1, \alpha_2) + \beta_1(r_1, r_2) + \beta_3[(s_1, s_2) - (s_3, s_4)] = \mathbf{0}.\]

Analogously, by columns $2$ and $4$ of $M^\text{Sym}_\Psi$, we have that \[2q(\alpha_1, \alpha_2) + \beta_2(r_1, r_2) = \mathbf{0}.\]

Eliminating $(\alpha_1, \alpha_2)$, we obtain \[(\beta_1 q - \beta_2 p)(r_1, r_2) + \beta_3 q[(s_1, s_2) - (s_3, s_4)] = \mathbf{0}.\]

By Lemma \ref{weight-psi-independence}, vectors $(s_1, s_2)$ and $(s_3, s_4)$ are linearly independent. On the other hand, by Lemma \ref{proportions} vectors $(r_1, r_2)$ and $(s_1, s_2) + (s_3, s_4)$ are linearly dependent. Since vector $(r_1, r_2)$ is nonzero, from this it follows that vectors $(r_1, r_2)$ and $(s_1, s_2) - (s_3, s_4)$ are linearly independent. Consequently, \begin{align*} \beta_1 q - \beta_2 p &= 0 \text{ and }\\ \beta_3 q &= 0, \end{align*} and so $\beta_3 = 0$.

Analogously, by columns $5$, $6$, $7$, and $8$ of $M^\text{Sym}_\Psi$, and in view of $\beta_3 = 0$, we have that \begin{align*} 2q(\alpha_3, \alpha_4) + \beta_1[(r_5, r_6) - (r_7, r_8)] &= \mathbf{0} \text{ and }\\ 2p(\alpha_3, \alpha_4) + \beta_2[(r_5, r_6) + (r_7, r_8)] &= \mathbf{0}. \end{align*}

Eliminating $(\alpha_3, \alpha_4)$, we obtain \[(\beta_1 p - \beta_2 q)(r_5, r_6) - (\beta_1 p + \beta_2 q)(r_7, r_8) = \mathbf{0}.\]

By Lemma \ref{weight-psi-independence}, it follows that \begin{align*} \beta_1 p - \beta_2 q &= 0 \text{ and }\\ \beta_1 p + \beta_2 q &= 0, \end{align*} and so $\beta_1 = \beta_2 = 0$.

Thus all three of $\beta_1$, $\beta_2$, and $\beta_3$ must equal zero. On the other hand, since the first four rows of $M^\text{Sym}_\Psi$ are linearly independent, all four of $\alpha_1$, $\alpha_2$, $\alpha_3$, and $\alpha_4$ must equal zero as well. Therefore, the rows of $M^\text{Sym}_\Psi$ are linearly independent, as required. \end{proof}

\medbreak

Our discussion of the case $p \ge 2$ of Theorem \ref{flexibility-psi} is complete.

\section{Flexibility II} \label{flex-ii}

We go on to resolve the case $p = 1$ of Theorem \ref{flexibility-psi}. Thus, from this point on throughout the rest of this section, suppose that $p = 1$. Consequently, $q$ is any even positive integer.

When we attempt to apply the method of Section \ref{method} to $\Psi$ directly, we meet with significant difficulties as the number of rhombic classes of $L$ on $\Psi$ becomes very large. We get around this obstacle by embedding $\Psi$ into one larger set of cells and then applying our method to that larger set instead. (As we said we would in Section \ref{method}.)

We form our larger set of cells by deleting one cell from $\Phi$. The middle cell of any side of $\Phi$ would do; for concreteness, we choose the top. Formally, suppose that $\Phi$ is the grid $[1; 2q + 1] \times [1; 2q + 1]$, as in Section \ref{rigid}. Then we define $\Phi\text{*}$ to be the set of cells $\Phi \setminus \{(q + 1, 2q + 1)\}$.

\medbreak

\begin{theorem} Suppose that $p = 1$. Then the leaper framework of $L$ on $\Phi\text{*}$ is flexible. \label{flexibility-phi-star} \end{theorem}

\medbreak

\begin{figure}[p] \centering \includegraphics[scale=\figuresScale]{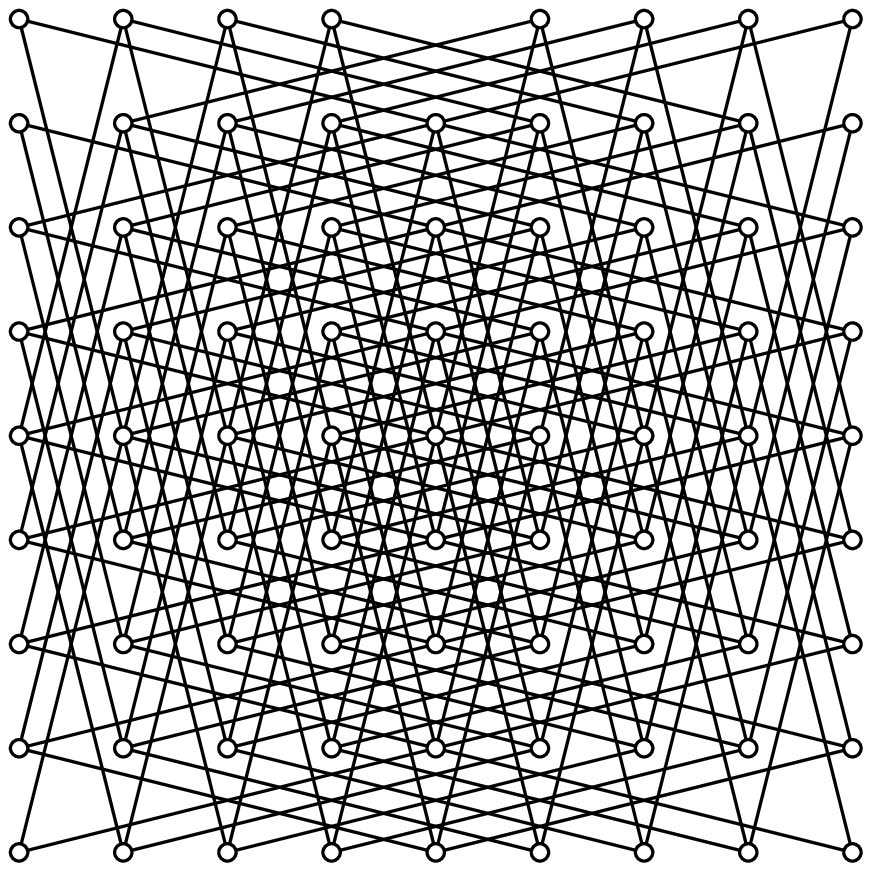} \caption{} \label{giraffe-phi-star} \end{figure}

\begin{figure}[p] \centering \includegraphics[scale=\figuresScale]{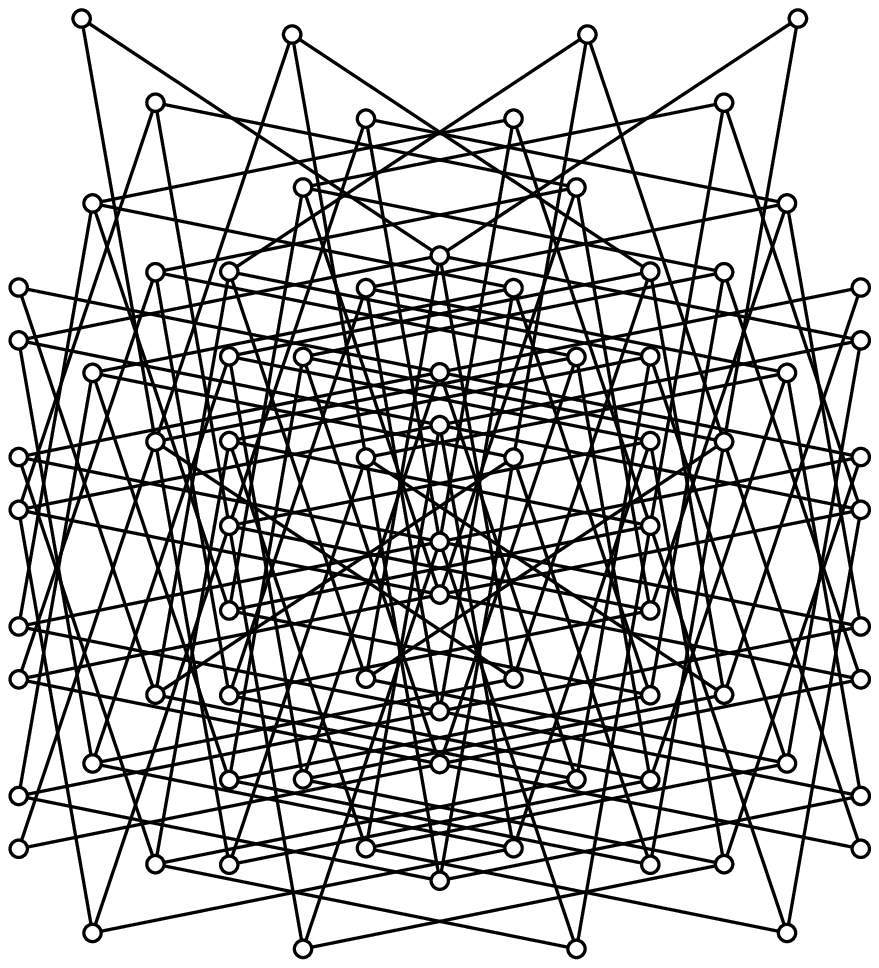} \caption{} \label{giraffe-phi-star-flex} \end{figure}

Intuitively, Theorem \ref{flexibility-phi-star} states that, when $p = 1$, the leaper framework of $L$ on $\Phi$ is balanced precariously on the precipice of flexibility. For example, Figure \ref{giraffe-phi-star} shows $\mathcal{F}_{\Phi\text{*}}$ for the giraffe, and Figure \ref{giraffe-phi-star-flex} shows the same framework in the process of flexing.

Clearly, the case $p = 1$ of Theorem \ref{flexibility-psi} follows by Theorem \ref{flexibility-phi-star} along the lines of the proof of Lemma \ref{monotonicity}.

Note that $\Phi\text{*}$ still has a vertical axis of symmetry, just as $\Psi$ does. However, it lacks all other kinds of symmetry which $\Psi$ possesses.

We proceed to retrace the steps of our proof of the case $p \ge 2$ of Theorem \ref{flexibility-psi}, all the while keeping careful track of what has to change and what may safely stay the same. Luckily, since in the present setting everything is determined by just one single parameter, many things which used to be complicated in Section \ref{flex-i} will now become much simpler.

We begin, as usual, with the rhombic classes of $L$ on $\Phi\text{*}$. We employ the generalisation of Lemma \ref{rhombic-class} to arbitrary sets of cells which we introduced in Section \ref{method}.

We consider slope $\frac{q}{1}$ first. Let $\Gamma'$ be the set of cells \[\Phi\text{*}\llbracket(1, q)\rrbracket = ([1; 2q] \times [1; q + 1]) \setminus \{(q, q + 1)\}.\]

By the generalised Lemma \ref{rhombic-class}, in order to describe the rhombic classes of $L$ on $\Phi\text{*}$ for slope $\frac{q}{1}$, it is enough to describe the connected components of graph $\mathcal{L}_{\Gamma'} \setminus \frac{q}{1}$.

Recall from Section \ref{slope} that, when $p = 1$, the forbidden-slope leaper graph $\mathcal{L}_\Gamma \setminus \frac{q}{p}$ is a Hamiltonian path of $L$ on $\Gamma$. Thus it is far from surprising that the deletion of one cell breaks this graph down into two connected components.

Explicitly, the cell set of one of these connected components is \begin{align*} K_\texttt{I} &= [\{(i, 2j - 1) \mid 1 \le i \le q \text{ and } 1 \le j \le \tfrac{1}{2}q + 1\} \setminus \{(q, q + 1)\}] \cup \mbox{ }\\ &\phantom{= \mbox{ }} \{(i, 2j) \mid q + 1 \le i \le 2q \text{ and } 1 \le j \le \tfrac{1}{2}q\}, \end{align*} and the cell set of the other one is \begin{align*} K_\texttt{II} &= \{(i, 2j) \mid 1 \le i \le q \text{ and } 1 \le j \le \tfrac{1}{2}q\} \cup \mbox{ }\\ &\phantom{= \mbox{ }} \{(i, 2j - 1) \mid q + 1 \le i \le 2q \text{ and } 1 \le j \le \tfrac{1}{2}q + 1\}. \end{align*}

\begin{figure}[p] \centering \includegraphics[scale=\figuresScale]{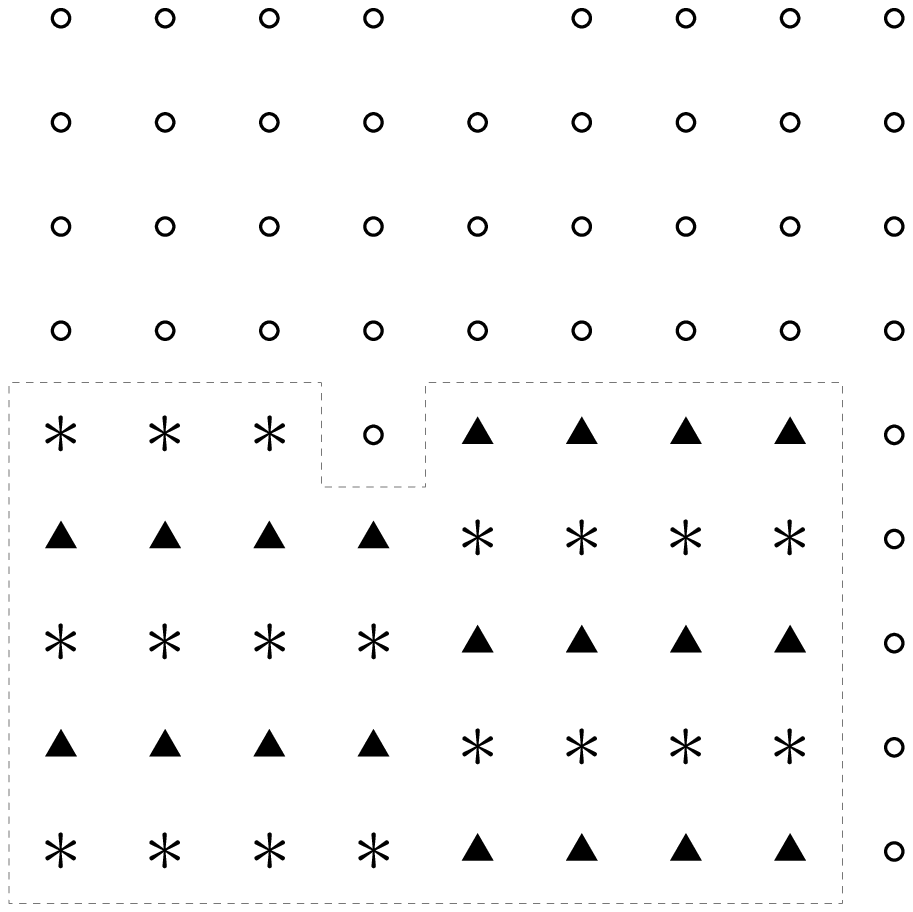} \caption{} \label{giraffe-phi-star-components-1} \end{figure}

For example, Figure \ref{giraffe-phi-star-components-1} shows $\Gamma'$, $K_\texttt{I}$, and $K_\texttt{II}$ for the giraffe. Observe that $|K_\texttt{I}| = q^2 + q - 1$ and $|K_\texttt{II}| = q^2 + q$.

We denote the two rhombic classes of $L$ on $\Phi\text{*}$ which correspond to $K_\texttt{I}$ and $K_\texttt{II}$ by $R_1$ and $R_2$, respectively. Moreover, as in Section \ref{inf-flex}, we denote the reflections of $R_1$ and $R_2$ across the axis of symmetry of $\Phi\text{*}$ by $R_3$ and $R_4$, respectively. Then $R_1$ and $R_2$ are the two rhombic classes of $L$ on $\Phi\text{*}$ for slope $\frac{q}{1}$, and $R_3$ and $R_4$ are the two rhombic classes of $L$ on $\Phi\text{*}$ for slope $-\frac{q}{1}$.

We continue with slope $\frac{1}{q}$, and we handle it analogously to slope $\frac{q}{1}$.

Let $\Gamma''$ be the set of cells \[\Phi\text{*}\llbracket(q, 1)\rrbracket = ([1; q + 1] \times [1; 2q]) \setminus \{(1, 2q)\}.\]

By the generalised Lemma \ref{rhombic-class}, in order to describe the rhombic classes of $L$ on $\Phi\text{*}$ for slope $\frac{1}{q}$, it is enough to describe the connected components of graph $\mathcal{L}_{\Gamma''} \setminus \frac{1}{q}$.

Once again, there are two of them. The cell set of one is \begin{align*} K_\texttt{V} &= \{(2i, q) \mid 1 \le i \le \tfrac{1}{2}q\} \cup \mbox{ }\\ &\phantom{= \mbox{ }} \{(2i - 1, 2q) \mid 2 \le i \le \tfrac{1}{2}q + 1\}, \end{align*} and the cell set of the other one is most conveniently described as \[K_\texttt{VI} = \Gamma'' \setminus K_\texttt{V}. \]

\begin{figure}[p] \centering \includegraphics[scale=\figuresScale]{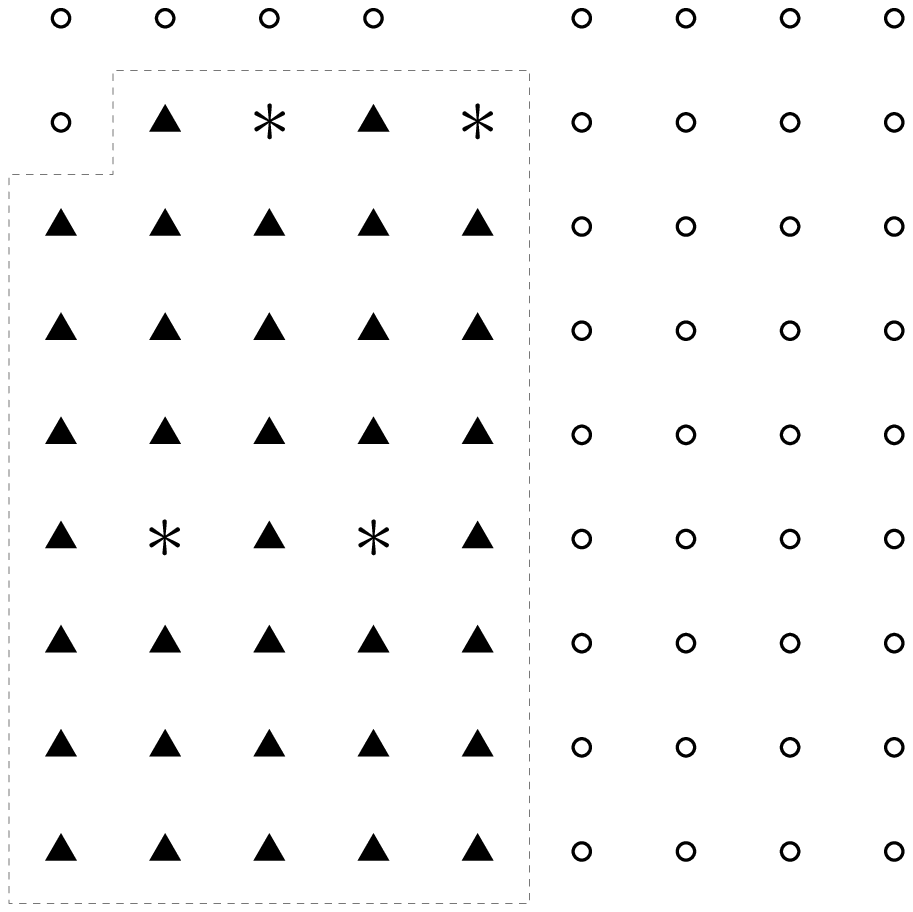} \caption{} \label{giraffe-phi-star-components-2} \end{figure}

For example, Figure \ref{giraffe-phi-star-components-2} shows $\Gamma''$, $K_\texttt{V}$, and $K_\texttt{VI}$ for the giraffe. Observe that $|K_\texttt{V}| = q$ and $|K_\texttt{VI}| = 2q^2 + q - 1$.

We denote the two rhombic classes of $L$ on $\Phi\text{*}$ which correspond to $K_\texttt{V}$ and $K_\texttt{VI}$ by $R_5$ and $R_6$, respectively. Once again as in Section \ref{inf-flex}, we also denote the reflections of $R_5$ and $R_6$ across the axis of symmetry of $\Phi\text{*}$ by $R_7$ and $R_8$, respectively. Then $R_5$ and $R_6$ are the two rhombic classes of $L$ on $\Phi\text{*}$ for slope $\frac{1}{q}$, and $R_7$ and $R_8$ are the two rhombic classes of $L$ on $\Phi\text{*}$ for slope $-\frac{1}{q}$.

Thus, echoing Lemma \ref{rhombic-psi}, there are a total of eight rhombic classes of $L$ on $\Phi\text{*}$, two per each slope.

We prove that $\mathcal{H}^\text{Weight}_{\Phi\text{*}}$ is of dimension three, and so also that the cyclic constraint space $\mathcal{H}^\text{Cyc}_{\Phi\text{*}}$ of $L$ on $\Phi\text{*}$ is of dimension six, mostly as in Section \ref{inf-flex}. One fine point is worth mentioning, however, as follows.

The subgrids $\Psi_\texttt{LR}$ and $\Psi_\texttt{LL}$ of $\Psi$ played an important role in our analysis of the dimensions of $\mathcal{H}^\text{Weight}_\Psi$ and $\mathcal{H}^\text{Cyc}_\Psi$. Their analogues on $\Phi\text{*}$ are the two sets of cells $\Phi\text{*}\llbracket(-1, q)\rrbracket$ and $\Phi\text{*}\llbracket(q, 1)\rrbracket$. Neither one of them is a rectangular grid, and so Lemmas \ref{single-slope} and \ref{nets-connectedness} do not apply to them anymore.

Thus, on $\Phi\text{*}$, we have to check by hand that the restriction of each net to $\Phi\text{*}\llbracket(-1, q)\rrbracket$ is nonempty and connected, and also that the number of connected components of graph $\mathcal{L}_{\Phi\text{*}\llbracket(q, 1)\rrbracket} \restriction -\frac{1}{q}$ works out as it should. Fortunately, both of these claims are fairly straightforward to verify.

We assign directions $d_1$, $d_2$, \ldots, $d_8$ to the rhombic classes of $L$ on $\Phi\text{*}$ just as in Section \ref{flex-i}. Then the analogue of Lemma \ref{weight-psi-reflection} continues to hold for $\Phi\text{*}$.

We construct the oriented cycle $C_{\Phi\text{*}}$ of $L$ on $\Phi\text{*}$ in the exact same way as the oriented cycle $C_\Psi$ of $L$ on $\Psi$ in Section \ref{flex-i}.

This time around, we simply calculate the weight $\mathbf{r}$ of $C_{\Phi\text{*}}$ explicitly. (Contrast this with what we did for the weight of $C_\Psi$ in Section \ref{flex-i}.) Direct counting yields \[\mathbf{r} = (-q,\; -q,\; 0,\; 0,\; \tfrac{1}{2}q,\; q^2 - \tfrac{1}{2}q + 1,\; \tfrac{1}{2}q,\; q^2 - \tfrac{1}{2}q - 1).\]

We define $C'_{\Phi\text{*}}$ to be the $90^\circ$ rotation of $C_{\Phi\text{*}}$ about point $(\frac{1}{2}q + 1, \frac{1}{2}q + 1)$. (We cannot define $C'_{\Phi\text{*}}$ to be the reflection of $C_{\Phi\text{*}}$ in the line $x = y$ anymore because this reflection contains cell $(q + 1, 2q + 1)$, which lies outside of $\Phi\text{*}$.)

Once again, we simply calculate the weight $\mathbf{s}$ of $C'_{\Phi\text{*}}$ explicitly. It works out to \[\mathbf{s} = (-\tfrac{1}{2}q^2 + 1,\; -\tfrac{1}{2}q^2,\; \tfrac{1}{2}q^2,\; \tfrac{1}{2}q^2 + 1,\; 0,\; 0,\; -1,\; -2q + 1).\]

The natural analogue of Lemma \ref{weight-psi-independence} as well as the literal statement of Lemma \ref{weight-psi-basis} for $\mathbf{r}$ and $\mathbf{s}$ both follow routinely, and the rest of the argument does not require any adjustments whatsoever.

This completes our proof of Theorem \ref{flexibility-phi-star}, and so also our proof of Theorem \ref{flexibility-psi}.

Analogues of Corollary \ref{one-degree-psi} and Theorem \ref{flexibility-psi-symmetry} hold on $\Phi\text{*}$ as well, as follows.

\medbreak

\begin{corollary} Suppose that $p = 1$. Then the leaper framework of $L$ on $\Phi\text{*}$ is flexible with precisely one degree of freedom. \label{one-degree-phi-star} \end{corollary}

\medbreak

\begin{theorem} Suppose that $p = 1$. Then the leaper framework of $L$ on $\Phi\text{*}$ preserves its axial symmetry as it flexes. \label{flexibility-phi-star-symmetry} \end{theorem}

\medbreak

The proofs are fully analogous.

With this, our discussion of Theorem \ref{flexibility-psi} is complete.

\section{Further Work} \label{further}

By way of a conclusion, we survey a number of open problems.

To begin with: On what grids is the leaper framework of $L$ rigid?

This is perhaps the central open problem in the study of leaper frameworks. Theorem \ref{square} solves it for square grids. Arbitrary rectangular grids, however, are much tougher to classify completely. Still, we expect that the method which we outlined in Section \ref{method} will be of help in this more general setting as well.

Lemma \ref{monotonicity} allows us to sketch out the rough shape of the answer well in advance, as follows.

Let $A$ be any grid larger than the $1 \times 1$ grid. We say that $A$ is \emph{minimally rigid} for $L$ if the leaper framework of $L$ is rigid on $A$ and flexible on all grids smaller than $A$, with the trivial exception of the $1 \times 1$ grid. For example, $\Phi$ is minimally rigid for $L$ by Theorems \ref{rigidity-phi} and \ref{flexibility-psi}. (Is $\Theta$, when $p \ge 2$?)

Let $\mathfrak{R}_L$ be the set of the sizes of all grids minimally rigid for $L$. For convenience, we are going to talk about $\mathfrak{R}_L$ as if it is a set of grids rather than a set of grid sizes. That is, we are going to say simply ``grid $A$ is in $\mathfrak{R}_L$'' instead of ``the size of grid $A$ is in $\mathfrak{R}_L$''.

By Lemma \ref{monotonicity}, the leaper framework of $L$ on $A$ is rigid if and only if $A$ is larger than or equal to some grid in $\mathfrak{R}_L$. In other words, $\mathfrak{R}_L$ immediately yields the complete classification that we seek.

Note that Knuth's complete classification of all grids on which the leaper graph of $L$ is connected is of the same overall shape. In essence, Knuth proved that the set of the sizes of all grids \emph{minimally connected} for $L$ is given by $\mathfrak{C}_L = \{(p + q) \times 2q, 2q \times (p + q)\}$.

It is easy to see that $\mathfrak{R}_L$ is finite. We mentioned already that $\Phi$ is in $\mathfrak{R}_L$; on the other hand, by Theorem \ref{rigidity-theta} and Corollary \ref{incompleteness}, there is more to $\mathfrak{R}_L$ than that at least when $p \ge 2$.

One more piece of low-hanging fruit is well worth the picking.

For all positive integers $n$, the leaper framework of $L$ on the grid of size $(p + q) \times n$ is flexible. The proof is largely similar to the proof of Lemma \ref{monotonicity}. We proceed by induction on $n$. The base case is clear. For the induction step, each joint in every newly added column is joined by at most two nonparallel bars to the rest of the framework. Therefore, the addition of one more column always preserves flexibility.

Consequently, each side of every grid in $\mathfrak{R}_L$ must be greater than or equal to $p + q + 1$.

We do not have much to say about $\mathfrak{R}_L$ beyond that. For example, we do not know if the lower bound of $p + q + 1$ is attained for all leapers. We do not even know if the size of $\mathfrak{R}_L$ is bounded from above over all free leapers $L$ (as is the case with connectedness) or not.

One step towards the long-term goal of working out $\mathfrak{R}_L$ could be to determine if, on rectangular grids, the conditions of Theorem \ref{flexibility-method} are necessary as well as sufficient. Or, more ambitiously, if global rigidity with distinctness, infinitesimal rigidity, and continuous rigidity coincide for all leaper frameworks on rectangular grids, as they do for the leaper framework of $L$ on $\Phi$.

Most of our analysis applies just as well to a similar question suggested by Section \ref{slope}: On what grids are the forbidden-slope leaper graphs of $L$ connected?

For concreteness, we only consider the case when the forbidden slope is $\frac{q}{p}$. All other cases can be obtained from it by reflection and rotation.

An analogue of Lemma \ref{monotonicity} holds here as well. (As we hinted at, briefly, in our discussion of Corollary \ref{rigidity-phi-stronger}.) Note that, in this instance, it matters quite a lot that we compare grids by means of the relation ``smaller than'' instead of ``fits inside''. The former forbids rotation and the latter allows it. By contrast, the same distinction does not matter at all in the setting of flexibility and rigidity.

Let, then, $A$ be \emph{minimally forbidden-slope connected} for $L$ if the restriction of $\mathcal{L} \setminus \frac{q}{p}$ is connected on $A$ and disconnected on all grids smaller than $A$, with the trivial exception of the $1 \times 1$ grid. For example, $\Gamma$ is minimally forbidden-slope connected for $L$. Indeed, the case $p = 1$ follows because $\mathcal{L}_\Gamma \setminus \frac{q}{p}$ is then a Hamiltonian path of $L$ on $\Gamma$, and the case $p \ge 2$ follows by Lemma \ref{three-slopes-components}. When $p \ge 2$, $\Lambda'$ and the $90^\circ$ rotation of $\Lambda''$ are minimally forbidden-slope connected for $L$ as well. The proof is not difficult.

Moreover, let $\mathfrak{S}_L$ be the set of the sizes of all grids minimally forbidden-slope connected for $L$. Then, by the appropriate analogue of Lemma \ref{monotonicity}, $\mathcal{L}_A \setminus \frac{q}{p}$ is connected if and only if $A$ is larger than or equal to some grid in $\mathfrak{S}_L$. Thus, once again, $\mathfrak{S}_L$ immediately yields the complete classification that we seek.

Of course, working out $\mathfrak{S}_L$ is an interesting question in its own right. Quite apart from that, however, it will likely be of great help in determining $\mathfrak{R}_L$ as well.

Theorem \ref{flexibility-psi-symmetry} raises one natural question: Does every flexible leaper framework on a rectangular grid admit a flexion which preserves all of its symmetries?

So does Theorem \ref{flexibility-phi-star}: Suppose that $A$ is minimally rigid for $L$. What is the least number of joints (or bars) that we need to remove from the leaper framework of $L$ on $A$ in order to make it flexible? For example, what is this number for the leaper framework of $L$ on $\Phi$?

There are also the higher-dimensional leapers to consider.

Let $\mathcal{S}$ be any multiset of nonnegative integers and let $N$ be the size of $\mathcal{S}$. The $\mathcal{S}$-leaper $L^\mathcal{S}$ lives on the $N$-dimensional integer lattice $\mathbb{Z}^N$, and two points $(x'_1, x'_2, \ldots, x'_N)$ and $(x''_1, x''_2, \ldots, x''_N)$ of this lattice are joined by an edge of $L^\mathcal{S}$ if and only if the multisets $\mathcal{S}$ and $\{|x'_1 - x''_1|, |x'_2 - x''_2|, \ldots, |x'_N - x''_N|\}$ coincide.

The definitions of a leaper graph and a leaper framework extend to higher dimensions in an obvious way. We denote the leaper graph of $L^\mathcal{S}$ on $\mathbb{Z}^N$ by $\mathcal{L}^\mathcal{S}$, the leaper graph of $L^\mathcal{S}$ on the subset $C$ of $\mathbb{Z}^N$ by $\mathcal{L}^\mathcal{S}_C$, and the leaper framework of $L^\mathcal{S}$ on $C$ by $\mathcal{F}^\mathcal{S}_C$. Just as with two-dimensional leapers, we say that $L^\mathcal{S}$ is free if $\mathcal{L}^\mathcal{S}$ is connected, and, as far as flexibility and rigidity are concerned, we only consider the leaper frameworks of free higher-dimensional leapers.

Let $\mathcal{S}_+$ be the multiset of all nonzero elements of $\mathcal{S}$ and let $N_+$ be the size of $\mathcal{S}_+$. It makes sense to think of leapers $L^{\mathcal{S}_+}$ and $L^\mathcal{S}$ as two different forms of the same leaper. Intuitively, when we let $L^{\mathcal{S}_+}$ loose into a space higher-dimensional than its natural habitat, it acts in the exact same way as the native $L^\mathcal{S}$. We say that $L^\mathcal{S}$ is \emph{proper} if all elements of $\mathcal{S}$ are nonzero, and \emph{improper} otherwise. If multisets $\mathcal{S'_+}$ and $\mathcal{S}''_+$ coincide, then we say that leapers $L^{\mathcal{S}'}$ and $L^{\mathcal{S}''}$ are \emph{equivalent}.

In \cite{SW} and \cite{W}, Solymosi and White study the leaper frameworks of higher-dimensional leapers equivalent to the knight. They prove that if $N \ge 3$ and $L^\mathcal{S}$ is equivalent to the knight, then the leaper framework of $L^\mathcal{S}$ is infinitesimally rigid on the $N$-dimensional grid of size $4 \times 4 \times \cdots \times 4$, and thus also rigid on this as well as all larger grids.

Their result generalises without any trouble as follows. Suppose that $n \ge 2$, $N_+ \ge 2$, and the leaper framework of $L^\mathcal{S_+}$ is infinitesimally rigid on the $N_+$-dimensional grid of size $n \times n \times \cdots \times n$. Then the leaper framework of $L^\mathcal{S}$ is infinitesimally rigid on the $N$-dimensional grid of size $n \times n \times \cdots \times n$, and thus also rigid on this as well as all larger grids.

In particular, the above claim and Corollary \ref{rigidity-phi-stronger}, taken together, tell us that if $L^\mathcal{S}$ is equivalent to a free proper two-dimensional leaper, then its leaper framework is rigid on all sufficiently large $N$-dimensional grids.

More generally, suppose that $L^\mathcal{S}$ is free. Is its leaper framework rigid on all sufficiently large $N$-dimensional grids?

It would be very interesting to see an analogue of Theorem \ref{square} for higher-dimensional leapers. And, of course, a complete classification would be even better. A lot of the ideas and techniques of the present work generalise to higher dimensions in a natural way; however, in all likelihood, many novel insights will be necessary as well.

This concludes our discussion of the flexibility and rigidity of leaper frameworks.

\end{document}